\documentclass{amsart} 

\usepackage[all]{xy}
\usepackage{amssymb} 
\usepackage{amsmath} 
\usepackage{amscd} 
\usepackage{latexsym} 
\usepackage{amsthm} 
\usepackage{verbatim} 
\usepackage{graphics} 
\usepackage{graphicx}
\usepackage[small,nohug,heads=littlevee]{diagrams}
\diagramstyle[labelstyle=\scriptstyle]

\let\cal\mathcal 
\def\Ascr{{\cal A}} 
\def\Bscr{{\cal B}} 
\def\Cscr{{\cal C}} 
\def\Dscr{{\cal D}} 
\def\Escr{{\cal E}}

\def\Iscr{{\cal I}} 
\def\Jscr{{\cal J}} 
 
\def\Lscr{{\cal L}} 
\def\Mscr{{\cal M}} 
\def\Nscr{{\cal N}} 
\def\Oscr{{\cal O}} 
\def\Pscr{{\cal P}} 
\def\Qscr{{\cal Q}} 
\def\Rscr{{\cal R}} 
\def\Sscr{{\cal S}} 
 
\def\Uscr{{\cal U}} 
\def\Vscr{{\cal V}} 
\def\Wscr{{\cal W}} 
\def\Xscr{{\cal X}} 
\def\Yscr{{\cal Y}} 
\def\Zscr{{\cal Z}} 

\let\blb\mathbb

\def\CC{{\blb C}}

\def\II{{\blb I}}
 
\def \NN{{\blb N}} 
\def \PP{{\blb P}}

\def\XX{{\blb X}}
\def\YY{{\blb Y}}
\def \ZZ{{\blb Z}} 
\def\1{{\blb{1}}}

\def\cd{\operatorname {cd}}

\def\Cl{\operatorname{Cl}} 
 
\def\coh{\mathop{\text{\upshape{coh}}}} 
\def\coker{\operatorname {coker}}

\def\D{\operatorname{D}}

\def\div{\operatorname{div}}

\def\Ext{\operatorname {Ext}} 
 
\def\Filt{\operatorname {Filt}} 
\def\filt{\operatorname {filt}}

\def\gkdim{\operatorname {GKdim}} 
\def\Gl{\operatorname {Gl}} 
 
\def\gldim{\operatorname {gldim}} 
\def\gr{\operatorname{gr}}

\def\grmod{\operatorname {grmod}} 
\def\GrMod{\operatorname {GrMod}} 
\def\Hilb{\operatorname {Hilb}} 
 
\def\Hom{\operatorname {Hom}} 
 
\def\id{{\operatorname {id}}} 
\def\Id{\operatorname{id}} 
\def\im{\operatorname {im}}

\def\ker{\operatorname {ker}}

\def\length{\mathop{\text{length}}} 
\def\Lie{\mathop{\text{Lie}}}

\def\LotimesU{\overset{{\bf L}}{\otimes}_{U}}

\def\LotimesGamma{\overset{{\bf L}}{\otimes}_{\Gamma}}
\def\lr{\longrightarrow} 

\def\mod{\operatorname{mod}} 
\def\Mod{\operatorname{Mod}}

\def\P2q{\operatorname {{\blb P}^{2}_{q}}}

\def\pd{\operatorname {pd}}

\def\Pic{\operatorname {Pic}} 
\def\PP{\operatorname {\blb P}} 
\def\pr{\mathop{\text{pr}}\nolimits}

\def\Qcoh{\operatorname {Qcoh}}

\def\quot{/\!\!/} 
\def\r{\rightarrow} 
 
\def\rank{\operatorname {rank}} 
\def\red{\operatorname {red}}
\def\Rees{\operatorname {Rees}}
\def\Rep{\operatorname {Rep}} 
 
\def\Res{\operatorname{Res}} 
\def\RHom{\operatorname {RHom}} 
\def\rk{\operatorname {rk}} 
\def\sh{\operatorname{sh}}

\def\Tails{\operatorname {Tails}}  
\def\tails{\operatorname {tails}}  
 
\def\Tor{\operatorname {Tor}}

\def\tr{\operatorname {Tr}}

\def\un{\operatorname {--}}

\DeclareMathOperator{\Aut}{Aut}

\DeclareMathOperator{\Ind}{Ind}

\DeclareMathOperator{\Proj}{Proj}

\DeclareMathOperator{\tors}{tors} 
\DeclareMathOperator{\Tors}{Tors}

\newtheorem*{warning}{Warning} 

\newtheorem{theoremd}{Theorem}  

\newtheorem{lemma}{Lemma}[section] 
\newtheorem{proposition}[lemma]{Proposition} 
\newtheorem{theorem}[lemma]{Theorem} 
\newtheorem{corollary}[lemma]{Corollary} 
 
\newtheorem{convention}[lemma]{Convention}

\theoremstyle{definition} 
\newtheorem{example}[lemma]{Example} 
\newtheorem{definition}[lemma]{Definition}

\newtheorem*{Acknowledgements}{Acknowledgements}

{ 
 
\newtheorem{step}{Step} 
 
\newtheorem{case}{Case}

}

\theoremstyle{remark} 
\newtheorem{remark}[lemma]{Remark} 
 
\newtheorem*{notation}{Notation}

\newdimen\uboxsep \uboxsep=1ex 
\def\uboxn#1{\vtop to 0pt{\hrule height 0pt depth 0pt\vskip\uboxsep 
\hbox to 0pt{\hss #1\hss}\vss}} 

\def\uboxs#1{\vbox to 0pt{\vss\hbox to 0pt{\hss #1\hss} 
\vskip\uboxsep\hrule height 0pt depth 0pt}}

\numberwithin{equation}{section}

\newcommand{\ra}{\rightarrow}


\author{Koen De Naeghel and Nicolas Marconnet} 

\address{Koen De Naeghel \\Departement WNI\\Hasselt University\\ Agoralaan gebouw D \\ B-3590 
Diepenbeek (Belgium).}
\email[K. De Naeghel]{koen.denaeghel@uhasselt.be}
\address{Nicolas Marconnet \\ 
Departement Wiskunde en Informatica\\Universiteit Antwerpen\\Middelheimlaan 1\\ B-2020 Antwerp (Belgium).}
\email[N. Marconnet]{
nicolas.marconnet@ua.ac.be} 

\date{January 5, 2006} 

\title[Ideals of cubic algebras and an invariant ring of the Weyl algebra]{Ideals of cubic algebras and an invariant ring of the Weyl algebra} 

\begin{document}

\begin{abstract}
We classify reflexive graded right ideals, up to isomorphism and shift, of generic cubic three dimensional Artin-Schelter regular algebras. We also determine the possible Hilbert functions of these ideals. These results are obtained by using similar methods as for quadratic Artin-Schelter algebras \cite{DV1,DV2}. In particular our results apply to the enveloping algebra of the Heisenberg-Lie algebra from which we deduce a classification of right ideals of the invariant ring $A_{1}^{\langle \varphi \rangle}$ of the first Weyl algebra $A_1 = k\langle x,y \rangle/(xy - yx - 1)$ under the automorphism $\varphi(x) = -x$, $\varphi(y) = -y$.
\end{abstract}

\maketitle

\setcounter{tocdepth}{1} 
\tableofcontents

\section{Introduction} 

Let $k$ be an algebraically closed field of characteristic zero and consider the first Weyl algebra $A_1 = k\langle x,y\rangle/(yx-xy-1)$. It is well-known that $A_1$ is a simple ring, however it is a natural question to describe one-sided ideals of $A_1$. In 1994 Cannings and Holland \cite{CH1} classified right $A_1$-ideals by means of the adelic Grassmannian. A few years later Wilson \cite{wilson} found a relation between the adelic Grassmannian and Calogero-Moser spaces, to obtain the following result (as it was formulated by Berest and Wilson) 
\begin{theorem} \cite{BW0} \label{ref-1.1-0} 
Put $G = \Aut(A_1)$ and let $R(A_1)$ be the set of isomorphism classes of 
right $A_1$-ideals. Then the orbits of the natural $G$-action on 
$R(A_1)$ are indexed by the set of non-negative integers $\NN$, and the orbit corresponding 
to $n\in \NN$ is in natural bijection with the $n$'th Calogero-Moser 
space 
\begin{equation*} 
C_n = \{\XX,\YY \in M_n(\CC)\mid \rk(\YY\XX-\XX\YY-\II)\leq1\} /\Gl_n(\CC) 
\end{equation*} 
where  $\Gl_n(\CC)$ acts  by simultaneous conjugation $(gXg^{-1},gYg^{-1})$. In particular $C_n$ is a connected smooth affine variety of dimension $2n$.
\end{theorem} 
In \cite{BW} Berest and Wilson gave a new proof of Theorem \ref{ref-1.1-0} using noncommutative algebraic geometry \cite{AZ,VW}, by picking up an idea of Le Bruyn \cite{LeBruyn}. Let us briefly recall how this is done. The idea of Le Bruyn was to relate right $A_1$-ideals with graded right ideals of the homogenized Weyl algebra $H = \CC \langle x,y,z\rangle/(zx-xz,zy-yz,yx-xy-z^2)$. Such a relation comes from the fact that $A_1$ coincides with $H[z^{-1}]_0$, the degree zero part of the localisation of $H$ at the powers of the central element $z$. Describing $R(A_1)$ then becomes equivalent to describing certain objects on a noncommutative projective scheme $\P2q = \Proj H$ in the sense of Artin and Zhang \cite{AZ} i.e. $\Proj H$ is the quotient of the abelian category of finitely generated graded right $H$-modules by the Serre subcategory of finite dimensional modules. Objects on $\P2q$ may be used to define moduli spaces, just as in the commutative case.

We next observe there are many more noncommutative projective schemes like $\Proj H$. The homogenized Weyl algebra $H$ belongs to a class of algebras which has turned out to be very fruitful: the so-called Artin-Schelter regular algebras $A$ of global dimension $d$ \cite{AS}. See Section \ref{Preliminaries} for preliminary definitions. For $d \leq 3$ a complete classification is known \cite{ATV1,ATV2,Steph1,Steph2}. They are all noetherian domains of Gelfand-Kirillov dimension $d \leq 3$ and may be considered as noncommutative analogues of the polynomial rings. In case $d = 3$ and $A$ is generated in degree one it turns out \cite{AS} there are two possibilities for such an algebra $A$. Either there are three generators $x,y,z$ and three quadratic relations (we say $A$ is quadratic) or two generators $x,y$ and two cubic relations ($A$ is cubic). If $A$ is quadratic then $A$ is Koszul and has Hilbert series $h_A(t) = 1 + 3t + 6t^2 + 10t^3 + \dots = 1/(1-t)^3$, and we may think of $\P2q = \Proj(A)$ as a noncommutative projective plane. In case $A$ is cubic then $h_A(t) = 1 + 2t + 4t^2 + 6t^3 + \dots = 1/(1-t)^2(1-t^2)$, we then write $(\PP^1 \times \PP^1)_q = \Proj(A)$ which we think of as a noncommutative quadric surface. The generic class of quadratic and cubic Artin-Schelter regular algebras are usually called type A-algebras \cite{AS}, in which case the relations are respectively given by
\begin{equation*} 
\left\{ 
\begin{array}{l} 
ayz + bzy + cx^{2} =0 \\ 
azx + bxz + cy^{2}=0\\ 
axy + byx + cz^{2}=0 
\end{array} \right. 
\text{ and }
\left\{ 
\begin{array}{l} 
ay^{2}x + byxy + axy^{2} + cx^{3} = 0  \\ 
ax^{2}y + bxyx + ayx^{2} + cy^{3} = 0
\end{array} \right. 
\end{equation*} 
where $a,b,c \in k$ are generic scalars. As shown in \cite{ATV1} a three dimensional Artin-Schelter regular algebra $A$ generated in degree one is completely determined by a triple $(E,\sigma,j)$ where either
\begin{itemize}
\item
$j: E \cong \PP^2$ if $A$ is quadratic, resp. $j: E \cong \PP^1 \times \PP^1$ if $A$ is cubic; or
\item
$j:E \hookrightarrow \PP^2$ is an embedding of a divisor $E$ of degree three if $A$ is quadratic, resp.
$j:E \hookrightarrow \PP^1 \times \PP^1$ where $E$ is a divisor of bidegree $(2,2)$ if $A$ is cubic 
\end{itemize}
and $\sigma \in \Aut(E)$. In the first case we say $A$ is linear, otherwise $A$ is called elliptic. If $A$ is of type A and the divisor $E$ is a smooth elliptic curve (this is the generic case) then we say $A$ is of generic type A. In that case $\sigma$ is a translation on $E$. Quadratic three dimensional Artin-Schelter regular algebras of generic type A are also called three dimensional Sklyanin algebras.  

In \cite{DV1} Van den Bergh and the first author showed how to extend the ideas in \cite{BW,LeBruyn} to obtain a classification of reflexive graded right ideals of generic quadratic Artin-Schelter regular algebras. This has been extended in the PhD thesis of the first author.
\begin{theorem} \cite{DV1,D} \label{ref-1.2-2} 
Assume $k$ is uncountable. Let $A$ be an elliptic quadratic Artin-Schelter regular algebra for which $\sigma$ has infinite order. Let $R(A)$ be the set of reflexive graded right $A$-ideals, considered up to isomorphism and shift of grading. There exist smooth locally closed varieties $D_n$ of dimension $2n$ such that $R(A)$ is naturally in bijection with $\coprod_{n \in \NN} D_n$. If in addition $A$ is of generic type A i.e. $A$ is a three dimensional Sklyanin algebra then the varieties $D_n$ are affine. 
\end{theorem} 
In case $A$ is a three dimensional Sklyanin algebra Theorem \ref{ref-1.2-2} holds without the hypothesis $k$ is uncountable.

A result similar to Theorem \ref{ref-1.2-2} was proved by Nevins and Stafford \cite{NS} for all quadratic three dimensional Artin-Schelter regular algebras $A$. In addition they showed $D_n$ is an open subset in a projective variety $\Hilb_{n}(\P2q)$ of dimension $2n$, parameterizing graded right $A$-ideals of projective dimension one up to isomorphism and shift of grading. Thus $\Hilb_{n}(\P2q)$ is the analog of the classical Hilbert scheme of points on 
$\PP^{2}$. Furthermore $\Hilb_n(\PP^2_q)$ is connected, proved by Nevins and Stafford \cite{NS} for almost all $A$ using deformation-theoretic methods and relying on the commutative case $A = k[x,y,z]$. In \cite{DV2} Van den Bergh and the first author determined the Hilbert series of objects in $\Hilb_{n}(\P2q)$ and deduced an intrinsic proof for the connectedness of $\Hilb_n(\PP^2_q)$ for all quadratic Artin-Schelter regular algebras.

In this paper we apply the methods used in \cite{DV1,DV2} to obtain similar results for cubic Artin-Schelter regular algebras. Most of these results may also be found in the PhD thesis of the first author \cite{D}. 

Let $A$ be a cubic Artin-Schelter regular algebra and let $R(A)$ denote the set of reflexive graded right $A$-ideals, considered up to isomorphism and shift of grading. Define $N = \{(n_e,n_o) \in \NN^2 \mid n_e - (n_e - n_o)^2 \geq 0 \}$. The main result in this paper is the following analogue of Theorem \ref{ref-1.2-2}.
\begin{theoremd} \label{theorem1}
Assume $k$ is uncountable. Let $A$ be an elliptic cubic Artin-Schelter regular algebra for which $\sigma$ has infinite order. Then for $(n_e,n_o) \in N$ there exists a smooth locally closed variety $D_{(n_e,n_o)}$ of dimension $2(n_e - (n_e - n_o)^2)$ such that $R(A)$ is in natural bijection with $\coprod_{(n_e,n_o) \in N} D_{(n_e,n_o)}$. If in addition $A$ is of generic type A then the varieties $D_{(n_e,n_o)}$ are affine. 
\end{theoremd}
Theorem \ref{theorem1} will follow from Theorem \ref{ref-5.5.5-65cubic} and Theorem \ref{ref-5.5.4-59cubic} below, where a description of the appearing varieties $D_{(n_e,n_o)}$ is given. In particular $D_{(0,0)}$ is a point and $D_{(1,1)}$ is the complement of $C = E_{\red}$ under a natural embedding in $\PP^1 \times \PP^1$, see Corollary \ref{smallinvcubic}. In fact $D_{(n_e,n_o)}$ is a point whenever $n_e = (n_e - n_o)^2$. 

A crucial part of the proof of Theorem \ref{theorem1} consists in showing that the spaces $D_{(n_e,n_o)}$ are actually nonempty for $(n_e,n_o) \in N$.
In contrast to quadratic Artin-Schelter regular algebras \cite{DV1,NS} this is not entirely straightforward. We will prove this by characterizing the appearing Hilbert series for objects in $R(A)$. In a very similar way as in \cite{DV2} for quadratic Artin-Schelter regular algebras, we show in Section \ref{Hilbert series of torsion free rank one modules} that the Hilbert series of graded right $A$-ideals of projective dimension one are characterised by so-called {\em Castelnuovo polynomials} \cite{Davis} $s(t)=\sum_{i=0}^{n}s_it^i \in \ZZ[t]$ which are by definition of the form 
\begin{equation*} 
s_0=1,s_1=2,\ldots,s_{\sigma-1}=\sigma \mbox{ and } s_{\sigma-1}\geq s_{\sigma}\geq s_{\sigma+1}\geq \cdots \geq 0
\end{equation*} 
for some integer $\sigma \geq 0$. We refer to $\sum_{i}s_{2i}$ as the  {\em even weight} of $s(t)$ and $\sum_{i}s_{2i+1}$ as the {\em odd weight} of $s(t)$.
\begin{example} 
$s(t) = 1 + 2t + 3t^{2} + 4t^{3} + 5t^{4} + 5t^{5} + 3t^{6} + 2t^{7} + t^{8} + t^{9} + t^{10} + t^{11}$ is a Castelnuovo polynomial of even weight $14$ and odd weight $15$. The corresponding Castelnuovo diagram is (where the even columns are black) \\ 

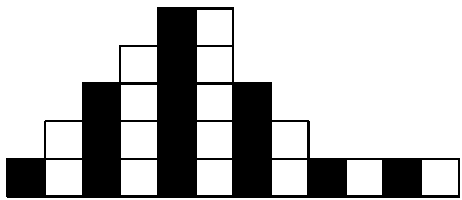
\unitlength 1mm
\begin{picture}(90.00,25.00)(0,0)

\linethickness{0.15mm}
\put(30.00,0.00){\line(0,1){5.00}}

\linethickness{0.15mm}
\put(30.00,5.00){\line(1,0){5.00}}

\linethickness{0.15mm}
\put(35.00,5.00){\line(0,1){5.00}}

\linethickness{0.15mm}
\put(35.00,10.00){\line(1,0){5.00}}

\linethickness{0.15mm}
\put(40.00,10.00){\line(0,1){5.00}}

\linethickness{0.15mm}
\put(40.00,15.00){\line(1,0){5.00}}

\linethickness{0.15mm}
\put(45.00,15.00){\line(0,1){5.00}}

\linethickness{0.15mm}
\put(45.00,20.00){\line(1,0){5.00}}

\linethickness{0.15mm}
\put(50.00,20.00){\line(0,1){5.00}}

\linethickness{0.15mm}
\put(50.00,25.00){\line(1,0){5.00}}

\linethickness{0.15mm}
\put(55.00,0.00){\line(0,1){25.00}}

\linethickness{0.15mm}
\put(55.00,15.00){\line(1,0){5.00}}

\linethickness{0.15mm}
\put(60.00,10.00){\line(0,1){5.00}}

\linethickness{0.15mm}
\put(60.00,10.00){\line(1,0){5.00}}

\linethickness{0.15mm}
\put(65.00,0.00){\line(0,1){10.00}}

\linethickness{0.15mm}
\put(65.00,5.00){\line(1,0){20.00}}

\linethickness{0.15mm}
\put(85.00,0.00){\line(0,1){5.00}}

\linethickness{0.15mm}
\put(30.00,0.00){\line(1,0){55.00}}

\linethickness{0.15mm}
\put(50.00,0.00){\line(0,1){20.00}}

\linethickness{0.15mm}
\put(45.00,0.00){\line(0,1){15.00}}

\linethickness{0.15mm}
\put(40.00,0.00){\line(0,1){10.00}}

\linethickness{0.15mm}
\put(35.00,0.00){\line(0,1){5.00}}

\linethickness{0.15mm}
\put(35.00,5.00){\line(1,0){30.00}}

\linethickness{0.15mm}
\put(40.00,10.00){\line(1,0){20.00}}

\linethickness{0.15mm}
\put(45.00,15.00){\line(1,0){10.00}}

\linethickness{0.15mm}
\put(50.00,20.00){\line(1,0){5.00}}

\linethickness{0.15mm}
\put(60.00,0.00){\line(0,1){10.00}}

\linethickness{0.15mm}
\put(70.00,0.00){\line(0,1){5.00}}

\linethickness{0.15mm}
\put(75.00,0.00){\line(0,1){5.00}}

\linethickness{0.15mm}
\put(80.00,0.00){\line(0,1){5.00}}

\linethickness{0.15mm}
\put(55.00,25.00){\line(1,0){5.00}}

\linethickness{0.15mm}
\put(60.00,15.00){\line(0,1){10.00}}

\linethickness{0.15mm}
\put(55.00,20.00){\line(1,0){5.00}}

\linethickness{0.15mm}
\put(60.00,15.00){\line(1,0){5.00}}

\linethickness{0.15mm}
\put(65.00,10.00){\line(0,1){5.00}}

\linethickness{0.15mm}
\put(65.00,10.00){\line(1,0){5.00}}

\linethickness{0.15mm}
\put(70.00,5.00){\line(0,1){5.00}}

\linethickness{0.15mm}
\put(85.00,5.00){\line(1,0){5.00}}

\linethickness{0.15mm}
\put(90.00,0.00){\line(0,1){5.00}}

\linethickness{0.15mm}
\put(85.00,0.00){\line(1,0){5.00}}

\linethickness{0.15mm}
\put(30.00,0.00){\rule{5.00\unitlength}{5.00\unitlength}}

\linethickness{0.15mm}
\put(40.00,0.00){\rule{5.00\unitlength}{15.00\unitlength}}

\linethickness{0.15mm}
\put(50.00,0.00){\rule{5.00\unitlength}{25.00\unitlength}}

\linethickness{0.15mm}
\put(60.00,0.00){\rule{5.00\unitlength}{15.00\unitlength}}

\linethickness{0.15mm}
\put(70.00,5.00){\rule{5.00\unitlength}{0.00\unitlength}}

\linethickness{0.15mm}
\put(70.00,0.00){\rule{5.00\unitlength}{5.00\unitlength}}

\linethickness{0.15mm}
\put(80.00,0.00){\rule{5.00\unitlength}{5.00\unitlength}}

\end{picture}
\end{example}
Denote $X = \Proj(A) = (\PP^1 \times \PP^1)_q$. Write $\Hilb_{(n_e,n_o)}(X)$ for the groupoid of the torsion free graded right $A$-modules $I$ of projective dimension one for which
\begin{eqnarray*}
\dim_k A_m - \dim_k I_m =
\left\{ 
\begin{array}{ll}
n_e & \text{ for $m$ even} \\
n_o & \text{ for $m$ odd}
\end{array}
\right.
\text{ for } m \gg 0
\end{eqnarray*} 
(in particular $I$ has rank one, see \S\ref{Normalized line bundles}). Any graded right $A$-ideal $I$ of projective dimension one admits an unique shift of grading $I(d)$ for which $I(d) \in \Hilb_{(n_e,n_o)}(X)$. Writing $R_{(n_e,n_o)}(A)$ for the full subcategory of $\Hilb_{(n_e,n_o)}(X)$ consisting of reflexive objects we have $R(A) = \coprod R_{(n_e,n_o)}(A)$, and in the setting of Theorem \ref{theorem1} the isoclasses of $R_{(n_e,n_o)}(A)$ are in natural bijection with the points of the variety $D_{(n_e,n_o)}$. In Section \ref{Hilbert series of torsion free rank one modules} below we prove the following analog of \cite[Theorem A]{DV2}. 
\begin{theoremd} \label{theorem2}
Let $A$ be a cubic Artin-Schelter regular algebra. There is a bijective correspondence between Castelnuovo polynomials $s(t)$ of even weight $n_e$ and odd weight $n_o$ and Hilbert series $h_{I}(t)$ of objects $I$ in $\Hilb_{(n_e,n_o)}(X)$, given by 
\begin{equation*} 
h_{I}(t) = \frac{1}{(1-t)^{2}(1-t^2)} - \frac{s(t)}{1-t^2} 
\end{equation*}
Moreover if $A$ is elliptic for which $\sigma$ has infinite order this correspondence restricts to Hilbert series $h_I(t)$ of objects $I$ in $R_{(n_e,n_o)}(A)$.
\end{theoremd}
By shifting the rows in a Castelnuovo diagram in such a way they are left aligned one sees that the number of Castelnuovo diagrams of even weight $n_e$ and odd weight $n_o$ is equal to the number of partitions $\lambda$ of $n_e + n_o$ with distinct parts, with the additional property that by putting a chessboard pattern on the Ferrers diagram of  $\lambda$ the number of black squares is equal to $n_e$ and the number of white squares is equal to $n_o$. Anthony Henderson pointed out to us this number is given by the number of partitions of $n_e - (n_e - n_o)^2$. Alternatively, see \cite{DM}. In particular the varieties $D_{(n_e,n_o)}$ in Theorem \ref{theorem1} are nonempty whenever $(n_e,n_o) \in N$.
\begin{remark}
In Appendix \ref{ref-B-67} we have included the list of Castelnuovo polynomials $s(t)$ of even weight $n_e \leq 3$ and odd weight $n_o \leq 3$, as well as some associated data.
\end{remark} 
As there exists no commutative cubic Artin-Schelter regular algebra $A$ it seems difficult to compare our results with the commutative situation. However if $A$ is a linear cubic Artin-Schelter regular algebra then $\Proj(A)$ is equivalent with the category of coherent sheaves on the quadric surface $\PP^1 \times \PP^1$. In Section \ref{Ideals of linear cubic Artin-Schelter regular algebras} we discuss how the (classical) Hilbert scheme of points $\Hilb(\PP^1 \times \PP^1)$ parameterizes the objects in $\coprod_{(n_e,n_o) \in N} \Hilb_{(n_e,n_o)}(X)$ with the groupoid $\Hilb_{(n_e,n_o)}(X)$ as defined above.
\begin{remark}
For cubic Artin-Schelter regular algebras $A$ we expect a similar treatment as in \cite{NS} to show $\Hilb_{(n_e,n_o)}(X)$ is a smooth projective variety of dimension $2(n_e - (n_e - n_o)^2)$. The authors are convinced that using the same methods as in the proof of \cite[Theorem B]{DV2} will lead to a proof that $\Hilb_{(n_e,n_o)}(X)$ is connected, hence also $D_{(n_e,n_o)}$ (for elliptic $A$ for which $\sigma$ has infinite order). We hope to come back on this in a subsequent paper.
\end{remark}
As an application, consider the enveloping algebra of the Heisenberg-Lie algebra 
\begin{equation*}
\begin{aligned}
H_{c} = k \langle x,y,z \rangle / (yz - zy,xz - zx,xy - yx - z) = 
k \langle x,y \rangle / ([y,[y,x]],[x,[x,y]]) 
\end{aligned}
\end{equation*}
where $[a,b] = ab-ba$. The graded algebra $H_c$ is a cubic Artin-Schelter regular algebra. Consider the localisation $\Lambda = H_c[z^{-1}]$ of $H_c$ at the powers of the central element $z=xy-yx$ and its subalgebra $\Lambda_0$ of elements of degree zero. It was shown in \cite{ATV2} that $\Lambda_0 = A_1^{\langle \varphi \rangle}$, the algebra of invariants of the first Weyl algebra $A_1 = k\langle x,y\rangle/(xy-yx-1)$ under the automorphism $\varphi$ defined by $\varphi(x)=-x$, $\varphi(y)=-y$. In Section \ref{Ideals of an invariantring of the first Weyl algebra} we deduce a classification of right ideals of $A_1^{\langle \varphi \rangle}$.
\begin{theoremd} \label{theorem3}
The set $R(A_1^{\langle \varphi \rangle})$ of isomorphism classes of 
right $A_1^{\langle \varphi \rangle}$-ideals is in natural bijection with the points of $\coprod_{(n_e,n_o) \in N}D_{(n_e,n_o)}$ where 
\begin{multline*}
D_{(n_e,n_o)} = \{ (\XX,\YY,\XX',\YY') \in M_{n_{e} \times n_{o}}(k)^2 \times M_{n_{o} \times n_{e}}(k)^2 \mid 
\YY'\XX - \XX'\YY = \II \text{ and }   \\ 
\rank (\YY\XX'-\XX\YY'-\II)\leq 1 \}/\Gl_{n_e}(k)\times\Gl_{n_o}(k)  
\end{multline*}
is a smooth affine variety $D_{(n_e,n_o)}$ of dimension $2(n_e - (n_e - n_o)^2)$. 
\end{theoremd}
Note $\Gl_{n_e}(k)\times \Gl_{n_o}(k)$ acts by conjugation $(g\XX h^{-1},g \YY h^{-1},h \XX' g^{-1},h \YY' g^{-1})$.
Comparing with Theorem \ref{ref-1.1-0} it would be interesting to see if the orbits of $R(A_1^{\langle \varphi \rangle})$ under the automorphism group $\Aut(A_1^{\langle \varphi \rangle})$ are in bijection to $D_{(n_e,n_o)}$.

Finally, in Section \ref{filtrationscubic} we describe the elements of $R(A)$ by means of filtrations.
\begin{theoremd} \label{theorem4}
Assume $k$ is uncountable. Let $A$ be an elliptic cubic Artin-Schelter regular algebra and assume $\sigma$ has infinite order. Let $I \in R(A)$. Then there exists an $m \in \NN$ together with a filtration of reflexive graded right $A$-modules of rank one 
\[
I_0 \supset I_1 \supset \cdots \supset I_m = I
\]
with the property that up to finite length modules the 
quotients are shifted conic modules i.e. modules of the form $A/bA$ where $b \in A$ has degree two. Moreover $I_0$ 
admits a minimal resolution of the form
\begin{equation} \label{commall}
0 \r A(-c-1)^{c} \r A(-c)^{c+1} \r I_0 \r 0
\end{equation}
for some integer $c \geq 0$, and $I_0$ is up to isomorphism uniquely determined by $c$.
\end{theoremd} 
If $A$ is linear it follows from Proposition \ref{commfreecubic} below that every reflexive graded right ideal of $A$ admits a resolution of the form \eqref{commall} (up to shift of grading). Hence Theorem \ref{theorem4} is trivially true for linear cubic Artin-Schelter regular algebras.
\begin{remark}
In case $A$ is of generic type A or $A = H_c$ is the enveloping algebra of the Heisenberg-Lie algebra, Theorems \ref{theorem1} and 4 are proved without the hypothesis $k$ is uncountable. 
\end{remark}
Most results in this paper are obtained {\em mutatis mutandis} as for quadratic algebras in \cite{DV1,DV2} and to some extend \cite{BW,LeBruyn,NS}. However at some points the situation for cubic algebras is more complicated.
\begin{Acknowledgements}
Both authors were funded by the European Research Training Network Liegrits, which made it possible for the authors to work together. We are therefore grateful to J. Alev, R. Berger, T. Levasseur, B. Keller and F. Van Oystaeyen. We thank M. Van den Bergh for showing us his preprint \cite{ncquadrics}. We would like to thank R. Berger and A. Henderson  for their discussions and comments. 
\end{Acknowledgements}

\section{Preliminaries} \label{Preliminaries}

In this section we give some definitions and results used in this paper. These are collected from \cite{AS,ATV1,ATV2,AV,AZ,DV1,MS,NVO1,VdBSt}. Alternatively, see \cite[Chapter 1]{D}.

Throughout we work over an algebraically closed field $k$ of characteristic zero. We will use
\begin{convention}\label{XyUvW} 
Whenever $\mathrm{XyUvw}(\cdots)$ denotes an abelian category then \linebreak $\mathrm{xyuvw}(\cdots)$ denotes the full subcategory of $\mathrm{XyUvw}(\cdots)$ consisting of noetherian objects. 
\end{convention}

\subsection{Graded algebras and modules}

Let $A=\oplus_{i\in \NN} A_i$ be a connected graded $k$-algebra. We write $\GrMod(A)$ for the category of graded right $A$-modules with morphisms the $A$-module homomorphisms of degree zero. Since $\GrMod(A)$ is an abelian  category with enough injective objects we may define the functors $\Ext_{A}^{n}(M,-)$ on $\GrMod(A)$ as the right derived functors of $\Hom_{A}(M,-)$. It is convenient to write $\underline{\Hom}_{A}(M,-)=\bigoplus_{d \in \ZZ}^{}{\Hom_{A}(M,N(d))}$ and
\[ 
\underline{\Ext}_{A}^{n}(M,N) := \bigoplus_{d \in \ZZ}^{}{\Ext_{A}^{n}(M,N(d))}.
\] 
Let $M$ be a graded right $A$-module. We use the notation $M_{\geq n} = \oplus_{d \geq n}^{}{M_{d}}$ and $M_{\leq n} = \oplus_{d \leq n}^{}{M_{d}}$ for all $n \in \ZZ$. 
For any integer $n$, define $M(n)$ as the graded $A$-module equal to $M$ with its original $A$ action, but which is graded by $M(n)_{i} = M_{n + i}$. We refer to the functor $M \mapsto M(n)$ as the $n$-th {\it shift functor}.

The $k$-dual of a $k$-vector space $V$ is denoted by $V' = \Hom_{k}(V,k)$. The graded dual of a graded right $A$-module $M$ is $M^{\ast} = \underline{\Hom}_{A}(M,A)$ and $M$ is said to be {\em reflexive} if $M^{\ast \ast} = M$. We also write $(-)'$ for the functor on graded $k$-vector spaces which sends $M$ to its Matlis dual $M' = \underline{\Hom}_{k}(M,k) =
\oplus_{n}M'_{-n}$.

\subsection{Tails}

Let $A$ be a noetherian connected graded $k$-algebra. We denote by $\tau$ the functor sending a graded right $A$-module to the sum of all its finite dimensional submodules. Denote by $\Tors(A)$ the full subcategory of $\GrMod(A)$ consisting of all modules $M$ such that $\tau M = M$ and write $\Tails(A)$ for the quotient category $\GrMod(A)/\Tors(A)$. Objects in $\Tails(A)$ will be denoted by script letters like $\Mscr$. The shift functor on $\GrMod(A)$ induces an automorphism $\sh : \Mscr \mapsto \Mscr(1)$ on $\Tails(A)$. 

We write $\pi : \GrMod(A) \ra \Tails(A)$ for the (exact) quotient functor and $\Oscr = \pi A$. The functor $\pi$ has a right adjoint $\omega$. It is well-known that $\pi\circ\omega=\id$ and $\omega = \underline{\Hom}_{\Tails(A)}(\Oscr,-)$.  

When there is no possible confusion we write $\Hom$ instead of $\Hom_{A}$ and \linebreak
$\Hom_{\Tails(A)}$. The context will make clear in which category we work. 

If $\Mscr \in \Tails(A)$ then $\Hom(\Mscr, -)$ is left exact and since $\Tails(A)$ has enough injectives \cite{AZ} we may define its right derived functors $\Ext^{n}(\Mscr, -)$. We also use the notation 
\[ 
\underline{\Ext}^{n}(\Mscr,\Nscr) := \bigoplus_{d \in 
\ZZ}^{}{\Ext^{n}(\Mscr,\Nscr(d))} 
\] 
and we set $\underline{\Hom}(\Mscr,\Nscr) = \underline{\Ext}^{0}(\Mscr,\Nscr)$. 

Convention \ref{XyUvW} fixes the meaning of $\grmod(A)$, $\tors(A)$ and $\tails(A)$. It is easy to see $\tors(A)$ consists of the finite dimensional graded $A$-modules. Furthermore $\tails(A)=\grmod(A)/\tors(A)$. For $M \in \grmod(A)$ we have
\begin{equation} \label{ref-equation1.2} 
\underline{\Ext}^{n}(\pi M,\pi N) \cong \lim_{\lr} \underline{\Ext}^{n}_{A}(M_{\geq m},N)
\end{equation} 
for all $N \in \GrMod(A)$. If $M$ and $N$ are both finitely generated, then \eqref{ref-equation1.2} implies 
\[ 
\pi M \cong \pi N \mbox{ in } \tails(A) \hspace{0.25cm} \Leftrightarrow 
\hspace{0.25cm} M_{\geq n} \cong 
N_{\geq n} 
\mbox{ in $\grmod(A)$ for some $n \in \ZZ$} 
\] 
explaining the word ``tails''. The right derived functors of $\tau$ are given by $R^{i}\tau = \lim \underline{\Ext}_{A}^{i}(A/A_{\geq n},-)$ and for $M \in \GrMod(A)$ there is an exact sequence 
\begin{eqnarray} \label{ref-equation1.3} 
0 \ra \tau M \ra M \ra \omega \pi M \ra R^1\tau M \ra 0. 
\end{eqnarray} 



\subsection{Projective schemes}

Let $A$ be a noetherian graded $k$-algebra. As suggested by Artin and Zhang \cite{AZ}, define the (polarized) projective scheme $\Proj A$ of $A$ as the triple $(\Tails(A),\Oscr,\sh)$. In analogy with classical projective schemes we use the notation $\coh(X) := \tails(A)$, $\Qcoh(X) := \Tails(A)$. We write $\Ext^{i}_{X}(\Mscr,\Nscr)$ for $\Ext^{i}_{\Tails(A)}(\Mscr,\Nscr)$. We define the {\it cohomology groups} of $\Mscr \in \Qcoh(X)$ by $H^{n}(X,\Mscr) := \Ext^{n}_{X}(\Oscr_{X},\Mscr)$. We refer to the graded right $A$-modules 
\[ 
\underline{H}^{n}(X,\Mscr) := \bigoplus_{d \in \ZZ}^{}{H^{n}(X,\Mscr(d))} 
\] 
as the {\it full cohomology modules} of $\Mscr$. The {\it cohomological dimension} of $X$ is defined as $\cd X := \max \{ n \in \NN \mid  H^{n}(X, -) \neq 0 \}$. Clearly for $i > \cd X$ we have $\Ext^{i}_X(\Mscr,\Nscr) = 0$ for all $\Mscr,\Nscr \in \Qcoh(X)$.

\subsection{Hilbert series and Grothendieck group}
The {\em Hilbert series} of a graded $k$-vector space $V$ having finite dimensional components is the formal power series
\[ 
h_{V}(t) = \sum_{i = - \infty}^{+ \infty}{(\dim_{k}V_{i})t^{i}} \in\ZZ((t)). 
\] 
Let $A$ be a noetherian connected graded $k$-algebra. Then the Hilbert series $h_{M}(t)$ of $M \in \grmod(A)$ makes sense since $A$ is right noetherian. Note $h_{k}(t) = 1$ and $h_{M(l)}(t) = t^{-l}h_{M}(t)$. Assume moreover $A$ has finite global dimension. We denote by $\pd M$ the projective dimension of $M \in \grmod(A)$. Given a projective resolution of $M \neq 0$
\[ 
0 \r P^{r} \r \dots \r P^{1} \r P^{0} \r M \r 0 
\] 
we have $h_{M}(t) = \sum_{i=0}^{r}{(-1)^{i}h_{P^{i}}(t)}$. Since $A$ is connected, left bounded graded right $A$-modules are projective if and only if they are free hence isomorphic to a sum of shifts of $A$. So if we write $P^{i} = \oplus_{j=0}^{r_{i}}{A(-l_{ij})}$ we obtain 
\begin{align} \label{ref-equation1.4}
q_{M}(t) := h_{M}(t)h_A(t)^{-1} = \sum_{i=0}^{r}(-1)^{i}\sum_{j=0}^{r_{i}}t^{l_{ij}} \in \ZZ[t,t^{-1}]
\end{align}
where $q_{M}(t)$ is the so-called \emph{characteristic polynomial of $M$}. Note $q_{M(l)} = t^{-l}q_{M}(t)$, $q_{A}(t) = 1$ and $q_{k}(t) = h_{A}(t)^{-1}$. 

\medskip

The {\em Grothendieck group} $K_{0}(\Cscr)$ of an abelian category $\Cscr$ is the abelian group generated by all objects of $\Cscr$ (we write $[A] \in K_{0}(\Cscr)$ for $A \in \Cscr$) and for which we define $[A] - [B] + [C] = 0$ whenever there is an exact sequence $0 \ra A \ra B \ra C \ra 0$ in $\Cscr$.
Assume furthermore $\Cscr$ is $k$-linear and $\Ext$-finite with finite global dimension. Then the following map defines a bilinear form on $K_{0}(\Cscr)$ called the \emph{Euler form} for $\Cscr$ 
\begin{align} \label{Eulerform}
\chi: K_{0}(\Cscr) \times K_{0}(\Cscr) \ra \ZZ: ([A],[B]) \mapsto \chi(A,B) = \sum_i (-1)^i\dim_{k} \Ext^i_\Cscr(A,B).
\end{align}
Put $X = \Proj A$. We will write $K_{0}(X)$ for the Grothendieck group $K_{0}(\coh(X))$ of $\coh(X)$. The shift functor on $\coh(X)$ induces a group automorphism 
\begin{align*}
\sh: K_{0}(X) \ra K_{0}(X): [\Mscr] \mapsto [\Mscr(1)].
\end{align*}
We may view $K_{0}(X)$ as a $\ZZ[t,t^{-1}]$-module with $t$ acting as the shift functor $\sh^{-1}$. 
In \cite{MS} it was shown that $K_{0}(X)$ may be described in terms of Hilbert series. 
\begin{theorem} \cite[Theorem 2.3]{MS} \label{MoriSmith} 
Let $A$ be a noetherian connected graded $k$-algebra of finite global dimension and set $X = \Proj A$. Then there is an isomorphism of $\ZZ[t,t^{-1}]$-modules
\begin{equation} \label{isoMS}
\begin{aligned}
\theta : K_{0}(X) & \xrightarrow{\cong} \ZZ[t,t^{-1}]\left/(q_{k}(t))\right. \\
[\Mscr] & \mapsto \overline{q_{M}(t)} \quad \text{ where $M \in \grmod(A)$, $\Mscr = \pi M$.}
\end{aligned}
\end{equation}
In particular $[\Oscr(n)]$ is sent to $t^{-n}$. 
\end{theorem}

\subsection{Cubic three dimensional Artin-Schelter regular algebras} \label{Cubic three dimensional Artin-Schelter regular algebras}

We now come to the definition of regular algebras, introduced by Artin and Schelter \cite{AS} in 1987. They may be considered as noncommutative analogues of polynomial rings.
\begin{definition} \cite{AS} \label{defAS}
A connected graded $k$-algebra $A$ is called an {\it Artin-Schelter 
regular algebra} (or AS-algebra for short) {\it of dimension $d$}  
if it has the following properties: 
\begin{enumerate} 
\item[(i)] $A$ has finite global dimension $d$; 
\item[(ii)]$A$ has polynomial growth 
i.e. there are positive real 
numbers $c,e$ such that $\dim_{k}A_{n} \leq cn^{e}$ for all positive integers $n$; 
\item[(iii)] $A$ is Gorenstein, meaning there is an integer $l$ for which 
\[ 
\underline{\Ext}_{A}^{i}(k_{A},A) \cong 
\left \{ 
\begin{array}{ll} 
_{A}k(l) & \mbox{ if $i = d$,}\\ 
0 & \mbox{ otherwise} 
\end{array} 
\right. 
\] 
where $l$ is called the {\it Gorenstein parameter} of $A$. 
\end{enumerate} 
\end{definition}
If $A$ is commutative then the condition (i) already implies $A$ is isomorphic to a polynomial ring $k[x_{1},\dots,x_{n}]$ with some positive grading, if the grading is standard then $n = l$. 

The Gorenstein property determines the full cohomology modules of $\Oscr$.
\begin{theorem} \cite{AZ} \label{fullcohomology} 
Let $A$ be a noetherian AS-algebra of dimension $d = n + 1$ and let $X = \Proj A$. Let $l$ denote the Gorenstein parameter of $A$. Then $\cd X = n$ and the full cohomology modules of $\Oscr = \pi A$ are given by 
\[ 
\underline{H}^{i}(X,\Oscr) \cong 
\left\{ 
\begin{array}{ll} 
A & \mbox{ if } i = 0 \\ 
0 & \mbox{ if } i \neq 0,n \\ 
A'(l) & \mbox{ if } i = n 
\end{array} 
\right. 
\] 
\end{theorem} 
There exists a complete classification for Artin-Schelter regular algebras of dimension $d \leq 3$, see \cite{AS,ATV1,ATV2,Steph1,Steph2}. They are all left and right noetherian domains with Hilbert series of a weighted polynomial ring $k[x_1,\dots,x_d]$. We will be concerned with the case where $d = 3$ and $A$ is generated in degree one. Then there are two possibilities \cite{AS}, either 
\begin{itemize}
\item
$k_{A}$ has a minimal resolution of the form
\[
0 \ra A(-3) \ra A(-2)^{3} \ra A(-1)^{3} \ra A \ra k_{A} \ra 0
\]
thus $A$ has three generators and three defining homogeneous relations in degree two. Hence $A$ is Koszul and the Gorenstein parameter is $l = 3$. We then refer to $A$ as a \emph{quadratic} AS-algebra.
\item
$k_{A}$ has a minimal resolution of the form
\begin{equation} \label{minimalres}
0 \ra A(-4) \ra A(-3)^{2} \ra A(-1)^{2} \ra A \ra k_{A} \ra 0
\end{equation}
thus  $A$ has two generators and two defining homogeneous relations in degree three. In this case $l = 4$. These algebras are called \emph{cubic} AS-algebras.
\end{itemize}
In this article we will restrict ourselves to cubic AS-algebras $A$ and we denote $X = \Proj A$. From \eqref{minimalres} we easily deduce
the Hilbert series of $A$
\begin{equation} \label{HilbseriesA}
h_{A}(t) = \frac{1}{(1-t)^{2}(1-t^{2})}
\end{equation}
which means that for all integers $n$
\begin{eqnarray*}
\dim_{k}A_{n} = 
\left\{
\begin{array}{ll}
(n+2)^{2}/4 & \text{ if $n \geq 0$ is even,} \\
(n+1)(n+3)/4 & \text{ if $n \geq 0$ is odd.}
\end{array}
\right.
\end{eqnarray*} 
The Hilbert series of the Veronese subalgebra $A^{(2)} = k \oplus A_{2} \oplus A_{4} \oplus \dots$ of $A$ is the same as the Hilbert series of the commutative ring $k[x_{0},x_{1},x_{2},x_{3}]\left/(x_{0}x_{1} - x_{2}x_{3})\right.$ which is the homogeneous coordinate ring of a quadratic surface (quadric) in $\PP^{3}$. Since \cite{V} $\Tails(A) \cong \Tails(A^{(2)})$ and $\cd X = 2$ we therefore think of $\Proj A$ as a {\em quantum quadric}, a noncommutative analogue of the quadric surface $\PP^{1} \times \PP^{1}$. See also \cite{MS}.

\begin{example} \label{exampleHeis}
Consider the first Weyl algebra
\[
A_1 = k\langle x,y\rangle/(xy-yx-1)
\]
and introduce a third variable $z$ of degree two which commutes with $x$ and $y$ and makes the relation $xy-yx-1$ homogeneous. We obtain the enveloping algebra of the Heisenberg-Lie algebra
\begin{equation} \label{relenvweyl}
\begin{aligned}
H_{c} & = k \langle x,y,z \rangle / (yz - zy,xz - zx,xy - yx - z) \\
& = k \langle x,y \rangle / (y^{2}x - 2yxy + xy^{2},x^{2}y - 2xyx + yx^{2}) \\
& = k \langle x,y \rangle / ([y,[y,x]],[x,[x,y]]) 
\end{aligned}
\end{equation}
It is easy to verify that $H_c$ is a cubic AS-algebra. We refer to $H_{c}$ as the {\em enveloping algebra} for short.
\end{example}
\begin{example} \label{typeA} 
The generic cubic AS-algebras are the so-called type A-algebras \cite{AS}, they are of the form $k\langle x,y \rangle/(f_{1},f_{2})$ where $f_{1},f_{2}$ are the cubic equations
\begin{equation} \label{eqntypeA}
\left\{ 
\begin{array}{l} 
f_{1} = ay^{2}x + byxy + axy^{2} + cx^{3}  \\ 
f_{2} = ax^{2}y + bxyx + ayx^{2} + cy^{3}
\end{array} \right. 
\end{equation} 
with $(a:b:c) \in \PP^{2}\setminus F$ where $F = \{ (a:b:c) \in \PP^{2} \mid  a^{2} = b^{2} = c^{2} \} \, \cup$ \linebreak $\{ (0:1:0),(0:0:1) \}$.
A generic subclass is given by the more restrictive condition $(a:b:c) \in \PP^{2} \setminus F'$ where 
\begin{align*}
F' & = \{ (a:b:c) \in \PP^{2} \mid abc = 0 \text{ or } b^2 = c^2 \text{ or } 4b^2c^2 = (4a^2 - b^2 - c^2)^2 \}.
\end{align*}
We will refer to cubic AS-algebras $A$ of type A for which $(a:b:c) \in \PP^{2} \setminus F'$ as {\em generic type} A. The particular choice of $F'$ will become clear in Example \ref{ellipticsmooth} below.
\end{example}
\begin{remark}
Thus the enveloping algebra $H_{c}$ from Example \ref{exampleHeis} is a cubic AS-algebra of type A for which $(a:b:c) = (1:-2:0)$ in \eqref{eqntypeA}. However since $abc = 0$, $H_{c}$ is not of generic type A.
\end{remark}

\subsection{Geometric data associated to cubic AS-algebras}
\label{Geometric data associated to a cubic AS-algebra}

As before, let $A$ be a cubic AS-algebra. As shown in \cite{ATV1,ATV2} $A$ is completely determined by geometric data $(E,\sigma,\Lscr)$ where $E$ is a divisor on $\PP^1 \times \PP^1$, $\sigma \in \Aut(E)$ and $\Lscr$ is a line bundle on $E$. In this section we briefly recall this correspondence. 

We start with writing the relations of $A$ as 
\begin{eqnarray} \label{relations}
\left(
\begin{array}{c}
f_{1} \\
f_{2}
\end{array}
\right)
=
M_{A} \cdot
\left(
\begin{array}{c}
x \\
y 
\end{array}
\right)
\end{eqnarray}
where $M_{A} = (m_{ij})$ has entries $m_{ij} \in A_{2}$. We consider the multilinearizations $\tilde{f}_1$ and $\tilde{f}_2$ of the relations. Let $\Gamma \subset \PP^1\times \PP^1 \times \PP^1$ be the locus of common zeroes of $\tilde{f}_1$ and $\tilde{f}_2$ (with its scheme structure). Define the projection
\begin{eqnarray*}
\pr_{12}: \PP^1\times \PP^1 \times \PP^1 \ra \PP^1 \times \PP^1: (q_1,q_2,q_3)\mapsto (q_1,q_2).
\end{eqnarray*}
It turns out \cite{ATV1} that $\Gamma$ is the graph of an automorphism $\sigma$ of $E=\pr_{12}(\Gamma)$ and there are two distinguished cases:
either  $E = \PP^1\times \PP^1$ in which case we call $A$ \emph{linear}, or $E$ is a divisor of bidegree $(2,2)$ in $\PP^1\times \PP^1$. In the latter case we say $A$ is \emph{elliptic}. Then $\sigma$ is of the form $\sigma(q_1,q_2)=(q_2,f(q_1,q_2))$ for some map $f:E\ra \PP^1$.
\begin{example} \label{HeisC}
Consider the enveloping algebra $H_c$. Then $E$ is given by $((x_0:y_0),(x_1:y_1)) \in \PP^1\times \PP^1$ satisfying the relation $(x_0y_1-x_1y_0)^2$, thus $E$ is the double diagonal on $\PP^1\times \PP^1$ i.e. the points $((x:y),(x+\epsilon:y+\epsilon))$ such that $\epsilon^2=0$. Thus $H_c$ is elliptic. One further computes $\sigma((x:y),(x+\epsilon:y+\epsilon))=((x+\epsilon:y+\epsilon),(x+2\epsilon:y+2\epsilon))$, whence $\sigma \in \Aut(E)$ has infinite order.
\end{example}
\begin{example} \label{ellipticsmooth}
Let $A$ be a cubic AS-algebra of type A. Then $E$ is given by all $((x_0:y_0),(x_1:y_1)) \in \PP^1\times \PP^1$ for which
\[
(c^{2}-b^{2})x_{0}y_{0}x_{1}y_{1} 
+ ax_{0}^{2}(cx_{1}^{2} - by_{1}^{2})
+ ay_{0}^{2}(cy_{1}^{2} - bx_{1}^{2}) = 0
\]
whence $A$ is elliptic. In particular $E$ is smooth unless $abc = 0$ or $b^2 = c^2$ or $4b^2c^2 = (4a^2 - b^2 - c^2)^2$, i.e. $E$ is smooth if and only if $A$ is of generic type A. In this case $\sigma$ is given by a translation under the group law of the elliptic curve $E$.
\end{example}
Let $j$ be the inclusion $E\hookrightarrow \PP^1\times \PP^1$ and put $\Oscr_E(1) := j^{\ast}\Oscr_{\PP^1}(1)$. Associated to the geometric data $(E,\sigma,\Oscr_E(1))$ is a so-called twisted homogeneous coordinate ring $B=B(E,\sigma,\Oscr_E(1))$, see \cite{ATV1,ATV2,AV} or the construction below. If $A$ is linear then $A\cong B$. If $A$ is elliptic then there exists a canonical normal element $g \in A_{4}$ such that the quotient ring $A/gA$ is isomorphic to the twisted homogeneous coordinate ring $B = B(E, \sigma, \Oscr_{E}(1))$. 

As in \cite{DV2} the fact that $A$ may be linear or elliptic or $E$ may be non-reduced presents notational problems and difficulties. We therefore define $C = E_{\red}$ if $A$ is elliptic and let $C$ be a $\sigma$ invariant line $\{p\}\times \PP^1$ in $\PP^{1} \times \PP^{1}$ (where $p$ is a point in $\PP^1$). The geometric data $(E,\sigma,\Oscr_{E}(1))$ then restricts to geometric data $(C,\sigma_{C},\Oscr_{C}(1))$. In the elliptic case, writing $E = \sum_i n_i C_i$ where $C_i$ are the irreducible components of the support of $E$ we have $C = E_{\red} = \sum_i C_i$ and the irreducible components $C_i$ of $C$ form a single $\sigma$-orbit. 

It may occur that the order of $\sigma$ is different from the order of $\sigma_{C}$, being the restriction of $\sigma$ to $C$. For instance when $A=H_c$ is the  enveloping algebra, $\sigma$ has infinite order but $\sigma_{C}$ is the identity.
\begin{warning}
To simplify further expressions we write $(C,\sigma,\Oscr_{C}(1))$ for the triple \linebreak $(C,\sigma_{C},\Oscr_{C}(1))$. Below we will often assume $\sigma$ has infinite order. By this we will always mean the automorphism $\sigma$ in the geometric data $(E,\sigma,\Oscr_{E}(1))$ has infinite order and {\em not} the restriction of $\sigma$ to $C$.
\end{warning}
We will now recall the construction of the twisted homogeneous coordinate ring $B(C,\sigma,\Oscr_{C}(1))$. To simplify notations we will write $\Lscr = \Oscr_{C}(1)$ and we denote the auto-equivalence $\sigma_\ast(- \otimes_C \Lscr)$ by $- \otimes \Lscr_\sigma$. It is easy to check \cite[(3.1)]{VdBSt} that for $n \geq 0$ we have $\Mscr \otimes (\Lscr_\sigma)^{\otimes n} = \sigma_{\ast}^{n}\Mscr \otimes_{C} \sigma_{\ast}^{n} \Lscr \otimes_{C} \sigma_{\ast}^{n-1} \Lscr \otimes_{C} \dots \otimes_{C} \sigma_{\ast}\Lscr$
and since $(- \otimes\Lscr_{\sigma})^{-1} = 
\sigma^{\ast}(-)\otimes_{C}\Lscr^{-1}$ one may compute a similar expression for $\Mscr \otimes (\Lscr_\sigma)^{\otimes n}$ for $n < 0$. For $\Mscr\in \Qcoh(X)$ put $\Gamma_\ast(\Mscr)=\oplus_{n \geq 0} \Gamma(C,\Mscr\otimes 
(\Lscr_\sigma)^{\otimes n})$ and $D = B(C,\sigma,\Lscr) \overset{\text{def}}{=} \Gamma_\ast(\Oscr_C)$. Now $D$ has a natural ring structure and $\Gamma_\ast(\Mscr)$ is a right $D$-module. In \cite[\S 5]{ATV2} it is shown that there is a surjective morphism $A \ra D: a \rightarrow \overline{a}$ of graded $k$-algebras whose kernel is generated by a normalizing element $h$. In the elliptic case $h$ divides $g$ and $D$ is a prime ring. However $D$ may not be a domain since $C$ may have multiple components $C_i$. 

By analogy with the commutative case we may say $\Proj A$ contains $\Proj D$ as a ``closed'' subscheme. Though the structure of $\Proj A$ is somewhat obscure, $\Proj D$ is well understood. Indeed, it follows from \cite{AV,AZ} that the functor $\Gamma_\ast : \Qcoh(C)\r \GrMod(D)$ defines an equivalence $\Qcoh(C)\cong \Tails(D)$. The inverse of this equivalence and its composition with $\pi:\GrMod(D)\r \Tails(D)$ are both denoted by $\widetilde{(-)}$.
\begin{notation}
It will be convenient below to let the shift functor $-(n)$ on $\coh(C)$ be the one obtained from the equivalence $\coh(C)\cong \tails(D)$ and \emph{not} the one coming from the embedding $j: E \hookrightarrow \PP^1 \times \PP^1$.
\end{notation}

\subsection{Dimension, multiplicity and linear modules}

Let $A$ be a cubic AS-algebra and let $0 \neq M \in \grmod(A)$. As shown in \cite{ATV2} we may compute the Gelfand-Kirillov 
dimension $\gkdim M$ as the order of the pole of $h_{M}(t)$ at $t = 1$. We sometimes refer to $\gkdim M$ as the {\em dimension} of $M$. Note \eqref{HilbseriesA} implies $\gkdim A = 3$. As usual, a module $M \in \grmod(A)$ is {\em Cohen-Macaulay} if $\underline{\Ext}_{A}^{i}(M,A) = 0$ for $i \neq 3 - \gkdim M$, or equivalently $\pd M = 3- \gkdim M$. We then denote $M^{\vee} = \underline{\Ext}_{A}^{\pd M}(M,A)$. Define for $0 \neq M \in \grmod(A)$ 
\begin{equation} \label{defmultiplicity}
e_{n}(M) = \lim_{t \ra 1}(1-t)^{n}h_{M}(t).
\end{equation}
We have $e_{n}(M) \geq 0$ and $e_{n}(M) = 0$ if and only if $\gkdim M < n$. We define $\rank M = e_{3}(M)/e_3(A)$. For $n = \gkdim M$ we put $e(M) = e_{n}(M)$ and call this the {\em multiplicity} of $M$. Thus $e(M)$ is the first nonvanishing coefficient of the expansion of $h_{M}(t)$ in powers of $1-t$. 
For $0 \neq \Mscr \in \tails(A)$ we put $\dim \Mscr = \gkdim M - 1$, $e(\Mscr) = e(M)$, $\rank \Mscr = \rank M$
where $M \in \grmod(A)$, $\pi M = \Mscr$. 

An object in $\grmod(A)$ or $\tails(A)$ is said to be {\em pure} if it contains no subobjects of strictly smaller dimension, and it is called {\em critical} if every proper quotient has lower dimension.
Note $A$ is critical and for a critical $A$-module $M$ we have $\Hom_{A}(M,M) = k$, see \cite[Proposition 2.30]{ATV2}. We will often use 
\begin{lemma} \label{purepure}
\begin{enumerate}
\item
If $M \in \grmod(A)$ is pure (resp. critical) then $\pi M \in \tails(A)$ is pure (resp. critical). 
\item
If $\Mscr \in \tails(A)$ is pure (resp. critical) then $\Mscr = \pi M$ for some pure (resp. critical) object in $\grmod(A)$.
\item
Let $M,N \in \grmod(A)$ (resp. $\tails(A)$) be of the same dimension and assume $M$ is critical and $N$ is pure. Then every non-zero morphism in $\Hom(M,N)$ is injective.
\end{enumerate}
\end{lemma}
\begin{proof}
Elementary, see for example \cite[Lemma 1.9.3]{D}.
\end{proof}
A \emph{linear module} of dimension $d$ over $A$ is a cyclic graded right $A$-module generated in degree zero with Hilbert series $(1-t)^{-d}$. If $d = 1$ such a module is called a {\em point module}. Assume furthermore $A$ is elliptic. Then the map $p \mapsto \Gamma_{\ast}(\Oscr_p)$ defines a bijection between the points of $C$ (hence the closed points of $E$) and the point modules over $A$. We will denote $N_p = \Gamma_{\ast}(\Oscr_p)$ and $\Nscr_p = \pi N_p$. In particular all point modules $N_p$ over $A$ are $D$-modules i.e. $N_p \, h = 0$. We have $N_p(1)_{\geq 0} = N_{\sigma p}$ and a minimal resolution for $N_p$ is of the form
\begin{equation} \label{respoint}
0\ra A(-3) \ra A(-2)\oplus A(-1) \ra A \ra N_p \ra 0.
\end{equation}
Point modules are critical modules of GK-dimension zero. In case $\sigma$ has infinite order the converse is also true, up to shift of grading \cite{ATV2}.


We refer to linear modules of GK-dimension 2 as \emph{conic modules}. They are of the form $Q = A/vA$ where $v\in A_2$. A minimal resolution for $Q$ is of the form
\begin{equation} \label{resconic}
0\ra A(-2) \ra A \ra Q \ra 0.
\end{equation}
The name ``conic module" is justified from the fact that the multilinearization of $v\in A_2$ determines a curve in $\PP^3$ which is the intersection of an hyperplane and a quadric (the embedding of $\PP^1 \times \PP^1$).  

Finally, a {\em line module} over $A$ is a quotient module $A/uA$ with $u \in A_1$. There is a bijective correspondence between line modules $A/uA$ and lines $\{u=0\} \times \PP^1$ in $\PP^1\times \PP^1$. A line module $S = A/uA$ has a minimal resolution of the form
\begin{equation} \label{resline}
0\ra A(-1) \ra A \ra S \ra 0.
\end{equation}
Clearly point, line and conic modules over $A$ are Cohen-Macaulay modules.

\subsection{Serre duality}

Let $A$ be a cubic AS-regular algebra. As shown in \cite{DV1} there is a graded automorphism $\phi$ of $A$, passing to an automorphism $(-)_{\phi}$ on $\Tails(A)$, for which there are natural isomorphisms
\begin{equation} \label{Serre}
\Ext^i_{D^b(\coh(X))}(\Mscr,\Nscr)\cong \Ext^{n-i}_{D^b(\coh(X))}(\Nscr,\Mscr_{\phi}(-4))' 
\end{equation}
for all objects $\Mscr, \Nscr$ of the bounded derived category $D^b(\coh(X))$ of $\coh(X)$. However we will look for an algebra $\widehat{A}$ for which $\GrMod(A) \cong \GrMod(\widehat{A})$ and for which \eqref{Serre} takes a more simple form. 

As in \eqref{relations} we write the relations $f$ of $A$ as $f = M_A x$. With a suitable choice of the relations $f$ we have $x^t M_A = (Q_A f)^t$ for some invertible matrix $Q_A$ with scalar entries, see \cite[Theorem 1.5]{AS}. It now turns out there exists a Zhang-twist \cite{Zhang} $A^{\tau}$ of $A$ for which $Q_{A^{\tau}}$ is the identity matrix. This was pointed out to us by M. Van den Bergh, see also \cite{ncquadrics}. By \cite{Zhang} we have $\GrMod(A) \cong \GrMod(A^{\tau})$ and $\Tails(A) \cong \Tails(A^{\tau})$ where $(\pi A)(n)$ is sent to $(\pi A^{\tau})(n)$.

If $A$ is of type A then writing the relations as in \eqref{eqntypeA} yields $Q_A = \Id$ whence we may put $A = A^{\tau}$. 
\begin{convention} \label{conventionSerre}
Whenever $A$ is a cubic Artin-Schelter algebra we replace $A$ with a Zhang-twist $A^{\tau}$ for which $Q_{A^{\tau}}$ is the identity matrix. 
\end{convention}
\begin{remark}
We are allowed to use Convention \ref{conventionSerre} in this paper since we will only specify to linear or elliptic algebras for which we often require $\sigma$ has infinite order (but these properties are invariant under Zhang-twisting) and we will not rely on the relations of $A$ except for specific relations \eqref{eqntypeA} for algebras of type A and in particular \eqref{relenvweyl} for the enveloping algebra.
\end{remark}
Using Convention \ref{conventionSerre} we see \eqref{Serre} takes a particularly simple form.
\begin{theorem} \label{Serreduality} (Serre duality)
Let $\Mscr$, $\Nscr \in D^b(\coh(X))$. Then there are natural isomorphisms 
\[
\Ext^i_{D^b(\coh(X))}(\Mscr,\Nscr)\cong \Ext^{n-i}_{D^b(\coh(X))}(\Nscr,\Mscr(-4))'.
\]
\end{theorem}
Below we refer to Theorem \ref{Serreduality} as Serre duality on $X$.

\section{From reflexive ideals to normalized line bundles}

Throughout $A$ will be a cubic AS-algebra as defined in \S\ref{Cubic three dimensional Artin-Schelter regular algebras}. We will use the notations from the previous section, so we write $X = \Proj A$, $\Qcoh(X) = \Tails(A)$, $\coh(X) = \tails(A)$, $\pi A = \Oscr$.

In this section our first aim is to relate reflexive $A$-modules with certain objects on $X$ (so-called vector bundles). Any shift of such a reflexive module remains reflexive and in the rank one case we will normalize this shift. The corresponding objects in $\coh(X)$ will be called normalized line bundles. A helpful tool will be the choice of a suitable basis of the Grothendieck group $K_0(X)$. At the end of this section we will compute partially the cohomology of these normalized line bundles.

\subsection{Reflexive modules and vector bundles}
An object $M \in \grmod(A)$ is {\em torsion free} if $M$ is pure of maximal GK-dimension three. Recall $M$ is called reflexive if $M^{\ast \ast} = M$. Similarly an object $\Mscr \in \coh(X)$ is {\em torsion free} if $\Mscr$ is pure of maximal dimension two. An object $\Mscr \in \coh(X)$ is called \emph{reflexive} (or a \emph{vector bundle} on $X$) if $\Mscr = \pi M$ for some reflexive $M \in \grmod(A)$. We refer to a vector bundle of rank one as a \emph{line bundle}. We will need the following lemmas, see \cite[Lemma 3.4]{DV1} and \cite[Proposition 3.4.1 and Corollary 3.4.2]{DV2}.
\begin{lemma} \label{reflexive}
Let $\Mscr \in \coh(X)$. Then $\Mscr$ is a vector bundle on $X$ if and only if $\Mscr$ is torsion free and $\Ext^{1}_{X}(\Nscr,\Mscr) = 0$ for all $\Nscr \in \coh(X)$ of dimension zero.
\end{lemma}
\begin{lemma} \label{equivalences}
The functors $\pi$ and $\omega$ define inverse equivalences between the full subcategories of $\grmod(A)$ and $\coh(X)$ with objects
\[
\{ \text{torsion free objects in }  \grmod(A) \text{ of projective dimension one}  \} 
\]
and
\[
\{\text{torsion free objects in } \coh(X) \}
\] 
Moreover this equivalence restricts to an equivalence between the full subcategories of $\grmod(A)$ and $\coh(X)$ with objects
\[
\{\text{reflexive objects in }\grmod(A) \} \quad \text{ and } \quad \{\text{vector bundles on } X \}.
\]
\end{lemma}
In this paper we are interested in torsion free rank one modules of projective dimension one, or more restrictively, reflexive modules of rank one. Every graded right ideal of $A$ is a torsion free rank one $A$-module. The following proposition shows that, up to shift of grading, the converse is also true.
\begin{proposition}
Let $0\neq I \in \grmod(A)$ be torsion free of rank one. Then there is an integer $n$ such that $I(-n)$ is isomorphic to a graded right ideal of $A$.
\end{proposition}
\begin{proof}
By $\gkdim I = 3$, \cite[Theorem 4.1]{ATV2} implies $I^{\ast}=\underline{\Hom}_A(I,A) \neq 0$. Thus $(I^{\ast})_n = \Hom_A(I(-n),A)\neq 0$ for some integer $n$. By Lemma 
\ref{purepure} we are done. 
\end{proof}
\begin{remark} 
The set of all graded right ideals is probably too large to describe, as for any ideal $I$ we may construct numerous other closely related ideals by taking the kernel of any surjective map to a module of GK-dimension zero. 
We will restrict to graded ideals of projective dimension one (or more restrictively reflexive rank one modules). For such modules $M$ we have $\underline{\Ext}^1_A(k,M) = 0$ and therefore $M$ cannot appear as the kernel of such a surjective map.
\end{remark}

\subsection{The Grothendieck group of $X$} \label{The Grothendieck group of $X$}

In this part we describe a natural $\ZZ$-module basis for the Grothendieck group $K_{0}(X)$ and determine the matrix representation of the Euler form $\chi$ with respect to this basis. 
To do so, it is convenient to start with a different basis of $K_{0}(X)$, corresponding to the standard basis of $\ZZ[t,t^{-1}]/(q_{k}(t))$ under the isomorphism of Theorem \ref{MoriSmith}, and perform a base change afterwards.
\begin{proposition} \label{Grothendieck1}
The set $\Bscr = \{[\Oscr],[\Oscr(-1)],[\Oscr(-2)],[\Oscr(-3)]\}$ is a $\ZZ$-module basis of $K_{0}(X)$. The matrix representations with respect to the basis $\Bscr$ of the shift automorphism $\sh$ and the Euler form $\chi$ for $K_0(X)$ are given by
\begin{eqnarray*} 
\quad \quad m(\sh)_{\Bscr} = 
\left(
\begin{array}{rrrr}
2 & 1 & 0 & 0 \\
0 & 0 & 1 & 0 \\
-2 & 0 & 0 & 1 \\
1 & 0 & 0 & 0
\end{array}
\right), \quad m(\chi)_{\Bscr} = 
\left(
\begin{array}{rrrr}
1 & 0 & 0 & 0 \\
2 & 1 & 0 & 0 \\
4 & 2 & 1 & 0 \\
6 & 4 & 2 & 1
\end{array}
\right).
\end{eqnarray*}
\end{proposition}
\begin{proof}
Let $\theta$ denote the isomorphism \eqref{isoMS} of Theorem \ref{MoriSmith}. Since $q_{A(-l)}(t) = t^{l}$ we have $\theta[\Oscr(-l)] = \overline{t^{l}}$ for all integers $l$. As $\{ \overline{1},\overline{t},\overline{t^{2}},\overline{t^{3}} \}$ is a $\ZZ$-module basis for $\ZZ[t,t^{-1}]/(q_{k}(t)) = \ZZ[t,t^{-1}]/(1-t)^{2}(1-t^{2})$ we deduce $\Bscr$ is a basis for $K_0(X)$.

By $\sh[\Oscr(l)] = [\Oscr(l+1)]$ we find the last three columns of $m(\sh)_{\Bscr}$. Applying the exact functor $\pi$ to \eqref{minimalres} yields the exact sequence
\[
0 \ra \Oscr(-4) \ra \Oscr(-3)^{2} \ra \Oscr(-1)^{2} \ra \Oscr \ra 0
\]
from which we deduce $[\Oscr(1)] = 2[\Oscr] - 2[\Oscr(-2)] + [\Oscr(-3)]$, giving the first column of $m(\sh)_{\Bscr}$. Finally, Theorem \ref{fullcohomology} implies for all integers $l$
\begin{eqnarray*}
\quad \quad \quad \quad \chi(\Oscr, \Oscr(l)) = \dim_{k}A_{l} + \dim_{k}A_{-l-4} = 
\left\{
\begin{array}{ll}
(l+2)^{2}/4 & $ if $ l $ is even$ 
\\
(l+1)(l+3)/4 & $ if $ l $ is odd$
\end{array} 
\right.
\end{eqnarray*}
which allows one to compute the matrix $m(\chi)_{\Bscr}$. This ends the proof.
\end{proof}
\begin{proposition} \label{basis2}
Let $P$ be a point module, $S$ a line module and $Q$ a conic module over $A$. Denote the corresponding objects in $\coh(X)$ by $\Pscr$, $\Sscr$ and $\Qscr$. Then $\Bscr' = \{ [\Oscr],[\Sscr],[\Qscr],[\Pscr]\}$ is a $\ZZ$-module basis of $K_{0}(X)$, which does not depend on the particular choice of $S$, $Q$ and $P$. The matrix representations with respect to the basis $\Bscr'$ of the shift automorphism $\sh$ and the Euler form $\chi$ for $K_0(X)$ are given by
\begin{eqnarray} 
\quad \quad m(\sh)_{\Bscr'} = 
\left(
\begin{array}{rrrr}
1 & 0 & 0 & 0 \\
-1 & -1 & 0 & 0 \\
1 & 1 & 1 & 0 \\
1 & 1 & 1 & 1
\end{array}
\right), \quad 
m(\chi)_{\Bscr'} = 
\left(
\begin{array}{rrrr}
1 & 1 & 1 & 1 \\
-1 & 0 & -1 & 0 \\
-3 & -1 & -2 & 0 \\
1 & 0 & 0 & 0
\end{array}
\right).
\end{eqnarray}
\end{proposition}
\begin{proof}
Easy by Proposition \ref{Grothendieck1}, base change and equations \eqref{respoint}-\eqref{resline}.
\end{proof}
From now on we fix such a $\ZZ$-module basis $\{[\Oscr],[\Sscr],[\Qscr],[\Pscr]\}$ of $K_{0}(X)$. For any object $\Mscr \in \coh(X)$ we may write 
\begin{equation} \label{expK0}
[\Mscr] = r[\Oscr] + a[\Sscr] + b[\Qscr] + c[\Pscr]
\end{equation}
Writing $\Mscr = \pi M$ where $M \in \grmod(A)$, equation \eqref{expK0} also follows directly from Theorem \ref{MoriSmith} and the fact that we may expand the characteristic polynomial $q_M(t)$ of $M$ as
\[
q_M(t) = r + a(1-t) + b(1-t^2) + c(1-t^2)(1-t) + f(t)(1-t^2)(1-t)^2
\] 
for some integers $a,b,c \in \ZZ$ and some Laurent polynomial $f(t)\in \ZZ[t,t^{-1}]$. Combining \eqref{ref-equation1.4} and \eqref{HilbseriesA} yields
\begin{equation} \label{hilbertseries1}
h_{M}(t)=\frac{r}{(1-t)^{2}(1-t^{2})}+\frac{a}{(1-t)(1-t^{2})}+\frac{b}{(1-t)^{2}}+\frac{c}{1-t} + f(t).
\end{equation}
Note $r = \rank M = \rank \Mscr$. By computing the powers of the matrix $m(\sh)_{\Bscr'}$ in Proposition \ref{basis2} we deduce for any integer $l$
\begin{equation} \label{shiftcubic}
\begin{aligned}
\left[\Mscr(2l) \right] & = r[\Oscr] + a[\Sscr] + (lr + b)[\Qscr]  +(l((l+1)r + a + 2b) + c)[\Pscr] \\
\left[\Mscr(2l-1)\right] & = r[\Oscr] - (r+a)[\Sscr] + (lr + a + b)[\Qscr]  + (l(lr + a + 2b) - b + c)[\Pscr]
\end{aligned}
\end{equation}

\subsection{Normalized line bundles} \label{Normalized line bundles}

Any shift $l$ of a torsion free rank one graded right $A$-module $I$ gives rise to a torsion free rank one object $\Iscr(l) = \pi I(l)$ on $X$. We will now normalize this shift. Our choice is motivated by 
\begin{proposition} \label{normalized}
Let $I \in \grmod(A)$, set $\Iscr = \pi I$ and write $[\Iscr] = r[\Oscr] + a[\Sscr] + b[\Qscr] + c[\Pscr]$.
Then the following are equivalent.
\begin{enumerate}
\item
There exist integers $n_{e}, n_{o}$ such that for $l \gg 0$ we have
\begin{eqnarray*}
\dim_{k}A_{l} - \dim_{k}I_{l} = 
\left\{
\begin{array}{ll}
n_{e} & \text{ if $l$ is even,} \\
n_{o} & \text{ if $l$ is odd.} 
\end{array}
\right.
\end{eqnarray*}
\item
The Hilbert series of $I$ is of the form
\[
h_I(t) = 
h_A(t) - \frac{s(t)}{1-t^2} 
\]
for a Laurent polynomial $s(t) \in \ZZ[t,t^{-1}]$.
\item
$I$ has rank one and $a = -2b$.
\end{enumerate}
If these conditions hold then $s(1) = b-2c$, $s(-1) = b$ and $n_{e} = b-c$, $n_{o} = -c$.
\end{proposition}
\begin{proof}
By \eqref{hilbertseries1} we may write
\[
h_{I}(t) = \frac{r}{(1-t)^{2}(1-t^2)} + \frac{a + b(1+t)}{(1-t)(1-t^2)} + \frac{c(1+t) + f(t)(1-t^2)}{1-t^2} 
\]
for some $f(t) \in \ZZ[t,t^{-1}]$. Thus the second and the third statement are equivalent, and in that case $s(t) = b - c(1+t) - f(t)(1-t^2)$. Moreover, for $l \gg 0$ we obtain
\begin{eqnarray}
\dim_{k}A_{l} - \dim_{k}I_{l} = 
\left\{
\begin{array}{ll}
(1-r)(l+2)^{2}/4 - a(l/2+1) - b(l+1) - c & \text{\hspace{-0.4cm} for $l$ even} \\
(1-r)(l+1)(l+3)/4 - a(l+1)/2 - b(l+1) - c & \text{\hspace{-0.3cm} for $l$ odd} 
\end{array}
\right. \nonumber
\end{eqnarray}
from which we deduce the equivalence of (1) and (3), proving what we want.
\end{proof} 

We will call a torsion free rank one object in $\grmod(A)$ {\em normalized} if it satisfies the equivalent conditions of Proposition \ref{normalized}. Similarly, a torsion free rank one object $\Iscr$ in $\coh(X)$ is {\em normalized} if $[\Iscr]$ is of the form
\[
[\Iscr] = [\Oscr] - 2b[\Sscr] + b[\Qscr] + c[\Pscr]
\]
for some integers $b,c \in \ZZ$. We refer to $(n_{e},n_{o}) = (b-c,-c)$ as the {\em invariants} of $I$ and $\Iscr$ and call $n_e$ the {\em even invariant} and $n_o$ the {\em odd invariant} of $I$ and $\Iscr$. We will prove in Theorem \ref{cohomology} below that $n_e$ and $n_o$ are actually positive and characterize the appearing invariants $(n_e,n_o)$ in Section \ref{Hilbert series of torsion free rank one modules}. 
\begin{lemma}
Let $I \in \grmod(A)$ be torsion free of rank one and set $\Iscr = \pi I$. Then there is a unique integer $d$ for which $I(d)$ (and hence $\Iscr(d)$) is normalized.
\end{lemma}
\begin{proof}
Easy by \eqref{shiftcubic}.
\end{proof}
By Lemma \ref{equivalences} the functors $\pi$ and $\omega$ define inverse equivalences between the full subcategories of $\grmod(A)$ and $\coh(X)$ with objects
\begin{multline*} \label{Hilb}
\Hilb_{(n_e,n_o)}(X) := \{\text{normalized torsion free rank one objects in $\grmod(A)$} \\
\text{of projective dimension one and invariants $(n_e,n_o)$}\} 
\end{multline*}
and
\[ 
\{\text{normalized torsion free rank one objects in $\coh(X)$ with invariants $(n_e,n_o)$} \}. 
\] 
\begin{remark}
We expect $\coprod_{(n_e,n_o)}\Hilb_{(n_e,n_o)}(X)$ to be the correct generalization of the usual Hilbert scheme of points on $\PP^1 \times \PP^1$. In case $A$ is linear then \linebreak $\coprod_{(n_e,n_o)}\Hilb_{(n_e,n_o)}(X)$ coincides with the Hilbert scheme of points on $\PP^1 \times \PP^1$, see \S\ref{Hilbert scheme of points for the quadric surface} below. 
\end{remark}
This equivalence restricts to an equivalence between the full subcategories of \linebreak $\grmod(A)$ and $\coh(X)$ with objects
\begin{multline*}
R_{(n_e,n_o)}(A) := \{\text{normalized reflexive rank one objects in $\grmod(A)$} \\  
\text{with invariants $(n_e,n_o)$} \} 
\end{multline*}
and
\[ 
\Rscr_{(n_e,n_o)}(X) := \{\text{normalized line bundles on $X$ with invariants $(n_e,n_o)$} \}.
\] 
We obtain a natural bijection between the set $R(A)$ of reflexive rank one graded right $A$-modules considered up to isomorphism and shift, and the isomorphism classes in the categories $\coprod_{(n_e,n_o)}R_{(n_e,n_o)}(A)$ and $\coprod_{(n_e,n_o)}\Rscr_{(n_e,n_o)}(X)$.
\begin{remark}
It is easy to see that the categories $R_{(n_e,n_o)}(A)$ and $\Rscr_{(n_e,n_o)}(X)$ are groupoids, i.e. all non-zero morphisms are isomorphisms.
\end{remark}

\subsection{Cohomology of normalized line bundles}

The next theorem describes partially the cohomology of normalized line bundles.
\begin{theorem}  \label{cohomology} 
Let $\Iscr \in \coh(X)$ be torsion free of rank one and normalized i.e. 
\[
[\Iscr] = [\Oscr] -2(n_{e} - n_{o})[\Sscr] + (n_{e} - n_{o})[\Qscr] - n_{o}[\Pscr]
\]
for some integers $n_e,n_o$. Assume $\Iscr$ is not isomorphic to $\Oscr$. Then 
\begin{enumerate} 
\item 
$H^{0}(X,\Iscr (l)) = 0 \mbox{ for } l \leq 0 $ \\ 
$H^{2}(X,\Iscr (l)) = 0 \mbox{ for } l \geq -3$ \\
$H^j(X,\Iscr(l))=0$ for $j\ge 3$ and for all integers $l$
\item 
$\chi (\Oscr,\Iscr(l)) = \left\{
\begin{array}{ll}
(l+2)^{2}/4 - n_{e} & \text{ if $l \in \ZZ$ is even } \\
(l+1)(l+3)/4 - n_{o} & \text{ if $l \in \ZZ$ is odd } \\
\end{array}
\right.$ 
\item 
$\dim_{k} H^{1}(X,\Iscr) = n_e - 1$ \\ 
$\dim_{k} H^{1}(X,\Iscr(-1)) = n_o$ \\ 
$\dim_{k} H^{1}(X,\Iscr(-2)) = n_e$ \\
$\dim_{k} H^{1}(X,\Iscr(-3)) = n_o$
\end{enumerate} 
As a consequence, $n_e > 0$ and $n_o \geq 0$. If $\Iscr$ is a line bundle i.e. $\Iscr \in \Rscr_{(n_e,n_o)}(X)$ then we have in addition
\[ 
H^{2}(X, \Iscr (-4)) = 0 \text{ and } \dim_{k} H^{1}(X,\Iscr(-4)) = n_e - 1.
\]
\end{theorem}
\begin{proof}
That $H^j(X,\Iscr(l))=0$ for $j\ge 3$ and for all integers $l$ follows from $\cd X = 2$, see Theorem \ref{fullcohomology}. The rest of the first statement is proved in a similar way as \cite[Theorem 3.5(1)]{DV2}. See also the proof of the final statement below.

For the second part, compute $\chi (\Oscr,\Iscr(l))$ using \eqref{shiftcubic} and the matrix representation $m(\chi)_{\Bscr'}$ from Proposition \ref{basis2}. 

Combining the first two statements together with \eqref{Eulerform} yields the third part.

Finally, assume $\Iscr$ is reflexive. By Theorem \ref{Serreduality} (Serre duality) we have \linebreak $H^{2}(X, \Iscr (-4)) = \Ext^{2}_X(\Oscr,\Iscr(-4)) \cong \Hom_{X}(\Iscr,\Oscr)'$. Assume by contradiction there is a non-zero morphism $f: \Iscr \r \Oscr$. As $\Iscr$ is critical, $f$ is injective and we compute $[\coker f] = 2(n_{e} - n_{o})[\Sscr] - (n_{e} - n_{o})[\Qscr] + n_{o}[\Pscr]$. By \eqref{defmultiplicity} and \eqref{expK0}-\eqref{hilbertseries1} we deduce $e_1(\coker f) = 0$ hence $\dim \coker f = 0$. Note $\coker f \neq 0$ by the assumption $\Iscr \not\cong \Oscr$. Since $\Iscr$ is reflexive, $\Ext^{1}_{X}(\coker f, \Iscr) = 0$ thus the exact sequence $0 \r \Iscr \r \Oscr \r \coker f \r 0$ splits, contradicting the fact that $\Oscr$ is torsion free. 
\end{proof}

\begin{corollary} \label{invzero}
Let $I \in \grmod(A)$ be torsion free of rank one with invariants $(n_e,n_o)$. Then $(n_e,n_o) = (0,0)$ if and only if $\Iscr \cong \Oscr(d)$ for some integer $d$.
\end{corollary}
\begin{proof}
If $\Iscr \cong \Oscr(d)$ then $[\Iscr(-d)] = [\Oscr]$ hence $n_{e} = n_{o} = 0$. Assume conversely $(n_e,n_o)=(0,0)$. We may assume $\Iscr$ is normalized. If $\Iscr\not\cong \Oscr$ then Theorem \ref{cohomology} implies $n_e > 0$. Since $n_e = 0$ we obtain $\Iscr\cong \Oscr$ by contraposition. 
\end{proof} 
At this point one may be tempted to think there are two independent parameters $n_e,n_o \in \NN$ associated to an object in $\Hilb_{(n_e,n_o)}(X)$. However 
\begin{lemma} \label{dimensionExt1}
Let $I \in \grmod(A)$ be torsion free of rank one with invariants $(n_e,n_o)$ and write $\Iscr = \pi I$. Then $\dim_{k}\Ext^{1}_{X}(\Iscr,\Iscr) = 2(n_{e} - (n_{e} - n_{o})^{2}) \geq 0$.
\end{lemma}
\begin{proof}
We may clearly assume $\Iscr$ is normalized and by Proposition \ref{basis2} we easily find $\chi(\Iscr, \Iscr) = 1 - 2(n_e - (n_e - n_o)^2)$. As $\Iscr$ is critical we have $\Hom_{X}(\Iscr,\Iscr) = k$. Hence it will be sufficient to prove $\Ext^{2}_{X}(\Iscr,\Iscr) = 0$. Serre duality implies $\Ext^{2}_{X}(\Iscr,\Iscr) \cong \Hom_{X}(\Iscr,\Iscr(-4))'$. Thus assume by contradiction there is a non-zero morphism $f: \Iscr \r \Iscr(-4)$. Then $f$ is injective and using \eqref{shiftcubic} we have $[\Iscr(-4)] = [\Oscr] - 2(n_e - n_o)[\Sscr] + (n_e - n_o - 2)[\Qscr]  + (2 - n_o)[\Pscr]$ hence $[\coker f] = -2[\Qscr]  + 2[\Pscr]$. By \eqref{defmultiplicity} and \eqref{expK0}-\eqref{hilbertseries1} we deduce $e_2(\coker f) < 0$ which is absurd.
\end{proof} 
As a consequence if $\Hilb_{(n_e,n_o)}(X) \neq \emptyset$ for some integers $n_e,n_o$ then $n_e \geq 0$, $n_o\geq 0$ and $n_e - (n_e - n_o)^2 \geq 0$. The converse will be proved in the next section.

\section{Hilbert series of ideals and proof of Theorem \ref{theorem2}}
\label{Hilbert series of torsion free rank one modules}

Let $A$ be a quadratic or cubic AS-algebra and let $M$ be a torsion free graded right $A$-module of projective dimension one (so we do not require $M$ to have rank one). Thus $M$ has a minimal resolution of the form
\[
0 \r \oplus_i A(-i)^{b_i} \r \oplus_i A(-i)^{a_i} \r M \r 0
\]
where $(a_i),(b_i)$ are finitely supported sequences of non-negative integers. These numbers are called the {\em Betti numbers} of $M$. It is easy to see that the characteristic polynomial of $M$ is given by $q_M(t) = \sum_i(a_i - b_i)t^i$. So by \eqref{ref-equation1.4} the Betti numbers determine the Hilbert series of $M$, but the converse is not true. 

For quadratic $A$ the appearing Betti numbers were characterised in \cite[Corollary 1.5]{DV2}. The same technique as in \cite{DV2} may be used to obtain the same characterisation for cubic $A$. The result is
\begin{proposition} \label{prophilbert}
Let $(a_i),(b_i)$ be finitely supported sequences of non-negative integers. Let $a_{\sigma}$ be the lowest non-zero $a_i$ and put $r = \sum_{i}(a_i-b_i)$.
Then the following are equivalent.
\begin{enumerate}
\item
$(a_i),(b_i)$ are the Betti numbers of a torsion free graded right module of projective dimension one and rank $r$ over a quadratic AS-algebra,
\item
$(a_i),(b_i)$ are the Betti numbers of a torsion free graded right module of projective dimension one and rank $r$ over a cubic AS-algebra,
\item
$b_i = 0$ for $i \leq \sigma$ and $\sum_{i\leq l} b_i < \sum_{i < l}a_i$ for $l > \sigma$.
\end{enumerate}
Moreover if $A$ is elliptic and $\sigma$ has infinite order, these modules can be chosen to be reflexive.
\end{proposition}
Assume for the rest of Section \ref{Hilbert series of torsion free rank one modules} $A$ is a cubic AS-algebra. The previous proposition allows us to describe the Hilbert series of objects in $\Hilb_{(n_e,n_o)}(X)$. Recall from the introduction a Castelnuovo polynomial \cite{Davis} $s(t)=\sum_{i=0}^{n}s_it^i \in \ZZ[t]$ is by definition of the form 
\begin{equation} 
s_0=1,s_1=2,\ldots,s_{\sigma-1}=\sigma \mbox{ and } s_{\sigma-1}\geq s_{\sigma}\geq s_{\sigma+1}\geq \cdots \geq 0
\end{equation} 
for some integer $\sigma \geq 0$. We refer to $\sum_{i}s_{2i}$ as the  {\em even weight} of $s$ and $\sum_{i}s_{2i+1}$ as the {\em odd weight} of $s(t)$.
We may now prove Theorem \ref{theorem2}. 
\begin{proof}[Proof of Theorem \ref{theorem2}]
First, let us assume $I \in \Hilb_{(n_e,n_o)}(X)$ for some integers $n_e, n_o$. By Proposition \ref{normalized} we may assume that the Hilbert series of $I$ has the form 
\[
h_{I}(t) = \frac{1}{(1-t)^{2}(1-t^2)} - \frac{s(t)}{1-t^2} 
\]
for a Laurent polynomial $s(t) \in \ZZ[t,t^{-1}]$. We deduce $q_I(t)/(1-t) = h_{I}(t)(1-t)(1-t^2) = 1/(1-t) - s(t)(1-t)$. Writing $q_I(t) = \sum_i q_i t^i$ it is easy to see Proposition \ref{prophilbert}(3) is equivalent with 
\begin{eqnarray*}
\sum_{i \leq l}q_i
\left\{
\begin{array}{ll}
= 0 & \text{ for } l < \sigma \\
> 0 & \text{ for } l \geq \sigma
\end{array}
\right.
\end{eqnarray*}
from which we deduce $s(t)(1-t)$ is of the form
\[
s(t)(1-t) = 1 + t + t^2 + \dots + t^{\sigma-1} + d_{\sigma}t^{\sigma} + d_{\sigma+1}t^{\sigma+1} +\dots
\]
where $d_i \leq 0$ for $i \geq \sigma$. Multiplying by $1/(1-t) = 1 + t + t^2 + \dots$ shows this is equivalent to $s(t)$ being a Castelnuovo polynomial. According to Proposition \ref{normalized}, $(s(1)+s(-1))/2 = n_e$ and $(s(1)-s(-1))/2 = n_o$ thus $s(t)$ has even weight $n_e$ and odd weight $n_o$. 

The converse statement is easily checked.
\end{proof} 
As an application we may now prove nonemptyness for $R_{(n_e,n_o)}(A)$. As in the introduction we define 
\begin{equation} \label{defN}
N = \{(n_e,n_o) \in \NN^2 \mid n_e - (n_e - n_o)^2 \geq 0 \}.
\end{equation}
It is a simple exercise to check  
\begin{equation} \label{descrN}
N = \{(k^2+l, k(k+1)+l) \mid k,l \in \NN \} \cup \{((k+1)^2+l, k(k+1)+l) \mid k,l \in \NN \}.
\end{equation}
\begin{proposition} \label{nonemptycubic}
Let $n_e$,$n_o$ be any integers. Then $\Hilb_{(n_e,n_o)}(X)$ is nonempty if and only if $(n_e,n_o) \in N$. \\
If $A$ is elliptic and $\sigma$ has infinite order then $R_{(n_e,n_o)}(A)$ whence $\Rscr_{(n_e,n_o)}(X)$ is nonempty if and only if $(n_e,n_o) \in N$.
\end{proposition}
\begin{proof}
Assume $(n_e,n_o) \in N$. Due to Theorem \ref{theorem2} it will be sufficient to show there exists a Castelnuovo polynomial $s(t)$ for which the even resp. odd weight of $s(t)$ is equal to $n_e$ resp. $n_o$. Shifting the rows in any Castelnuovo diagram in such a way they are left aligned induces a bijective correspondence between Castelnuovo functions $s$ and partitions $\lambda$ of $n = s(1)$ with distinct parts. For any partition $\lambda$ we put a chess colouring on the Ferrers graph of $\lambda$, and write $b(\lambda)$ resp. $w(\lambda)$ for the number of black resp. white unit squares. The result follows from the known fact that there exists a partition $\lambda$ in distinct parts for which $b(\lambda) = n_e$ and $w(\lambda) = n_o$ if and only if $(n_e,n_o) \in N$. See for example \cite{DM}.
\end{proof}
For $(n_e,n_o) \in N$ there is an unique integer $l \geq 0$ with the property (see \eqref{descrN})
\begin{equation} \label{uniquel}
(n_e-l,n_o-l) \in N \text{ and }(n_e-l-1,n_o-l-1) \not\in N. 
\end{equation}
One verifies $(n_e-l',n_o-l') \not\in N$ for all $l' > l$. By \eqref{descrN} we distinguish \cite{DM}
\setcounter{case}{0}
\begin{case}
$(n_e-l,n_o-l) = (k^2,k(k+1))$ for $k \in \NN$. The Castelnuovo polynomial of an object in $\Hilb_{(n_e-l,n_o-l)}(X)$ is $s(t) = 1 + 2t + 3t^2 + \dots + (v-1)t^v + vt^{v+1}$ where $v$ is even. Thus the Castelnuovo diagram is triangular and ends with a white column.
\end{case}
\begin{case}
$(n_e-l,n_o-l) = ((k+1)^2,k(k+1))$ for $k \in \NN$. Then the Castelnuovo polynomial of an object in $\Hilb_{(n_e-l,n_o-l)}(X)$ is $s(t) = 1 + 2t + 3t^2 + \dots + (v-1)t^v + vt^{v+1}$ where $v$ is odd. The Castelnuovo diagram is triangular and ends with a black column.
\end{case}

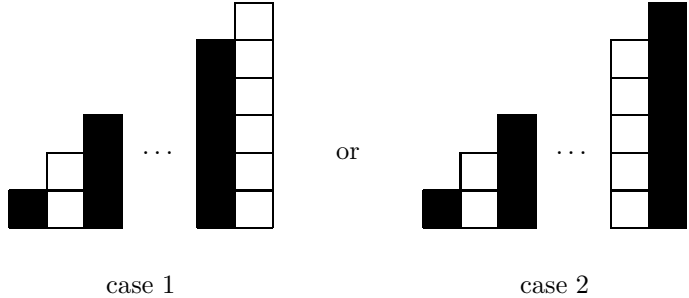
\unitlength 1mm
\begin{picture}(110.00,40.00)(0,0)

\linethickness{0.15mm}
\put(20.00,10.00){\line(1,0){5.00}}
\put(20.00,10.00){\line(0,1){5.00}}
\put(25.00,10.00){\line(0,1){5.00}}
\put(20.00,15.00){\line(1,0){5.00}}

\linethickness{0.15mm}
\put(25.00,15.00){\line(1,0){5.00}}
\put(25.00,15.00){\line(0,1){5.00}}
\put(30.00,15.00){\line(0,1){5.00}}
\put(25.00,20.00){\line(1,0){5.00}}

\linethickness{0.15mm}
\put(25.00,10.00){\line(1,0){5.00}}
\put(25.00,10.00){\line(0,1){5.00}}
\put(30.00,10.00){\line(0,1){5.00}}
\put(25.00,15.00){\line(1,0){5.00}}

\linethickness{0.15mm}
\put(30.00,20.00){\line(1,0){5.00}}
\put(30.00,20.00){\line(0,1){5.00}}
\put(35.00,20.00){\line(0,1){5.00}}
\put(30.00,25.00){\line(1,0){5.00}}

\linethickness{0.15mm}
\put(30.00,15.00){\line(1,0){5.00}}
\put(30.00,15.00){\line(0,1){5.00}}
\put(35.00,15.00){\line(0,1){5.00}}
\put(30.00,20.00){\line(1,0){5.00}}

\linethickness{0.15mm}
\put(30.00,10.00){\line(1,0){5.00}}
\put(30.00,10.00){\line(0,1){5.00}}
\put(35.00,10.00){\line(0,1){5.00}}
\put(30.00,15.00){\line(1,0){5.00}}

\linethickness{0.15mm}
\put(45.00,30.00){\line(1,0){5.00}}
\put(45.00,30.00){\line(0,1){5.00}}
\put(50.00,30.00){\line(0,1){5.00}}
\put(45.00,35.00){\line(1,0){5.00}}

\linethickness{0.15mm}
\put(45.00,25.00){\line(1,0){5.00}}
\put(45.00,25.00){\line(0,1){5.00}}
\put(50.00,25.00){\line(0,1){5.00}}
\put(45.00,30.00){\line(1,0){5.00}}

\linethickness{0.15mm}
\put(45.00,20.00){\line(1,0){5.00}}
\put(45.00,20.00){\line(0,1){5.00}}
\put(50.00,20.00){\line(0,1){5.00}}
\put(45.00,25.00){\line(1,0){5.00}}

\linethickness{0.15mm}
\put(45.00,15.00){\line(1,0){5.00}}
\put(45.00,15.00){\line(0,1){5.00}}
\put(50.00,15.00){\line(0,1){5.00}}
\put(45.00,20.00){\line(1,0){5.00}}

\linethickness{0.15mm}
\put(45.00,10.00){\line(1,0){5.00}}
\put(45.00,10.00){\line(0,1){5.00}}
\put(50.00,10.00){\line(0,1){5.00}}
\put(45.00,15.00){\line(1,0){5.00}}

\linethickness{0.15mm}
\put(50.00,35.00){\line(1,0){5.00}}
\put(50.00,35.00){\line(0,1){5.00}}
\put(55.00,35.00){\line(0,1){5.00}}
\put(50.00,40.00){\line(1,0){5.00}}

\linethickness{0.15mm}
\put(50.00,30.00){\line(1,0){5.00}}
\put(50.00,30.00){\line(0,1){5.00}}
\put(55.00,30.00){\line(0,1){5.00}}
\put(50.00,35.00){\line(1,0){5.00}}

\linethickness{0.15mm}
\put(50.00,25.00){\line(1,0){5.00}}
\put(50.00,25.00){\line(0,1){5.00}}
\put(55.00,25.00){\line(0,1){5.00}}
\put(50.00,30.00){\line(1,0){5.00}}

\linethickness{0.15mm}
\put(50.00,20.00){\line(1,0){5.00}}
\put(50.00,20.00){\line(0,1){5.00}}
\put(55.00,20.00){\line(0,1){5.00}}
\put(50.00,25.00){\line(1,0){5.00}}

\linethickness{0.15mm}
\put(50.00,15.00){\line(1,0){5.00}}
\put(50.00,15.00){\line(0,1){5.00}}
\put(55.00,15.00){\line(0,1){5.00}}
\put(50.00,20.00){\line(1,0){5.00}}

\linethickness{0.15mm}
\put(50.00,10.00){\line(1,0){5.00}}
\put(50.00,10.00){\line(0,1){5.00}}
\put(55.00,10.00){\line(0,1){5.00}}
\put(50.00,15.00){\line(1,0){5.00}}

\put(40.00,20.00){\makebox(0,0)[cc]{$\dots$}}

\linethickness{0.15mm}
\put(75.00,10.00){\line(1,0){5.00}}
\put(75.00,10.00){\line(0,1){5.00}}
\put(80.00,10.00){\line(0,1){5.00}}
\put(75.00,15.00){\line(1,0){5.00}}

\linethickness{0.15mm}
\put(80.00,15.00){\line(1,0){5.00}}
\put(80.00,15.00){\line(0,1){5.00}}
\put(85.00,15.00){\line(0,1){5.00}}
\put(80.00,20.00){\line(1,0){5.00}}

\linethickness{0.15mm}
\put(80.00,10.00){\line(1,0){5.00}}
\put(80.00,10.00){\line(0,1){5.00}}
\put(85.00,10.00){\line(0,1){5.00}}
\put(80.00,15.00){\line(1,0){5.00}}

\linethickness{0.15mm}
\put(85.00,20.00){\line(1,0){5.00}}
\put(85.00,20.00){\line(0,1){5.00}}
\put(90.00,20.00){\line(0,1){5.00}}
\put(85.00,25.00){\line(1,0){5.00}}

\linethickness{0.15mm}
\put(85.00,15.00){\line(1,0){5.00}}
\put(85.00,15.00){\line(0,1){5.00}}
\put(90.00,15.00){\line(0,1){5.00}}
\put(85.00,20.00){\line(1,0){5.00}}

\linethickness{0.15mm}
\put(85.00,10.00){\line(1,0){5.00}}
\put(85.00,10.00){\line(0,1){5.00}}
\put(90.00,10.00){\line(0,1){5.00}}
\put(85.00,15.00){\line(1,0){5.00}}

\linethickness{0.15mm}
\put(100.00,30.00){\line(1,0){5.00}}
\put(100.00,30.00){\line(0,1){5.00}}
\put(105.00,30.00){\line(0,1){5.00}}
\put(100.00,35.00){\line(1,0){5.00}}

\linethickness{0.15mm}
\put(100.00,25.00){\line(1,0){5.00}}
\put(100.00,25.00){\line(0,1){5.00}}
\put(105.00,25.00){\line(0,1){5.00}}
\put(100.00,30.00){\line(1,0){5.00}}

\linethickness{0.15mm}
\put(100.00,20.00){\line(1,0){5.00}}
\put(100.00,20.00){\line(0,1){5.00}}
\put(105.00,20.00){\line(0,1){5.00}}
\put(100.00,25.00){\line(1,0){5.00}}

\linethickness{0.15mm}
\put(100.00,15.00){\line(1,0){5.00}}
\put(100.00,15.00){\line(0,1){5.00}}
\put(105.00,15.00){\line(0,1){5.00}}
\put(100.00,20.00){\line(1,0){5.00}}

\linethickness{0.15mm}
\put(100.00,10.00){\line(1,0){5.00}}
\put(100.00,10.00){\line(0,1){5.00}}
\put(105.00,10.00){\line(0,1){5.00}}
\put(100.00,15.00){\line(1,0){5.00}}

\linethickness{0.15mm}
\put(105.00,35.00){\line(1,0){5.00}}
\put(105.00,35.00){\line(0,1){5.00}}
\put(110.00,35.00){\line(0,1){5.00}}
\put(105.00,40.00){\line(1,0){5.00}}

\linethickness{0.15mm}
\put(105.00,30.00){\line(1,0){5.00}}
\put(105.00,30.00){\line(0,1){5.00}}
\put(110.00,30.00){\line(0,1){5.00}}
\put(105.00,35.00){\line(1,0){5.00}}

\linethickness{0.15mm}
\put(105.00,25.00){\line(1,0){5.00}}
\put(105.00,25.00){\line(0,1){5.00}}
\put(110.00,25.00){\line(0,1){5.00}}
\put(105.00,30.00){\line(1,0){5.00}}

\linethickness{0.15mm}
\put(105.00,20.00){\line(1,0){5.00}}
\put(105.00,20.00){\line(0,1){5.00}}
\put(110.00,20.00){\line(0,1){5.00}}
\put(105.00,25.00){\line(1,0){5.00}}

\linethickness{0.15mm}
\put(105.00,15.00){\line(1,0){5.00}}
\put(105.00,15.00){\line(0,1){5.00}}
\put(110.00,15.00){\line(0,1){5.00}}
\put(105.00,20.00){\line(1,0){5.00}}

\linethickness{0.15mm}
\put(105.00,10.00){\line(1,0){5.00}}
\put(105.00,10.00){\line(0,1){5.00}}
\put(110.00,10.00){\line(0,1){5.00}}
\put(105.00,15.00){\line(1,0){5.00}}

\put(95.00,20.00){\makebox(0,0)[cc]{$\dots$}}

\linethickness{0.15mm}
\put(20.00,10.00){\rule{5.00\unitlength}{5.00\unitlength}}

\linethickness{0.15mm}
\put(30.00,10.00){\rule{5.00\unitlength}{15.00\unitlength}}

\linethickness{0.15mm}
\put(45.00,10.00){\rule{5.00\unitlength}{25.00\unitlength}}

\linethickness{0.15mm}
\put(75.00,10.00){\rule{5.00\unitlength}{5.00\unitlength}}

\linethickness{0.15mm}
\put(85.00,10.00){\rule{5.00\unitlength}{15.00\unitlength}}

\linethickness{0.15mm}
\put(105.00,10.00){\rule{5.00\unitlength}{30.00\unitlength}}

\put(65.00,20.00){\makebox(0,0)[cc]{or}}

\put(41.25,2.50){\makebox(0,0)[cc]{}}

\put(37.50,2.50){\makebox(0,0)[cc]{case 1}}

\put(92.50,2.50){\makebox(0,0)[cc]{case 2}}

\end{picture} \\
The next proposition shows that not only the Hilbert series but also the Betti numbers of an object in $\Hilb_{(n_e-l,n_o-l)}(X)$ are fully determined.
\begin{proposition} \label{minimalresunique}
Let $(n_e,n_o) \in N$ and let $l \geq 0$ be as in \eqref{uniquel}. 
Let $I_0 \in \Hilb_{(n_e-l,n_o-l)}(X)$. Then $I_0$ has a minimal resolution of the form
\[
0 \r A(-c-1)^{c} \r A(-c)^{c+1} \r I_0 \r 0  
\]
where 
\begin{eqnarray*} 
c =
\left\{
\begin{array}{ll} 
2k & \text{if $(n_e-l,n_o-l)=(k^2,k(k+1))$} \\
2k+1 & \text{if $(n_e-l,n_o-l)=((k+1)^2,k(k+1))$}
\end{array}
\right.
\end{eqnarray*}
\end{proposition}
\begin{proof}
By Proposition \ref{prophilbert} and same arguments as in the proof of Theorem \ref{theorem2}. 
\end{proof}
\begin{remark} \label{minimalresuniqueremark}
In the notations of the previous proposition one may compute \linebreak $\dim_k \Ext^{1}_A(I_0,I_0) = 0$ which indicates that up to isomorphism $\Hilb_{(n_e-l,n_o-l)}(X) = R_{(n_e-l,n_o-l)}(A)$ consist of only one object. See also \S\ref{Line bundles on $X$} below for linear $A$ and the proof of Theorem \ref{theorem4} in Section \ref{filtrationscubic} for generic elliptic $A$.
\end{remark}

\section{Ideals of linear cubic Artin-Schelter regular algebras} \label{Ideals of linear cubic Artin-Schelter regular algebras}

In this section we let $A$ be a linear cubic AS-algebra. As $\Tails(A)$ is equivalent to $\Qcoh(\PP^1 \times \PP^1)$ line bundles on $X = \Proj A$ are determined by line bundles on $\PP^1 \times \PP^1$. We will briefly recall the description of these objects which will lead to a characterisation of the set $R(A)$ of reflexive rank one modules over $A$, see Proposition \ref{commfreecubic}. We will end with a discussion on the Hilbert scheme of points.

Let $Y = \PP^1 \times \PP^1$ denote the quadric surface. Consider for any integers $m,n$ the canonical line bundle $\Oscr_Y(m,n) = \Oscr_{\PP^1}(m) \boxtimes \Oscr_{\PP^1}(n)$. It is well-known that the map $\Pic(Y) \r \ZZ \oplus \ZZ: \Oscr_Y(m,n) \mapsto (m,n)$ is a group isomorphism i.e. the objects $\Oscr_Y(m,n)$ are the only reflexive rank one sheaves on $\PP^1 \times \PP^1$. Note there are short exact sequences on $\coh(Y)$
\begin{equation} \label{sesinherited}
\begin{aligned}
& 0 \r \Oscr_Y(m,n-1) \r \Oscr_Y(m,n)^2 \r \Oscr_Y(m,n+1) \r 0 \\
& 0 \r \Oscr_Y(m-1,n) \r \Oscr_Y(m,n)^2 \r \Oscr_Y(m+1,n) \r 0 
\end{aligned}
\end{equation} 
for all integers $m,n$.

\subsection{Line bundles} \label{Line bundles on $X$}

As usual we put $X = \Proj A$ and $\Oscr_X = \Oscr$. In \cite{ncquadrics} it is shown there is an equivalence of categories $\Qcoh(Y) \cong \Qcoh(X)$ such that $\Oscr_{Y}(k,k)$ corresponds to $\Oscr_{X}(2k)$ and $\Oscr_{Y}(k,k+1)$ corresponds to $\Oscr_{X}(2k+1)$. See also \cite[\S 11.3]{VdBSt}.
Further, for any integers $m,n$ we denote the image of $\Oscr_Y(m,n)$ under the equivalence $\Qcoh(Y) \cong \Qcoh(X)$ as $\Oscr(m,n)$.
Clearly these objects $\Oscr(m,n) \in \coh(X)$ are the only line bundles on $X$. 

From \eqref{sesinherited} we compute the class of $\Oscr(m,n)$ in $K_0(X)$ 
\[
[\Oscr(m,n)] = [\Oscr] + (m-n)[\Sscr] + n[\Qscr] + n(m+1)[\Pscr]
\]
for all $m,n \in \ZZ$. Using \eqref{shiftcubic} we obtain
\begin{align*}
\Oscr(m,n)(2k) = \Oscr(m+k,n+k), \quad \Oscr(m,n)(2k+1) = \Oscr(n+k,m+k+1)   
\end{align*}
for all $m,n,k \in \ZZ$. By \eqref{shiftcubic} it is easy to see $\Oscr(m,n)(-m-n) = \Oscr(u,-u)$ is a normalized line bundle where
\begin{eqnarray*}
u =
\left\{
\begin{array}{ll} 
(m-n)/2 & \text{if $m-n$ is even} \\
(n-m-1)/2 & \text{if $m-n$ is odd}
\end{array}
\right.
\end{eqnarray*}
Since $[\Oscr(u,-u)] = [\Oscr] + 2u[\Sscr] - u[\Qscr] - u(u+1)[\Pscr]$ the invariants $(n_e,n_o)$ of $\Oscr(u,-u)$ are given by $(n_e,n_o) = (u^2,u(u+1))$. Either $k = u \geq 0$ or $k = -u-1 \geq 0$. These two possibilities correspond to Cases 1 and 2 of Section \ref{Hilbert series of torsion free rank one modules}. In particular $\Rscr_{(n_e,n_o)}(X)$ is nonempty if and only if $(n_e,n_o)$ is $(k^2,k(k+1))$ or $((k+1)^2,k(k+1))$ for some integer $k \geq 0$ and in that case $\Rscr_{(n_e,n_o)}(X) =$ \linebreak $\{ \Oscr(n_o - n_e,n_e - n_o) \}$.

Proposition \ref{minimalresunique} implies that a minimal resolution for $\Oscr(m,n)$ is of the form
\begin{align*} 
0 \r \Oscr(2n-1)^{m-n} \r \Oscr(2n)^{m-n+1} \r \Oscr(m,n) \r 0 & \quad \text{ if } m \geq n, \\
0 \r \Oscr(2m)^{n-m-1} \r \Oscr(2m+1)^{n-m} \r \Oscr(m,n) \r 0 & \quad \text{ if } m < n.
\end{align*}
We have shown
\begin{proposition} \label{commfreecubic}
Assume $A$ is linear and let $I \in \grmod(A)$ be a reflexive graded right ideal of $A$. Then $I$ has a minimal resolution of the form
\begin{equation} \label{smallres}
0 \r A(-c-1)^{c} \r A(-c)^{c+1} \r I(d) \r 0
\end{equation} 
for some integers $d$ and $c$. 
As a consequence $R_{(n_e,n_o)}(A) = \emptyset = \Rscr_{(n_e,n_0)}(X)$ unless 
$n_e = (n_e - n_o)^{2}$ i.e. $(n_e,n_o) = ((k+1)^2,k(k+1))$ or $(n_e,n_o) = (k^2,k(k+1))$ for some $k \in \NN$.
\end{proposition}

\subsection{Hilbert scheme of points} \label{Hilbert scheme of points for the quadric surface}

The Hilbert scheme of points for $Y = \PP^1 \times \PP^1$, which we will denote by $\Hilb(Y)$, parameterizes the torsion free rank one sheaves on $Y$ up to shifting. By the category equivalence $\Qcoh(Y) \cong \Qcoh(X)$ where $X = \Proj A$ we see $\Hilb(Y)$ also parameterizes the torsion free rank one objects on $X$ up to shifting. Let $\Iscr \in \coh(X)$ be such an object. Put $\Iscr = \pi I$ where $I \in \grmod(A)$. Thus $\Iscr^{\ast \ast} := \pi I^{\ast \ast}$ is a line bundle on $X$ of rank one hence $\Iscr^{\ast \ast} \cong \Oscr(m,n)$ for some integers $m,n$. By \cite[Corollary 4.2]{ATV2} there is an exact sequence
\[
0 \r \Iscr \r \Iscr^{\ast \ast} \r \Nscr \r 0
\]
where $\Nscr \in \coh(X)$ is a zero dimensional object of some degree $l \geq 0$. Since $\Nscr$ admits a filtration by point objects on $X$ we have $[\Nscr] = l[\Pscr]$. Also $\Iscr^{\ast \ast}(d) \cong O(u,-u)$ for some $d,u \in \ZZ$. Computing the class of $\Iscr(d)$ in $K_0(X)$ we find
\[
[\Iscr(d)] = [\Oscr] + 2u[\Sscr] - u[\Qscr] - \left( u(u+1) - l \right)[\Pscr]
\]
from which we deduce $\Iscr(d) \in \Hilb_{(n_e,n_o)}(X)$, as defined in \S\ref{Normalized line bundles}, where $(n_e,n_o) = (u^2 + l, u(u+1) + l)$. Again we separate
\setcounter{case}{0}
\begin{case}
$u \geq 0$. Put $k = u$. Then $(n_e,n_o) = (k^2+l, k(k+1)+l)$ where $k,l \in \NN$. 
\end{case}
\begin{case}
$u <0$. Put $k = -u-1$. Then $(n_e,n_o) = ((k+1)^2+l, k(k+1)+l)$ where $k,l \in \NN$. 
\end{case}
\begin{remark}
By the above discussion we may associate invariants $(n_e,n_o) \in N = \{ (n_e,n_o) \in \NN^{2} \mid n_e - (n_e - n_o)^2 \geq 0 \}$ to any object in $\Hilb(Y)$. Let $\Hilb_{(n_e,n_o)}(Y)$ denote the associated parameter space.  
The dimension of $\Hilb_{(n_e,n_o)}(Y)$ may be deduced as follows. Given $(n_e,n_o) \in N$ fixes $l \in \NN$ and $u \in \ZZ$ as above. The number of parameters to choose $\Oscr(u,-u)$ is zero. On the other hand, to choose a point in $\PP^1 \times \PP^1$ we have two parameters. Thus to pick a zero-dimensional subsheaf $\Nscr$ of degree $l$ we have $2l$ parameters since such $\Nscr$ admits a filtration of length $l$ in points of $\PP^1 \times \PP^1$. Hence the freedom of choice in a normalized torsion free rank one sheaf $\Iscr$ is $2l$. Hence $\dim \Hilb_{(n_e,n_o)}(Y) = 2l$. Since $l = n_e - (n_e - n_o)^{2}$ we have $\dim \Hilb_{(n_e,n_o)}(Y) = 2\left( n_e - (n_e - n_o)^{2} \right)$. 
\end{remark}

\section{Some results on line and conic objects}

In this section we gather some additional results on line objects and conic objects on quantum quadrics which will be used later on. These results are obtained by using similar techniques as in \cite{Ajitabh,ATV2}. 

Let $A$ be a cubic AS-algebra. We use the notations of \S\ref{Geometric data associated to a cubic AS-algebra}. In particular $(E,\sigma,\Oscr_{E}(1))$, $B = B(E,\sigma,\Oscr_{E}(1))$, $(C,\sigma,\Oscr_{C}(1))$, $D = B(C,\sigma,\Oscr_{C}(1)) = \Gamma_{\ast}(\Oscr_C)$, $g$ and $h$ will have their usual meaning. Recall the isomorphism of $k$-algebras $A/hA \xrightarrow{\cong} D:a \mapsto \overline{a}$. The dimension of objects in $\grmod(B)$, $\grmod(D)$ or $\tails(B)$, $\tails(D)$ will be computed in $\grmod(A)$ or $\tails(A)$.  We begin with
\begin{lemma} \label{incidence}
Let $w \in A_{d}$ for some integer $d \geq 1$ and put $W = A/wA$, $\Wscr = \pi W$. 
\begin{enumerate}
\item
Let $p \in C$. Then $\Hom_{X}(\Wscr,\Nscr_p) \neq 0$ if and only if $\overline{w}(p) = 0$. 
\item
$\dim_k \Hom_{X}(\Wscr,\Nscr_p) \leq 1$ for all $p \in C$.
\end{enumerate}
\end{lemma}
\begin{proof}
Firstly, if $f: W \r N_p$ is non-zero then $\pi f: \Wscr \r \Nscr_p$ is non-zero since $N_p$ is socle-free i.e. $\underline{\Hom}_A(k,N_p) = 0$. Conversely, $\Hom_{X}(\Wscr,\Nscr_p) \neq 0$ implies $\Hom_{A}(W,N_p) \neq 0$. Indeed, a non-zero map $g:\Wscr \r \Nscr_p$ yields a surjective map $\omega g: W \r (\omega \Nscr_p)_{\geq n}$ for $n \gg 0$. Now $(\omega \Nscr_p)_{\geq n} = N_{\sigma^np}(-n) \subset N_{p}$, which 
yields a non-zero map $W \r N_p$. 

So to prove the first statement it is sufficient to show $\Hom_{A}(W,N_p) \neq 0$ if and only if $\overline{w}(p) = 0$. This is proved in a similar way as \cite{Ajitabh}. For convenience we shortly repeat the arguments. Writing down resolutions for $W$, $N_p$ we see there is a non-zero map $f:W \r N_p$ if and only if we may find (non-zero) maps $f_0$, $f_1$ making the following diagram commutative
\begin{equation*}
\begin{diagram}[heads=LaTeX]
&& 0 & \rTo & A(-d) & \rTo^{w \cdot \quad \quad} & A & \rTo & W & \rTo 0 \\
&& && \dTo^{f_1} && \dTo^{f_0} && && && \\
0 & \rTo & A(-3) & \rTo & A(-1) \oplus A(-2) & \rTo & A & \rTo^{\theta} & N_p & \rTo 0 
\end{diagram}
\end{equation*}
The resolutions being projective, this is equivalent with saying there is a non-zero map $f_0$ such that $\theta \circ f_0 \circ w = 0$, i.e. $\overline{w}(p) = 0$. 

The second part is shown by applying $\Hom_{X}(-,\Nscr_p)$ to the short exact sequence $0 \r \Oscr(-d) \r \Oscr \r \Wscr \r 0$ and bearing in mind $\Hom_{X}(\Oscr,\Nscr_p) = k$.
\end{proof}
\begin{remark}
It follows from the first part of the previous lemma there exists at least one $p \in C$ for which $\Hom_{X}(\Wscr,\Nscr_p) \neq 0$. Moreover any such non-zero morphism is surjective since point objects are simple objects in $\coh(X)$.
\end{remark}

\subsection{Line objects}

Let $u = \lambda x + \mu y \in A_{1}$. Then $\overline{u} \in D_{1} = H^{0}(C,\Oscr_C(1))$ and a point $p = (p_1,p_2) \in C$ vanishes at $\overline{u}$ i.e. $\overline{u}(p) = 0$ if and only if $p_1 = (-\mu:\lambda) \in \PP^1$. 
We have shown
\begin{lemma} \label{uniqueline}
Let $p \in C$. There exists, up to isomorphism, a unique line object $\Sscr$ on $X$ for which $\Hom_{X}(\Sscr,\Nscr_p) \neq 0$.
\end{lemma}
In case $A$ is elliptic then $E$ is a divisor of bidegree $(2,2)$ which means that for $u \in A_1$ the line $\{ u = 0 \} \times \PP^{1}$ meets $C$ in at most two different points $p,q$. 

For general $A$ we call two different points $p,q \in C$ {\em collinear} if $\overline{l}(p) = \overline{l}(q) = 0$ for some global section in $\overline{l} \in H^{0}(C,\Oscr_C(1)) = D_{1}$. It follows from the previous discussion that $\pr_{1}p = \pr_{1}q$.

\subsection{Conic objects}
We now deduce
\begin{lemma} \label{threepoints}
Let $p,q,r$ be three distinct points in $C$. There exists, up to isomorphism, a unique conic object $\Qscr$ on $X$ for which $\Hom_{X}(\Qscr,\Nscr_p) \neq 0$, $\Hom_{X}(\Qscr,\Nscr_q) \neq 0$ and $\Hom_{X}(\Qscr,\Nscr_r) \neq 0$.
\end{lemma}
\begin{proof}
Due to Lemma \ref{incidence} it will be sufficient to prove there exists, up to scalar multiplication, a unique quadratic form $v \in A_2$ for which $\overline{v}(p) = \overline{v}(q) = \overline{v}(r) = 0$.

Writing $v = \lambda_1 x^2 + \lambda_2 xy + \lambda_3 yx + \lambda_4 y^2$ where $\lambda_i \in k$ and $p = ((\alpha:\beta),(\alpha':\beta')) \in C \subset \PP^1 \times \PP^1$, we see $\overline{v}(p) = 0$ if and only if $\lambda_1 \alpha \alpha' + \lambda_2 \alpha \beta' + \lambda_3 \beta \alpha' + \lambda_4 \beta \beta' = 0$. The condition $\overline{v}(p) = \overline{v}(q) = \overline{v}(r) = 0$ then translates to a system of three linear equations in $\lambda_1,\dots,\lambda_4$, which admits a non-trivial solution. Moreover, this solution is unique (up to scalar multiplication) unless all maximal minors are zero, which implies that at least two points of $p,q,r$ coincide.
\end{proof}
Subobjects of line objects on $X$ are shifted line objects \cite{ATV2}. We may prove a similar result for conic objects.
\begin{lemma} \label{subobjectscritical}
Let $\Qscr$ be a conic object and $p \in C$. Assume $\Hom_{X}(\Qscr,\Nscr_p) \neq 0$. 
\begin{enumerate}
\item
The kernel of a non-zero map $\Qscr \r \Nscr_p$ is a shifted conic object $\Qscr'(-1)$.
\item
Assume $A$ is elliptic and $\sigma$ has infinite order. If in addition $\Qscr$ is critical then all subobjects $\Qscr$ are shifted critical conic objects.
\end{enumerate}
\end{lemma}
\begin{proof}
Firstly, let $f$ denote such a non-zero map $\Qscr \ra \Nscr_p$. Since $\Nscr_p$ is simple, $f$ is surjective. Putting $\Qscr = \pi Q$ where $Q$ is a conic module over $A$ it is sufficient to show that the kernel of a surjective map $Q \ra (N_{p})_{\geq n}$ is of the form $Q'(-1)$ for some conic object $Q'$. This is done by taking the cone of the induced map between resolutions of $Q$ and $(N_{p})_{\geq n}$.

Secondly, as $\Qscr$ is critical, any quotient of $\Qscr$ has dimension zero and since $\sigma$ has infinite order such a quotient admits a filtration by shifted point objects on $X$, see \cite{ATV2}. By the first part this completes the proof.
\end{proof}
We will also need the dual statement of the previous result.
\begin{lemma} \label{extentionsconic}
Let $\Qscr$ be a conic object and $p \in C$. Assume $\Ext^{1}_{X}(\Nscr_p,\Qscr) \neq 0$.
\begin{enumerate}
\item
The middle term of a non-zero extension in $\Ext^{1}_{X}(\Nscr_p,\Qscr)$ is a shifted conic object $\Qscr'(1)$.
\item
Assume $A$ is elliptic and $\sigma$ has infinite order. Then any extension of $\Qscr$ by a zero dimensional object is a shifted conic object.
\end{enumerate}
\end{lemma}
\begin{proof}
Again the second statement is clear from the first one thus it suffices to prove the first part. Put $\Qscr = \pi Q$ where $Q$ is a conic module over $A$. Let $\Jscr$ denote the middle term of a non-trivial extension i.e. $0 \ra \Qscr \ra \Jscr \ra \Nscr_p \ra 0$. It is easy to see $\Jscr$ is pure and $\omega \Jscr \in \coh(X)$ has projective dimension one, see for example (the proof of) \cite[Proposition 3.4.1]{DV2}. Put $J = \omega \Jscr$. Application of $\omega$ gives a short exact sequence
\begin{equation} \label{sescriticals}
0 \ra Q \ra J \ra (N_{p})_{\geq n} \ra 0.
\end{equation}
Applying $\underline{\Hom}_{A}(-,A)$ on \eqref{sescriticals} yields $0 \ra J^{\vee} \ra Q^{\vee} \ra ((N_{p})_{\geq n})^{\vee} \ra 0$. As $((N_{p})_{\geq n})^{\vee}$ is a shifted point module and $Q^{\vee}$ is a shifted conic module it follows from Lemma \ref{subobjectscritical} that $J^{\vee}$ is also a shifted conic module. Hence the same is true for $J^{\vee \vee}$. Consideration of Hilbert series shows $J^{\vee \vee} = Q'(1)$ for some conic module $Q'$ over $A$. Since $\omega \Jscr$ is Cohen-Macaulay, \cite[Corollary 4.2]{ATV2} implies $\pi J^{\vee \vee} = \pi J = \Jscr$. This finishes the proof.
\end{proof}
\begin{remark}
Lemmas \ref{subobjectscritical} and \ref{extentionsconic} are in contrast with the situation for quadratic AS-algebras \cite[\S 4]{Ajitabh} where a non-zero map $A/vA \ra N_p$ (where $v \in A_2$ and $p \in C$) will yield an exact sequence $0 \ra Q'(-1) \ra A/vA \ra N_p \ra 0$ for which $Q'$ has a resolution of the form $0 \ra A(-1)^2 \ra A^2 \ra Q' \ra 0$.
\end{remark}
Let $\Zscr$ denote the full subcategory of $\coh(X)$ whose objects consist of zero dimensional objects of $\coh(X)$. $\Zscr$ is a Serre subcategory of $\coh(X)$, see for example \cite{Weibel}. We say $\Mscr, \Nscr \in \coh(X)$ are {\em equivalent up to zero dimensional objects} if their images in the quotient category $\coh(X)/\Zscr$ are isomorphic. We say $\Mscr$ and $\Nscr$ are {\em different modulo zero dimensional objects} if they are not equivalent up to zero dimensional objects. Using Lemmas \ref{subobjectscritical} and \ref{extentionsconic} one proves
\begin{lemma} \label{equivalentsubobject}
Assume $A$ is elliptic and $\sigma$ has infinite order. Then two critical conic objects on $X$ are equivalent up to zero dimensional objects if and only if they have a common subobject.
\end{lemma}
We now come to a key result which we will need in \S\ref{ref-5.7cubic} below.
\begin{lemma} \label{infinitecubic}
Assume $k$ is uncountable, $A$ is elliptic and $\sigma$ has infinite order. Let $p,p' \in C$ for which $p,p',\sigma p,\sigma p'$ are pairwise different and non-collinear. Then, modulo zero dimensional objects, there exist infinitely many critical conic objects $\Qscr$ for which $\Hom_{X}(\Qscr,\Nscr_p) \neq 0$ and  $\Hom_{X}(\Qscr,\Nscr_{p'}) \neq 0$. 
\end{lemma}
\begin{proof}
Write $p = ((\alpha_0:\beta_0),(\alpha_1:\beta_1)) \in C$. We prove the lemma in seven steps.
\setcounter{step}{0}
\begin{step}
Let $d \in \NN$ and let $\Qscr, \Qscr'$ be two critical objects for which $\Qscr'(-d) \subset \Qscr$. Then there is a filtration $\Qscr'(-d) = \Qscr_{d}(-d) \subset \Qscr_{d-1}(-d+1) \subset \dots \subset \Qscr_{1}(-1) \subset \Qscr_{0} = \Qscr$ where the $\Qscr_i$ are critical conic objects and the successive quotients are point objects on $X$. This follows from the proof of Lemma \ref{subobjectscritical}. 
\end{step}
\begin{step}
Up to isomorphism there are uncountably many conic objects $\Qscr$ on $X$ for which $\Hom_{X}(\Qscr,\Nscr_p) \neq 0$, $\Hom_{X}(\Qscr,\Nscr_{p'}) \neq 0$. See the proof of Lemma \ref{threepoints}.
\end{step}
\begin{step}
Let $\Ascr$ denote the set of isoclasses of critical conic objects $\Qscr$ for which $\Hom_{X}(\Qscr,\Nscr_p) \neq 0$, $\Hom_{X}(\Qscr,\Nscr_{p'}) \neq 0$. Then $\Ascr$ is an uncountable set. By the previous step it is sufficient to show there are only finitely many non-critical conic objects $\Qscr$ on $X$ for which $\Hom_{X}(\Qscr,\Nscr_p) \neq 0$, $\Hom_{X}(\Qscr,\Nscr_{p'}) \neq 0$. For such an object $\Qscr$ it is easy to see there exists an exact sequence 
\begin{equation} \label{sesconic}
0 \r \Sscr'(-1) \r \Qscr \r \Sscr \r 0
\end{equation}
for some line objects $\Sscr,\Sscr'$ on $X$. Also, $\dim_{k}\Ext^{1}_{X}(\Sscr,\Sscr'(-1)) \leq 1$ hence $\Qscr$ is, up to isomorphism, fully determined by $\Sscr$ and $\Sscr'$. We deduce from \eqref{sesconic} that either $\Hom_{X}(\Sscr,\Nscr_p) \neq 0$, $\Hom_{X}(\Sscr', \Nscr_{\sigma p}) \neq 0$ or $\Hom_{X}(\Sscr,\Nscr_{p'}) \neq 0$, $\Hom_{X}(\Sscr', \Nscr_{\sigma p'}) \neq 0$. By Lemma \ref{uniqueline} this means there are at most two non-critical conic objects for which $\Hom_{X}(\Qscr,\Nscr_p) \neq 0$ and $\Hom_{X}(\Qscr,\Nscr_{p'}) \neq 0$.
\end{step}
\begin{step}
Let $\Bscr \subset \Ascr$ denote the set of conic objects $\Qscr=\pi Q$ for which $Q$ is $h$-torsion free. 
Then $\Bscr$ is uncountable. Indeed, writing $Q = A/vA$ we find $Q$ is $h$-torsion free (meaning multiplication by $h$ is injective) unless $\overline{v}:\Oscr_{C}(-2) \r \Oscr_{C}$ is not injective. Hence $\overline{v} = 0$ which means that $v$ and $h$ have a common divisor. As $v$ is not a product of linear forms, $v$ divides $h$. Up to scalar multiplication there are only finitely many possibilities for such $v$.
\end{step}
\begin{step}
For any $\Qscr \in \Bscr$ there are, up to isomorphism, only finitely many points objects $\Nscr_p$ for which $\Hom_{X}(\Qscr,\Nscr_p) \neq 0$ or $\Ext^{1}_{X}(\Nscr_p,\Qscr(-1)) \neq 0$. To show this, write $\Qscr = \pi(A/vA)$ and $Q = A/vA$. Since $v$ does not divide $h$, it does not divide $g$ thus $Q$ is also $g$-torsion free. Thus $Q/gQ$ is a $B$-module of GK-dimension one so $(Q/gQ){\,\widetilde{}\,\,}$ is a finite dimensional $\Oscr_E$-module. Writing $\overline{v}_g$ for the image of $v$ in $B$ this implies there are only finitely many points $p \in E$ such that $\overline{v}_g(p) = 0$. By the same methods used in the proof of Lemma \ref{incidence} one may show there are finitely many point objects $\Nscr_p$ on $X$ for which $\Hom_{X}(\Qscr,\Nscr_p) \neq 0$.

For the second part, Serre duality implies $\Ext^{i}_{X}(\Nscr_p,\Qscr(-1)) \cong \Ext^{2-i}_{X}(\Qscr,\hat{\Nscr_p})'$ for $i = 0,1,2$ and a suitable point object $\hat{\Nscr_p}$ on $X$. By $\chi(\Qscr,\hat{\Nscr_p}) = 0$, Lemma \ref{incidence}(2) and the first part of Step 5 we are done.
\end{step}
\begin{step}
For any $\Qscr_i \in \Bscr$ and any integer $d \geq 0$ the following subset of $\Bscr$ is finite 
\[
\Vscr_d(\Qscr_i) = \{ \Qscr \in \Bscr \mid \Qscr'(-d) \subset \Qscr \text{ for a conic object $\Qscr'$ for which } \Qscr'(-d) \subset \Qscr_i \}
\]
We will prove this for $d = 1$, for general $d$ the same arguments may be used combined with Step 1. Let $\Qscr'(-1) \subset \Qscr_i$. Note $\Qscr' \in \Bscr$. Clearly any conic object $\Qscr$ on $X$ for which $\Qscr'(-1) \subset \Qscr$ holds is represented by an element of $\Ext^{1}_{X}(\Nscr_p,\Qscr'(-1))$ for some point object $\Nscr_p$, and two such conic objects $\Qscr$ are isomorphic if and only if the corresponding extensions only differ by a scalar. By Step 5 and its proof there are only finitely many such $\Qscr$, up to isomorphism. 
\end{step}
\begin{step}
There exist infinitely many critical conic objects $\Qscr_0, \Qscr_1, \Qscr_2, \dots$ for which $\Hom_{X}(\Qscr_i,\Nscr_p) \neq 0$, $\Hom_{X}(\Qscr_i,\Nscr_{p'}) \neq 0$ and $\Qscr_i$, $\Qscr_j$ do not have a common subobject for all $j < i$. Indeed, choose $\Qscr_0 \in \Bscr$ arbitrary and having $\Qscr_0, \Qscr_1, \dots, \Qscr_{i-1}$ we pick $\Qscr_i$ as an element of $\Bscr$ which does not appear in the countable subset $\bigcup_{d \in \NN, j<i}\Vscr_{d}(\Qscr_j)$. By Lemma \ref{subobjectscritical} subobjects of critical conic objects are shifted critical conic objects hence Step 7 follows.
\end{step}
Combining Step 7 with Lemma \ref{equivalentsubobject} completes the proof.
\end{proof}
\begin{remark} \label{remarkuncountable}
If $A$ is of generic type A for which $\sigma$ has infinite order then Lemma \ref{infinitecubic} may be proved alternatively by observing that for any conic object $\Qscr = \pi(A/vA)$ containing a shifted conic object $\Qscr' = \pi(A/v'A)$ we have
\[
\div(\overline{v'}) = (\sigma^{a}p) + (\sigma^{b}q) + (\sigma^{c}r) + (\sigma^{-a-b-c}r) \quad \text{ for some } a,b,c \in \ZZ
\]
where we have written $\div(\overline{v}) = (p) + (q) + (r) + (s)$ for the divisor of zeroes of $\overline{v} \in D_2$. This observation is proved by using similar methods as in \cite{Ajitabh}, see also \cite[Theorem 3.2]{AjVdB}.
Thus if $A$ is of generic type A we do not need the hypothesis $k$ is uncountable in Lemma \ref{infinitecubic}.
\end{remark}

\section{Restriction of line bundles to the divisor $C$} \label{Restriction of line bundles to the divisor $C$}

In this section $A$ is an elliptic cubic AS-regular algebra and $X = \Proj A$. We will need the following lemma which extends \cite[Lemma 2.8.2]{DV1} to $C = E_{\red}$.
\begin{lemma} \label{ref-2.8.2-18} 
Assume $A$ is elliptic. Let $M \in \grmod(A)$ be such that $M/Mh \in \tors(A)$. Then $\gkdim M \leq 1$. If $\sigma$ has infinite order then $M \in \tors(A)$. 
\end{lemma} 
\begin{proof} 
Let $c = \deg h$. Multiplication by $h$ induces an isomorphism $M_{n} \cong M_{n+c}$ for large $n$. Hence $\gkdim M \leq 1$. 

For the second part it suffices to prove $Mh = 0$. Assume by contradiction $Mh \neq 0$. Write $T \subset M$ for the submodule of $h$-torsion elements of $M$. Then $M'=M/T$ is a non-zero $h$-torsion free module of GK-dimension $\leq 1$. Write $(M'_h)_0$ for the degree zero part of the localization of $M'$ at the powers of $h$. We find that $(M'_h)_0$ is a finite dimensional representation of $(A_h)_0$. By \cite[Proposition 5.18]{ATV2} it is not difficult to see $(A_h)_0 = (A_g)_0$. If $\sigma$ has infinite order then $(A_g)_0$ is a simple ring \cite{ATV2}. In particular it has no finite dimensional representations. Thus $(M'_h)_0 = 0$ i.e. there is a positive integer $i$ for which $M'h^{i} = 0$. For such a minimal $i$ this implies $0\neq M''=M'h^{i-1} \subset M'$ satisfies $M''h=0$, a contradiction.
\end{proof} 
Define a map of noncommutative schemes $u: C \ra X$ by 
\begin{eqnarray*} 
\begin{array}{ll}
u^\ast \pi M =(M\otimes_A D){\,\widetilde{}\,\,} & \text{ for } M \in \GrMod(A), \\
u_\ast\Mscr =\pi(\Gamma_\ast(\Mscr)_A) & \text{ for } \Mscr \in \Qcoh(C) 
\end{array}
\end{eqnarray*} 
We refer to $u^\ast(\pi M)$ as the \emph{restriction} of $\pi M$ to $C$. Note $u_\ast$ is an exact functor. The following result is proved as in \cite[Propositions 4.3, 4.4]{DV1} using Lemma \ref{ref-2.8.2-18}. See also \cite{NS}.
\begin{proposition} \label{propvb}
Let $A$ be an elliptic cubic AS-regular algebra.
\begin{enumerate} 
\item
If $\mathcal{M}$ is a vector bundle in $X$ then $L_ju^{\ast}\mathcal{M}=0$ for $j>0$ and $u^{\ast}\mathcal{M}$ is a vector bundle on $C$.
\item Assume $\sigma$ has infinite order and let $\mathcal{M}\in D^b(\coh(X))$ for which $Lu^{\ast}\mathcal{M}$ is a vector bundle on $C$. Then $\mathcal{M}$ is a vector bundle on $X$.
\end{enumerate}
\end{proposition}
We will now assume for the rest of Section \ref{Restriction of line bundles to the divisor $C$} $A$ is of generic type A (see Example \ref{typeA}). Recall from Example \ref{ellipticsmooth} $E$ is a smooth divisor of bidegree $(2,2)$ on $\PP^1 \times \PP^1$. Since $E$ has arithmetic genus $1$ it is a smooth elliptic curve \cite{H}. In particular the canonical sheaf $\omega_E$ is isomorphic to $\Oscr_E$. Fixing a group law on $E$ the automorphism $\sigma$ is a translation by some element $\xi \in E$. Thus $\sigma p = p + \xi$. As $E$ is reduced the geometric data $(E,\sigma,\Oscr_{E}(1))$ and $(C,\sigma,\Oscr_{C}(1))$ coincide. 

We write $o$ for the origin of the group law of $C$. For $p,q,r \in C$ we have $p+q=r$ if and only if the divisors $(p) + (q)$ and $(r) + (o)$ on $C$ are linearly equivalent \cite[Chapter IV \S 4]{H}. We will use the notations $\det(\Escr):=\Lambda^{\rank(\Escr)}\Escr$ and $\deg(\Escr):=\deg(\det(\Escr))$ for vector bundles $\Escr \in \coh(C)$. The Riemann-Roch theorem states $\chi(\Oscr_C,\Escr) = \dim_{k}\Hom_{C}(\Oscr_C,\Escr) - \dim_{k}\Ext^{1}_{C}(\Oscr_C,\Escr) = \deg \Escr$ and Serre duality on $C$ is given by $\Ext^{1}_{C}(\Oscr_C,\Escr) \cong \Hom_C(\Escr,\Oscr_C)'$. According to \cite[Ex II. 6.11]{H} we have a group isomorphism
\begin{align*}
\Pic(C) \oplus \ZZ \r K_{0}(C): (\Oscr(D),r) \mapsto r[\Oscr_C] + \sum_i n_i [\Oscr_{p_i}]
\end{align*}
where we have written $D = \sum_i n_i(p_i) \in \Cl(C)$. The projection $K_0(C) \r \ZZ$ is given by the rank and the projection $K_0(C) \r \Pic(C)$ is given by the first Chern class $c_1$. If $\Escr$ is a vector bundle on $C$ then $c_1(\Escr) = \det \Escr$. We also have $c_1(\Oscr_q) = \Oscr_C(q)$ for $q \in C$. The homomorphism $\deg : \Pic(C) \r \ZZ$ assigns to a line bundle its degree. We will denote the composition $\deg \circ\, c_1$ also by $\deg$. If $\Uscr$ is a line bundle then $\deg [\Uscr] := \deg \Uscr$. If $F \in \coh(C)$ has finite length then $\deg [F] = \length F$ \cite[Ex II. 6.12]{H}. 

The functor $u^{\ast}$ induces a group homomorphism
\begin{align*}
u^{\ast}: K_0(X) \ra K_0(C): [\Mscr] \mapsto \sum_{j} (-1)^j [L_j u^{\ast} \Mscr] = [u^{\ast} \Mscr]-[L_1 u^{\ast} \Mscr]
\end{align*}
Recall the basis $\Bscr' = \{[\Oscr],[\Sscr],[\Qscr],[\Pscr]\}$ for $K_{0}(X)$ from \S\ref{The Grothendieck group of $X$}. The image of $\Bscr'$ under $u^{\ast}$ is computed in the following analogue of \cite[Proposition 4.3]{DV1}.
\begin{proposition} 
Assume $A$ is a cubic AS-algebra of generic type A. We have
\begin{align*} 
u^\ast[\Oscr]& =[\Oscr_C] \\ 
u^\ast[\Sscr]& =[\Oscr_p]+ [\Oscr_q] && \text{$p,q \in C$ arbitrary but collinear } \\ 
u^\ast[\Qscr]& =[\Oscr_p]+ [\Oscr_q]+ [\Oscr_{r}]+ [\Oscr_{s}] && \text{$p,q,r,s \in C$ arbitrary but $p+q+r+s = $} \\ & \text{\hspace{2.8cm} $2(p'+q'-\xi)$} && \text{for some collinear $p',q' \in C$} \\
u^\ast[\Pscr]& =[\Oscr_p]-[\Oscr_{\sigma^{-4}p}] &&\text{$p \in C$ arbitrary } 
\end{align*}
\end{proposition}
\begin{proof}
Since $A$ is $h$-torsion free we have $L_1u^\ast\Oscr = 0$ hence 
$u^{\ast}[\Oscr] = [u^{\ast}\Oscr] = [\Oscr_C]$. 
Second, write $S = A/aA$ for some $a \in A_{1}$ and $\Sscr = \pi S$. Application of $- \otimes_{A}D$ on the exact sequence $0 \r A(-1) \xrightarrow{a \cdot} A \r S \r 0$ gives
\[
0 \r \Tor^{A}_{1}(S,D) \r D(-1) \xrightarrow{\overline{a} \cdot} D \r S\otimes_{A}D \r 0 
\]
Since $D$ is a domain the middle map is injective. Hence $\Tor^{A}_{1}(S,D) = 0$ (thus $S$ is $h$-torsionfree) and therefore $L_{1}u^{\ast}\Sscr = 0$. Thus $u^{\ast}[\Sscr] = [u^{\ast}\Sscr]$. From \cite{ATV2} it follows that $u^{\ast}\Sscr = \Oscr_{L}$ where $L$ is the scheme-theoretic intersection of $C$ and the line $\{ a = 0 \} \times \PP^1$. Since $L$ consists of two points $p,q$ we obtain $[\Oscr_{L}] = [\Oscr_p]+ [\Oscr_q]$. By definition $p$ and $q$ are collinear points. This proves the second equality.

Third, same reasoning as above yields $u^{\ast}\Qscr$ is a finite dimensional $\Oscr_C$-module which corresponds to a divisor of degree four on $C$. Thus we may write $[u^{\ast}\Qscr] = [\Oscr_p] + [\Oscr_q] + [\Oscr_r] + [\Oscr_s]$ for some $p,q,r,s \in C$.  It is easy to see that $\Oscr_C(-2) = \sigma^{\ast}\Oscr_C(-1)\otimes_{C}\Oscr_C(-1)$. Since $\Oscr_C(-1) = \Oscr_{C}(-L)$ (from the previous part) we find that its pullback via $\sigma$ is equal to $\sigma^{\ast}\Oscr_{C}(-1) = \Oscr_{C}(-\sigma^{-1} L)$ hence
$\Oscr_C(-2) = \Oscr_{C}(-L-\sigma^{-1} L)$. This means the divisor of $u^{\ast}\Qscr$ is linearly equivalent to $L + \sigma^{-1} L$, which means they have the same sum under the group law of $C$. 

Finally we prove the fourth equation. Put $\Pscr = \Nscr_p$. Now $N_{p} \otimes_{A}D \cong N_{p}/N_{p}h = N_{p}$ thus $u^{\ast}\Nscr_{p} = \widetilde{N_{p}} = \Oscr_{p}$. Applying $N_{p} \otimes_{A}-$ to the short exact sequence of $A$-bimodules $0 \r A(-4) \xrightarrow{\cdot h} A \r D \r 0$
we get the exact sequence
\[
0 \r \Tor_{1}^{A}(N_{p},D) \r N_{p}(-4) \xrightarrow{\cdot \overline{h}} N_{p} \r N_{p}\otimes_{A}D \r 0
\]
As $\overline{h} = 0$ we find $\Tor_{1}^{A}(N_{p},D) = N_{p}(-4) = \left(N_{\sigma^{-4}p}\right)_{\geq 4}$. Thus $L_{1}u^{\ast}\Nscr_{p} =$ \linebreak $\Tor_{1}^{A}(N_{p},D){\,\widetilde{}\,\,}=\left(N_{\sigma^{-4}p}\right){\,\widetilde{}\,\,} = \Oscr_{\sigma^{-4}p}$. This ends the proof of the proposition.
\end{proof} 
Let $\Mscr \in \coh(X)$ and write $[\Mscr] = r[\Oscr] + a[\Sscr] + b[\Qscr] + c[\Pscr]$. By the previous proposition
\begin{equation} \label{rankdegreeC}
\rank u^\ast [\Mscr] = r = \rank \Mscr \quad \text{ and } \deg u^\ast [\Mscr] = 2a + 4b
\end{equation}
We deduce
\begin{proposition} \label{ref-4.3-39cubic} 
Let $A$ be a cubic AS-algebra of generic type A. 
\begin{enumerate} 
\item 
If $\Iscr$ is a line bundle on $X$ then $u^\ast \Iscr$ is a line bundle on $C$, and $\Iscr$ is normalized if and 
only if $\deg u^\ast\Iscr = 0$. 
\item 
If $\Iscr$ is a normalized line bundle on $X$ with invariants $(n_e,n_o)$ then 
\[
c_1(u^\ast\Iscr)=\Oscr_{C}((o)- (2(n_e+n_{o})\xi)) 
\] 
\end{enumerate} 
\end{proposition}
\begin{proof}
The first statement is immediate from the definition of a normalized line bundle on $X$. The second part results from a straightforward computation.
\end{proof}
\begin{corollary} \label{restrEcubic}
Let $A$ be a cubic AS-regular algebra of generic type A and assume $\sigma$ has infinite order. Then the category 
\[
\Rscr(X) = \coprod_{(n_e,n_o) \in N} \Rscr_{(n_e,n_o)}(X) = \{ \text{normalized line bundles on $X$} \}
\]
is equivalent to the full subcategory of $\coh(X)$ with objects 
\[
\{ \Mscr \in \coh(X) \mid u^{\ast}\Mscr \text{ is a line bundle on $C$ of degree zero} \}.
\]
\end{corollary}
\begin{proof}
Due to Proposition \ref{ref-4.3-39cubic} it is sufficient to prove that if $\Mscr \in \coh(X)$ for which $u^{\ast}\Mscr \in \coh(C)$ is a line bundle of degree zero, then $\Mscr$ is a normalized line bundle on $X$. Pick $M \in \grmod(A)$ for which $\pi M = \Mscr$. We may assume $M$ contains no subobject in $\tors(A)$. By Proposition \ref{propvb} and \eqref{rankdegreeC} it suffices to prove $Lu^{\ast}\Mscr = u^{\ast}\Mscr$ i.e. $L_{1}u^{\ast}\Mscr = 0$. 

It is sufficient to prove $M$ is torsion free, since it then follows that $M$ is $h$-torsion free whence $L_1 u^\ast \Mscr=\ker(M(-4)\xrightarrow{\cdot h} M){\,\widetilde{}\,\,} = 0$. So let us assume by contradiction $M$ is not torsion free. Let $T \subset M$ the maximal torsion submodule of $M$. Thus $0 \neq M/T$ is torsion free. Applying $u^{\ast}$ to $0 \r \pi T \r \Mscr \r \pi(M/T) \r 0$ then gives the exact sequence $0 \r u^{\ast}\pi T \r u^{\ast}\Mscr \r u^{\ast}\pi(M/T) \r 0$ on $C$. Since $u^{\ast}\Mscr$ is a line bundle on $C$, it is pure hence either $u^{\ast}\pi T$ is a line bundle or $u^{\ast}\pi T = 0$. 

If $u^{\ast}\pi T$ would be a line bundle then $u^{\ast}\pi (M/T) = (M/T\otimes_{A}D){\,\widetilde{}\,\,}$ has rank zero. 
Thus $M/T\otimes_{A}D \in \grmod(D)$ has GK-dimension $\leq 1$. 
But then $\gkdim M/T \leq 2$, a contradiction with the fact that $M/T \in \grmod(A)$ is non-zero and torsion free. Thus $u^{\ast} \pi T = 0$ i.e. $(T/hT){\,\widetilde{}\,\,} = 0$. This means $\pi(T/hT) = 0$ hence $T/hT \in \tors(A)$. By Lemma \ref{ref-2.8.2-18} we deduce $T \in \tors(A)$ thus $T = 0$ since $M$ contains no subobjects in $\tors(A)$. This ends the proof.
\end{proof}

\begin{remark} \label{RemarkHeis} 
Some of the results above may be generalized to other elliptic cubic AS-algebras. For example, if we consider the situation where $A = H_c$ is the enveloping algebra then one obtains the similar results
\begin{itemize}
\item
If $\Iscr$ is a line bundle on $X$ then $u^\ast\Iscr$ is a line bundle on $C = \Delta$ (the diagonal on $\PP^1 \times \PP^1$) and $L_1u^{\ast}\Iscr=0$. In addition $\Iscr$ is normalized if and only if $u^{\ast}\Iscr$ has degree zero, i.e. if and only if $u^{\ast}\Iscr \cong \Oscr_{\Delta}$ (since $\Pic(\Delta) \cong \ZZ$). 
\item
The category $\Rscr(X) = \coprod_{(n_e,n_o) \in N}\Rscr_{(n_e,n_o)}(X)$ is equivalent to the full subcategory of $\coh(X)$ with objects $\{ \Mscr \in \coh(X) \mid u^{\ast}\Mscr \cong \Oscr_{\Delta} \}$.
\end{itemize}
\end{remark}

\section{From normalized line bundles to quiver representations} \label{Quiver representations}

Throughout Section \ref{Quiver representations}, $A$ will be a cubic AS-algebra. From \S\ref{ref-5.3cubic} onwards we will furthermore assume $A$ is elliptic (and often restrict to the case where $\sigma$ has infinite order). We recycle the notations of Section \ref{Geometric data associated to a cubic AS-algebra} and write $u:C \r X$ for the map of noncommutative schemes as defined in Section \ref{Restriction of line bundles to the divisor $C$}.

\subsection{Generalized Beilinson equivalence} \label{ref-5.1cubic}

We set $\Escr = \Oscr(3) \oplus \Oscr(2) \oplus \Oscr(1) \oplus \Oscr$ and $U = \Hom_{X}(\Escr, \Escr) = \bigoplus_{i,j = 0}^{3}{\Hom_{X}(\Oscr(i),\Oscr(j))}$. The functor $\Hom_{X}(\Escr,-)$ from $\coh(X)$ to the category $\mod(U)$ of right $U$-modules extends to an equivalence $\RHom_{X}(\Escr,-)$ of bounded derived categories \cite{Bondal1}
\begin{equation} \label{Bondalcubic}
\begin{diagram}[heads=LaTeX] 
\D^{b}(\coh(X))
& \pile{
\rTo^{\RHom_{X}(\Escr,\un)} \\
\vspace{0.3cm} \\
\lTo_{\un \LotimesU \Escr}
}
& \D^{b}(\mod(U)) 
\end{diagram}
\end{equation}
where the inverse functor is given by $-\LotimesU \Escr$. For the classical case $X = \PP^n$ such an equivalence was found by Beilinson \cite{Beilinson}. We refer to \eqref{Bondalcubic} as \emph{generalized Beilinson equivalence}. For a non-negative integer $i$ this equivalence restricts to an equivalence \cite{Baer} between the full subcategories $\Xscr_{i} = \{ \Mscr \in \coh(X) \mid \Ext_{X}^{j}(\Escr,\Mscr) = 0 \text{ for } j \neq i \}$ and $\Yscr_{i} = \{ M \in \mod(U) \mid \Tor_{j}^{U}(M,\Escr) = 0 
\text{ for } j \neq i \}$. The inverse equivalences are given by $\Ext_{X}^{i}(\Escr, -)$ and $\Tor_{i}^{U}(-, \Escr)$. 

It is easy to see that $U \cong k\Gamma/(R)$ where $k\Gamma$ is the path algebra of the quiver $\Gamma$
\begin{equation} \label{ref-5.2-42cubic}
\begin{diagram}[heads=LaTeX] 
-3
& \pile{
\rTo^{X_{-3}} \\
\rTo^{Y_{-3}} 
}
& 
-2
& \pile{
\rTo^{X_{-2}} \\
\rTo^{Y_{-2}} 
}
& -1 & \pile{
\rTo^{X_{-1}} \\
\rTo^{Y_{-1}} 
}
& 0
\end{diagram}
\end{equation}
with relations $R$ reflecting the relations of $A$. If we write the relations of $A$ as \eqref{relations}
then the relations $R$ are given by 
\begin{eqnarray} \label{relationsGammacubic} 
\begin{pmatrix}
X_{-1} & Y_{-1}
\end{pmatrix}
\cdot M_{A}^{t}(X_{-2},Y_{-2},X_{-3},Y_{-3})
= 0
\end{eqnarray}
where $M_{A}^{t}(X_{-2},Y_{-2},X_{-3},Y_{-3})$ is obtained from the matrix $M_A^{t}$ by replacing $x^2$, $xy$, $yx$ and $y^2$ by $X_{-2}X_{-3},Y_{-2}X_{-3},X_{-2}Y_{-3}$ and $Y_{-2}Y_{-3}$. 

We write $\Mod(\Gamma)$ for the category of representations of the quiver $\Gamma$, where representations are assumed to satisfy the relations. For $M \in \Mod(\Gamma)$ we write $M_{i}$ for the $k$-linear space located at vertex $i$ of $\Gamma$ and $M(X_{i})$, $M(Y_i)$ for the linear maps corresponding to arrows $X_{i}$, $Y_{i}$ of $\Gamma$ ($i=-3,\dots,0$). We denote $S_i$ for the simple representation corresponding to $i$. 
Since the category $\Mod(\Gamma)$ of representations of $\Gamma$ is equivalent to the category of right $k\Gamma/(R)$-modules we deduce $\Mod(\Gamma)\cong \Mod(U)$. From now on we write $\Mod(\Gamma)$ instead of $\Mod(U)$. One verifies that the matrix representation of the Euler form $\chi: K_{0}(\Gamma) \times K_{0}(\Gamma) \ra \ZZ$ with respect to the basis $\{ S_{-3},S_{-2},S_{-1},S_{0} \}$ of $K_{0}(\Gamma)$ is given by 
\begin{equation} \label{ref-5.4-44cubic} 
\begin{pmatrix} 
1 & -2 & 0  & 2 \\ 
0 & 1  & -2 & 0  \\ 
0 & 0  & 1  & -2 \\
0 & 0  & 0  & 1
\end{pmatrix}. 
\end{equation} 

\subsection{Point, line and conic representations}

For further use we determine the representations of $\Gamma$ corresponding to point, line and conic objects on $X$. The following lemmas are proved in the same spirit as \cite[Lemmas 5.1.3 and 5.6.1]{DV1}.
\begin{lemma} \label{ref-5.1.3-47cubic} 
Let $p \in C$ and put $(\alpha_{i}: \beta_{i}) = \pr_1(\sigma^{i}p) \in \PP^1$. 
\begin{enumerate} 
\item  
$H^{j}(X,\Nscr_{p}(m)) = 0$ for all integers $m$ and $j > 0$. In particular $\Nscr_{p} \in \Xscr_{0}$. 
\item 
$\dim_{k} \left(\omega \Nscr_{p}\right)_{m} = 1$ for all $m$ and $\left(\omega \Nscr_{p}\right)_{\ge m}$ is a shifted point module for all integers $m$. In particular $\left(\omega \Nscr_{p}\right)_{\ge 0} = N_{p}$. 
\item 
$H^0(X,\Nscr_{p}(m))=\left(\omega \Nscr_{p}\right)_m$ for all integers $m$. 
\item 
Write $\RHom_{X}(\Escr,\Nscr_{p})=p$. Then $\underline{\dim}p = (1,1,1,1)$ and $p \in \mod(\Gamma)$ corresponds to the representation
\begin{equation*}
\begin{diagram}[heads=LaTeX]
k
& \pile{
\rTo^{\alpha_{-3}} \\
\rTo^{\beta_{-3}} 
}
& k
& \pile{
\rTo^{\alpha_{-2}} \\
\rTo^{\beta_{-2}} 
}
& k & \pile{
\rTo^{\alpha_{-1}} \\
\rTo^{\beta_{-1}} 
}
& k
\end{diagram}
\end{equation*} 
\end{enumerate} 
\end{lemma} 
\begin{lemma} \label{conicreps}
Let $n \geq 1$ be an integer, $w \in A_{n}$ and put $\Wscr = \pi(A/wA)$.
\begin{enumerate}
\item
$H^{1}(X,\Wscr(m)) \cong (A/Aw)_{-m-2}'$ for $m \leq -1$
\item
$H^j(X,\Wscr(m))=0$ for $m \leq -1$ and $j\neq 1$. In particular $\Wscr(-1) \in \Xscr_{1}$. 
\item
If $\eta\in A_{1}$ then the induced linear map $H^1(X,\Wscr(m))\xrightarrow{\cdot\eta} H^1(X,\Wscr(m+1))$ corresponds to $(\eta\cdot)'$ on $(A/Aw)'$. 
\item
Write $\RHom_{X}(\Escr,\Wscr(-1))=W[-1]$. Then 
\begin{eqnarray*}
\underline{\dim}W 
=
\left\{
\begin{array}{ll}
(2,1,1,0) & \text{ if } n = 1 \\
(3,2,1,0) & \text{ if } n = 2 \\
(4,2,1,0) & \text{ if } n > 2
\end{array}
\right.
\end{eqnarray*}
and $W \in \mod(\Gamma)$ corresponds to the representation
\begin{equation} \label{ref-5.10-67cubic}
\begin{diagram}[heads=LaTeX]
(A/Aw)_{2}'
& \pile{
\rTo^{(x\cdot)'} \\
\rTo^{(y\cdot)'} 
}
&
(A/Aw)_{1}'
& \pile{
\rTo^{(x\cdot)'} \\
\rTo^{(y\cdot)'} 
}
& k & \pile{
\rTo^{} \\
\rTo^{} 
}
& 0
\end{diagram}
\end{equation}
\end{enumerate}
\end{lemma}

\subsection{First description of $\Rscr_{(n_e,n_o)}(X)$} \label{ref-5.3cubic}

From now on we assume in Section \ref{Quiver representations} $A$ is an elliptic cubic AS-algebra. As in \eqref{defN} we put $N = \{(n_e,n_o) \in \NN^2 \mid n_e - (n_e - n_o)^2 \geq 0 \}$. Recall from \S\ref{Normalized line bundles} the set $R(A)$ of reflexive rank one graded right $A$-modules considered up to isomorphism and shift is in natural bijection with the isoclasses in the category $\coprod_{(n_e,n_o)\in N}\Rscr_{(n_e,n_o)}(X)$ where $\Rscr_{(n_e,n_o)}(X)$ is the full subcategory of $\coh(X)$ consisting of the normalized line bundles on $X$ with invariants $(n_e,n_o)$.  

Let $\Iscr$ be an object of $\Rscr_{(n_e,n_o)}(X)$, considered as a complex in $\D^{b}(\coh(X))$ of degree zero. Theorem \ref{cohomology} implies $\Iscr \in \Xscr_{1}$. Thus the image of this complex is concentrated in degree one i.e. $\RHom_{X}(\Escr, \Iscr) = M[-1]$ where $M = \Ext^{1}_{X}(\Escr, \Iscr)$ is a representation of $\Delta$. By functoriality, multiplication by $x,y \in A$ induces linear maps $M(X_{-i}),M(Y_{-i}): H^1(X,\Iscr(-i)) \r H^1(X,\Iscr(-i+1))$ hence $M$ is given by the following representation of $\Gamma$ 
\begin{equation*}
\begin{diagram}[heads=LaTeX]
H^1(X,\Iscr(-3))
& \pile{
\rTo^{M(X_{-3})} \\
\rTo^{M(Y_{-3})} 
}
& H^1(X,\Iscr(-2)) & \pile{
\rTo^{M(X_{-2})} \\
\rTo^{M(Y_{-2})} 
}
& H^1(X,\Iscr(-1)) & \pile{
\rTo^{M(X_{-1})} \\
\rTo^{M(Y_{-1})} 
}
& H^1(X,\Iscr)
\end{diagram}
\end{equation*}
We denote $\Cscr_{(n_e,n_o)}(\Gamma)$ for the image of $\Rscr_{(n_e,n_o)}(X)$ under the equivalence $\Xscr_{1} \cong \Yscr_{1}$. In an analogue way as \cite[Theorem  5.3.1]{DV1} we obtain
\begin{theorem} \label{firstdescrcubic}
Let $A$ be an elliptic cubic AS-algebra where $\sigma$ has infinite order. Let $(n_e,n_o) \in N \setminus \{(0,0)\}$. Then there is an equivalence of categories 
\begin{diagram}[heads=LaTeX]
\Rscr_{(n_e,n_o)}(X) 
& \pile{
\rTo^{\Ext_{X}^{1}(\Escr, \un)} \\
\vspace{0.3cm} \\
\lTo_{\Tor_{1}^{{{\Gamma}}}(\un, \Escr)}
}
& \Cscr_{(n_e,n_o)}(\Gamma) 
\end{diagram}
where 
\begin{multline} \label{cubicdescr}
\Cscr_{(n_e,n_o)}(\Gamma) = \{ M \in \mod(\Gamma) \mid \underline{\dim} M=(n_o,n_e,n_o,n_e-1) 
\text{ and } \\ 
\Hom_{\Gamma}(M,p) = 0, \Hom_{\Gamma}(p,M) = 0 \mbox{ for all } p \in C \}. 
\end{multline} 
\end{theorem} 
\begin{proof}
First, let $\Iscr$ be an object of $\Rscr_{(n_e,n_o)}(X)$ and write $M = \Ext^{1}_{X}(\Escr, \Iscr)$. That $\underline{\dim}M = (n_o,n_e,n_o,n_e-1)$ follows from Theorem \ref{cohomology}. Further, let $p \in C$ 
and as in Lemma \ref{ref-5.1.3-47cubic} we denote $p = \Hom_{X}(\Escr,\Nscr_{p})$. 

Lemma \ref{reflexive} implies $\Ext^{i}_{X}(\Iscr,\Nscr_{p}) = 0$ for $i \geq 1$ hence $\chi(\Iscr,\Nscr_{p}) = 1$ yields $k = \RHom_{X}(\Iscr,\Nscr_{p}) \cong \RHom_{\Gamma}(M[-1],p)$ where we have used \eqref{Bondalcubic}. In particular $\Hom_{\Gamma}(M,p) = H^{1}(\RHom_{\Gamma}(M[-1],p)) = 0$. Similarly $k[-2] = \RHom_{X}(\Nscr_{p},\Iscr) \cong \RHom_{\Gamma}(p,M[-1])$ whence $\Hom_{\Gamma}(p,M) = 0$. 

Conversely let $M \in \mod(\Delta)$ for which $\underline{\dim} M=(n_o,n_e,n_o,n_e-1)$ and \linebreak $\Hom_{\Gamma}(M,p) = \Hom_{\Gamma}(p,M) = 0$ for all $p \in C$. By Serre duality (Theorem \ref{Serreduality})
\begin{align*}
H^{2}(\RHom_{\Gamma}(M,p)) & = H^{2}(\RHom_{X}(M \LotimesGamma \Escr, \Nscr_{p})) \\
& \cong H^{0}(\RHom_{X}(\Nscr_{\sigma^{4}p},M \LotimesGamma \Escr)) = H^{0}(\RHom_{\Gamma}(\sigma^{4}p,M))
\end{align*}
Thus $\Hom_{\Gamma}(M,p) = \Ext^{2}_{\Gamma}(M,p) = 0$ for all $p \in C$. Now $\gldim \mod(\Gamma) = 2$ so we may compute $\dim_{k} \Ext^1_{\Gamma}(M,p)$ using the Euler form \eqref{ref-5.4-44cubic} on $\mod(\Gamma)$. We obtain $\chi(p,M) = -1$ hence $\Ext^1_{\Gamma}(M,p) = k$. In other words $\RHom_{\Gamma}(M[-1],p) = k$. 

Put $\Iscr = M[-1] \LotimesGamma \Escr$. By the generalized Beilinson equivalence \eqref{Bondalcubic}, $\Iscr \in D^{b}(\coh(X))$ and $k = \RHom_{\Gamma}(M[-1],p) = \RHom_{X}(\Iscr, \Nscr_{p})$, giving (by adjointness) $\RHom_{C}(Lu^{*}\Iscr, \Oscr_{p}) = k$. Thus $\RHom_{\Oscr_{C,p}}(Lu^{*}\Iscr_{p}, k) = k$ for all $p \in C$. From this we easily deduce $Lu^\ast \Iscr$ is a line bundle on $C$. 
Hence by Proposition \ref{propvb} we find $\Iscr$ is a vector bundle on $X$. In particular, $M \in \Yscr_{1}$. What is left to check is that $\Iscr$ is a normalized line bundle. 
The derived equivalence \eqref{Bondalcubic} gives rise to inverse group isomorphisms (see for example \cite[Proposition 3.2.3]{Baer})
\begin{align*}
\mu: K_{0}(X) \r K_{0}(\Gamma): & [\Nscr] \mapsto \sum_{i}(-1)^{i}[\Ext^{i}_{X}(\Escr,\Nscr)] \\
\nu : K_{0}(\Gamma) \r K_{0}(X): & [N] \mapsto \sum_{i}(-1)^{i}[\Tor_{i}^{\Gamma}(N,\Escr)]
\end{align*}
One verifies that the image of the basis $\Bscr' = \{ [\Oscr],[\Sscr],[\Qscr],[\Pscr] \}$ for $K_{0}(X)$ under $\mu$ is the $\ZZ$-basis $\{[S_{0}],-[S_{-3}]+[S_{0}],-2[S_{-3}]-[S_{-2}]+[S_{0}],[S_{-3}]+[S_{-2}]+[S_{-1}]+[S_{0}]\}$ for $K_{0}(\Gamma)$. This may be done by base change from $\Bscr'$ to the basis $\Bscr$ (see Propositions \ref{Grothendieck1} and \ref{basis2}) and using Theorem \ref{fullcohomology}.
Since $M \in \Yscr_{1}$ we have $\nu [M] = - [\Iscr]$. One now computes $[\Iscr] = [\Oscr] -2(n_{e} - n_{o})[\Sscr] + (n_{e} - n_{o})[\Qscr] - n_{o}[\Pscr]$. 
%
%
We conclude $\Iscr \in \Rscr_{(n_e,n_o)}(X)$, completing the proof.
\end{proof}
We may now parameterize the line bundles on $X$ for some low invariants.
\begin{corollary} \label{smallinvcubic}
Let $A$ be an elliptic cubic AS-algebra where $\sigma$ has infinite order.
\begin{enumerate}
\item
The category $\Cscr(1,0)$ consists of one object namely the simple object $S_{-2}$
\begin{eqnarray*}  
0
\begin{array}{c}
{\stackrel{0}{\longrightarrow}} \\
{\stackrel{0}{\longrightarrow}} 
\end{array}
k
\begin{array}{c}
{\stackrel{0}{\longrightarrow}} \\
{\stackrel{0}{\longrightarrow}} 
\end{array}
0
\begin{array}{c}
{\stackrel{0}{\longrightarrow}} \\
{\stackrel{0}{\longrightarrow}} 
\end{array}
0
\end{eqnarray*}
\item
The representations in $\Cscr(1,1)$ are the representations of $\Gamma$ of the form
\begin{eqnarray*} 
k
\begin{array}{c}
{\stackrel{\alpha}{\longrightarrow}} \\
{\stackrel{\beta}{\longrightarrow}} 
\end{array}
k
\begin{array}{c}
{\stackrel{\alpha'}{\longrightarrow}} \\
{\stackrel{\beta'}{\longrightarrow}} 
\end{array}
k
\begin{array}{c}
{\stackrel{0}{\longrightarrow}} \\
{\stackrel{0}{\longrightarrow}} 
\end{array}
0
\end{eqnarray*}
where $((\alpha:\beta),(\alpha':\beta')) \in (\PP^{1} \times \PP^{1}) - C$.
\end{enumerate}
\end{corollary}
\begin{proof}
This is easy to deduce from Theorem \ref{firstdescrcubic}.
\end{proof}

\subsection{Description of $\Rscr_{(n_e,n_o)}(X)$ for the enveloping algebra}

In this section we let $A$ be the enveloping algebra $H_c$. Thus $C = E_{\red}$ is the diagonal $\Delta$ on $\PP^1 \times \PP^1$. Recall from \S\ref{Geometric data associated to a cubic AS-algebra} that the restriction $\sigma_{\Delta}$ is the identity. Our proof of the next lemma is in the same spirit as the proof of \cite[Theorem 4.5(i)]{BW} for the homogenized Weyl algebra. 
\begin{lemma} \label{isoHeis}
Let $\Iscr \in \Rscr_{(n_{e},n_{o})}(X)$ for some $(n_e,n_o) \in N \setminus \{(0,0)\}$. Consider for any integer $m$ the linear map $M(Z_{m})$ induced by multiplication by $z = xy - yx$
\[
M(Z_{m}): H^{1}(X, \Iscr(-m)) \ra H^{1}(X,\Iscr(-m+2))
\]
Then $M(Z_{m})$ is surjective for $m < 4$ and injective for $m > 2$. 
\end{lemma}
\begin{proof}
Let $m$ be any integer and put $\Qscr = \pi(A/zA) = \pi D$. Then $u^{\ast}\Qscr = \Oscr_{\Delta}$. Applying $\Hom_{X}(-,\Iscr)$ to $0 \ra \Oscr_{X}(m-2) \stackrel{z}{\ra} \Oscr_{X}(m) \ra \Qscr(m) \ra 0$ yields
\[
\Ext^{1}_{X}(\Qscr(m),\Iscr) \ra H^{1}(X,\Iscr(-m)) \xrightarrow{M(Z_{m})} H^{1}(X,\Iscr(-m+2)) \ra \Ext^{2}_{X}(\Qscr(m),\Iscr).
\]
Furthermore Theorem \ref{Serreduality} (Serre duality) implies
\begin{align*}
\Ext^{1}_{X}(\Qscr(m),\Iscr) \cong  \Ext^{1}_{X}(\Iscr,\Qscr(m-4))', \,\,
\Ext^{2}_{X}(\Qscr(m),\Iscr) \cong \Hom_{X}(\Iscr,\Qscr(m-4))'.
\end{align*}
On the other hand since $\RHom_{X}(\Iscr,u_{\ast}\Oscr_{\Delta}(m-4)) \cong \RHom_{\Gamma}(Lu^{\ast}\Iscr, \Oscr_{\Delta}(m-4))$ by \eqref{Bondalcubic} and $Lu^{\ast}\Iscr=\Oscr_{\Delta}$ (see Remark \ref{RemarkHeis}) we derive
\[
\Hom_{X}(\Iscr,\Qscr(m-4)) = H^{0}(\Delta,\Oscr_{\Delta}(m-4)) = D_{m-4} = 0 \text{ for } m < 4
\]
and by Serre duality on $\Delta$
\[
\Ext^{1}_{X}(\Iscr,\Qscr(m-4)) = H^{1}(\Delta,\Oscr_{\Delta}(m-4)) \cong D_{-m+2}'  = 0 \text{ for } m > 2
\]
which completes the proof.
\end{proof}

\begin{theorem} \label{firstdescrHeis}
Let $A = H_c$ be the enveloping algebra. Let $(n_e,n_o) \in N \setminus \{(0,0)\}$. Define for any $M \in \mod(\Gamma)$ the linear maps
\begin{align*}
& M(Z_{-3}) = M(Y_{-2})M(X_{-3}) - M(X_{-2})M(Y_{-3}) \\
& M(Z_{-2}) = M(Y_{-1})M(X_{-2}) - M(X_{-1})M(Y_{-2})
\end{align*}
There is an equivalence of categories 
\begin{diagram}[heads=LaTeX]
\Rscr_{(n_e,n_o)}(X) 
& \pile{
\rTo^{\Ext_{X}^{1}(\Escr, \un)} \\
\vspace{0.3cm} \\
\lTo_{\Tor_{1}^{{{\Gamma}}}(\un, \Escr)}
}
& \Cscr_{(n_e,n_o)}(\Gamma) 
\end{diagram}
where 
\begin{multline} \label{Heisdescr}
\Cscr_{(n_e,n_o)}(\Gamma) = \{ M \in \mod(\Gamma) \mid \underline{\dim} M=(n_o,n_e,n_o,n_e-1) 
\text{ and } \\ 
\text{$M(Z_{-3})$ isomorphism and $M(Z_{-2})$ surjective} \}. 
\end{multline} 
\end{theorem}
\begin{proof}
Due to Theorem \ref{firstdescrcubic} it will be sufficient to prove that the descriptions \eqref{cubicdescr} \eqref{Heisdescr} coincide. One inclusion follows directly from Lemma \ref{isoHeis}, so let us assume $M \in \mod(\Gamma)$ for which $M(Z_{-3})$ is an isomorphism and $M(Z_{-2})$ is surjective. Let $p = ((\alpha:\beta),(\alpha:\beta)) \in \Delta$ and write $p \in \mod(\Gamma)$ for the corresponding representation of the quiver $\Gamma$. Let $\tau = (\tau_{-3},\tau_{-2},\tau_{-1},\tau_{0}) \in \Hom_{\Gamma}(p,M)$ be any morphism. Thus we have a commutative diagram in $\mod(k)$
\begin{equation*}
\begin{diagram}[heads=LaTeX]
k
& \pile{
\rTo^{\alpha} \\
\rTo^{\beta} 
}
& k
& \pile{
\rTo^{\alpha} \\
\rTo^{\beta} 
}
& k
& \pile{
\rTo^{\alpha} \\
\rTo^{\beta} 
}
& k
\\
\dTo^{\tau_{-3}}
&&
\dTo^{\tau_{-2}}
&&
\dTo^{\tau_{-1}}
&&
\dTo^{\tau_{0}} \\
M_{-3}
& \pile{
\rTo^{M(X_{-3})} \\
\rTo^{M(Y_{-3})} 
}
& M_{-2}
& \pile{
\rTo^{M(X_{-2})} \\
\rTo^{M(Y_{-2})} 
}
& M_{-1}
& \pile{
\rTo^{M(X_{-1})} \\
\rTo^{M(Y_{-1})} 
}
& M_{0}
\end{diagram}
\end{equation*}
We claim $\tau_{-3} = 0$. Assume by contradiction this is not the case. Since $M(Z_{-3})$ is an isomorphism we surely have $v = M(Z_{-3}) \tau_{-3}(1) \neq 0$. On the other hand, by the commutativity of the above diagram
\[
M(Z_{-3})\tau_{-3}(1) = \left( M(Y_{-2})M(X_{-3}) - M(X_{-2})M(Y_{-3})\right) \tau_{-3}(1) = \tau_{-1}(\alpha \beta - \beta \alpha) = 0
\]
leading to the desired contradiction. Thus $\tau_{-3} = 0$. It follows that $\tau = 0$. 

By similar arguments we may show for $\tau = (\tau_{-3},\tau_{-2},\tau_{-1},\tau_{0}) \in \Hom_{\Gamma}(M,p)$ we have $\tau_0=0$, which implies $\tau=0$. 
\end{proof}

\subsection{Restriction to a full subquiver} \label{ref-5.6cubic}

Let $A$ be an elliptic cubic AS-algebra. Although the description of $\Cscr_{(n_e,n_o)}(\Gamma)$ in Theorem \ref{firstdescrcubic} is quite elementary, it is not easy to handle. Similar as in \cite{DV1,LeBruyn} for quadratic AS-algebras we show representations in $\Cscr_{(n_e,n_o)}(\Gamma)$ are completely determined by the four leftmost maps.

Let $\Gamma^{0}$ be the full subquiver of $\Gamma$ consisting of the vertices $-3,-2,-1$ in \eqref{ref-5.2-42cubic}. Let $\Res : \Mod(\Gamma) \ra \Mod(\Gamma^{0})$ be the obvious restriction functor. $\Res$ has a left adjoint which we denote by $\Ind$. 
Note $\Res \circ \Ind = \Id$. If $M \in \Mod(\Gamma)$ we will denote $\Res M$ by $M^{0}$. 

In general, two objects $A$ and $B$ of an abelian category $\Cscr$ are called {\em perpendicular}, denoted by $A \perp B$, if $\Hom_{\Cscr}(A,B) = \Ext^{1}_{\Cscr}(A,B) = 0$. For an object $B \in \Cscr_{f}$ we define ${}^\perp B$ as the full subcategory of $\Cscr_{f}$ which objects are
\[
{}^\perp B = \{ A \in \Cscr_{f} \mid A \perp B \}.
\]
Repeating the arguments from the proof of \cite[Lemmas 5.1.2]{DV1} we have $M = \Ind \Res M$ for $M \in \mod(\Gamma)$ if and only if $M \perp S_{0}$. This means the functors $\Res$ and $\Ind$ define inverse equivalences \cite{Baer}
\begin{equation} \label{equivalenceresubic}
\begin{diagram}[heads=LaTeX]
\mod(\Gamma) \supset {}^\perp S_{0} & \pile{
\rTo^{\Res} \\
\vspace{0.3cm} \\
\lTo_{\Ind}
}
& \mod(\Gamma^{0})
\end{diagram}
\end{equation}
\begin{lemma} \label{ref-5.1.2-46cubic} 
Let $(n_e,n_o) \in N\setminus \{ (0,0) \}$. Then $\Cscr_{(n_e,n_o)}(\Gamma) \subset {}^\perp S_{0}$.
\end{lemma}
\begin{proof}
Similar as the proof of \cite[Lemmas 5.5.1 and 5.5.2]{DV1}. See also \cite{LeBruyn}.
\end{proof}
\begin{lemma} \label{pointperp}
Let $p \in C$ and $\Qscr$ be a conic object on $X$. Write $p = \Hom_{X}(\Escr, \Nscr_{p})$ and $Q = \Ext^{1}_{X}(\Escr, \Qscr(-1))$. Then $p \perp S_0$ and $Q \perp S_0$.
\end{lemma}
\begin{proof}
That $p,Q \in \mod(\Gamma)$ follows from Lemmas \ref{ref-5.1.3-47cubic} and \ref{conicreps}. By \eqref{Bondalcubic} we have $\Ext^{i}_{\Gamma}(p,S_0) \cong \Ext^{i}_{X}(\Nscr_p,\Oscr) = 0$ for $i = 0,1$. Similarly $\Ext^{i}_{\Gamma}(Q,S_0) =$ \linebreak $\Ext^{i-1}_{\Gamma}(Q[-1],S_0) \cong \Ext^{i-1}_{X}(\Qscr,\Oscr) = 0$ for $i = 0,1$. This proves what we want.
\end{proof}

\subsection{Stable representations} \label{ref-5.7cubic}

Our next objective is to show that the representations in $\Cscr_{(n_e,n_o)}(X)$ restricted to $\Gamma^{0}$ are stable. We first recall some generalities on (semi)stable quiver representations. For more details we refer to \cite{king, schofield, SVdB}.

Let $Q$ be a quiver without oriented cycles and let $\theta \in \ZZ^{Q_{0}}$ be a dimension vector. A representation $F$ of $Q$ is called $\theta$-{\em semistable} (resp. {\em stable}) if $\theta \cdot \underline{\dim} F = 0$ and $ \theta \cdot 
\underline{\dim} N \geq 0$ (resp. $> 0$) for every non-trivial subrepresentation $N$ of $F$. Here we denote ``$\cdot$'' for the standard scalar product on $\ZZ^{Q_0}$: $(\alpha_v)_v\cdot (\beta_v)_v=\sum_v \alpha_v\beta_v$. 
The full subcategory of $\theta$-semistable representations of $Q$ forms an exact abelian subcategory of $\mod(Q)$ 
in which the simple objects are precisely the stable representations.

It is known \cite[Corollary 1.1]{SVdB} that $F$ is semistable for some $\theta$ if and only there exists an object $G \in \mod(Q)$ for which $F \perp G$. The relation between $\theta$ and $\underline{\dim} G$ is such that the forms $- \cdot \theta$ and $\chi(-,\underline{\dim}G)$ are proportional. Associated to $G\in \mod(Q)$ there is a semi-invariant function $\phi_G$ on $\Rep_{\alpha}(Q)$ for which the set
\begin{equation} \label{ref-equation1.1} 
\{ F \in \Rep_{\alpha}(Q) \mid F \perp G \} 
\end{equation} 
coincides with $\{\phi_G\neq 0\}$. In particular \eqref{ref-equation1.1} is affine.
Put $Q = \Gamma^{0}$. The following lemmas are elementary. We omit the proofs.
\begin{lemma} \label{stablepointcubic}
Let $p\in C$. Then $\Res p \in \mod(\Gamma^{0})$ is $\theta$-stable for $\theta = (-1,0,1)$.
\end{lemma}
\begin{lemma} \label{stablesurjectivecubic}
Assume $0 \neq F,G \in \mod(\Gamma^{0})$ are $\theta$-semistable for $\theta = (-1,0,1)$. 
\begin{enumerate}
\item
If $G$ is $\theta$-stable then every non-zero map in $\Hom_{\Gamma^{0}}(F,G)$ is surjective.
\item
If $F$ is $\theta$-stable then every non-zero map in $\Hom_{\Gamma^{0}}(F,G)$ is injective.
\end{enumerate} 
\end{lemma}
\begin{proposition} \label{lemcubic}
Let $\Qscr = \pi(A/vA)$ be a conic object on $X$ where $v = \alpha x^2 + \beta xy + \gamma yx + \delta y^2 \in A_{2}$ and write $Q = \Ext^{1}_{X}(\Escr, \Qscr(-1)) \in \mod(\Gamma)$. Let $(n_e,n_o) \in N \setminus\{(0,0)\}$, $\Iscr \in \Rscr_{(n_e,n_o)}(X)$ and write $M = \Ext^{1}_{X}(\Escr, \Iscr) \in \mod(\Gamma)$. 
Then the following are equivalent:
\begin{enumerate}
\item
$M^{0} \perp Q^{0}$ 
\item
$\Hom_{\Gamma^{0}}(M^{0},Q^{0}) = 0$
\item
$\Hom_{X}(\Iscr, \Qscr(-1)) = 0$
\item
$\Iscr \perp \Qscr(-1)$
\item
The following linear map is an isomorphism
\begin{eqnarray} \label{eqf}
\begin{aligned}
f = & \alpha M^{0}(X_{-2})M^{0}(X_{-3}) + \beta M^{0}(Y_{-2})M^{0}(X_{-3}) \\
& + \gamma M^{0}(X_{-2})M^{0}(Y_{-3}) + \delta M^{0}(Y_{-2})M^{0}(Y_{-3}) : M_{-3} \ra M_{-1} 
\end{aligned}
\end{eqnarray}
\end{enumerate}
\end{proposition}
\begin{proof}
By definition $(1)$ implies $(2)$ and its converse is seen by $\chi(M^0,Q^0) = 0$. The equivalence $(2) \Leftrightarrow (3)$ follows from \eqref{Bondalcubic} as 
\begin{align*}
\Hom_{\Gamma^{0}}(M^{0},Q^{0}) & = \Hom_{\Gamma}(\Ind M^{0}, Q) = \Hom_{\Gamma}(M, Q) = H^{0}(\RHom_{\Gamma}(M,Q)) \\
& \cong H^{0}(\RHom_{X}(\Iscr,\Qscr(-1))) = \Hom_{X}(\Iscr, \Qscr(-1)).
\end{align*}
To prove $(3) \Rightarrow (4)$, as $\Iscr$ is a normalized line bundle with invariants $(n_e,n_o)$ we may write $[\Iscr] = [\Oscr] -2(n_{e} - n_{o})[\Sscr] + (n_{e} - n_{o})[\Qscr] - n_{o}[\Pscr]$. Furthermore \eqref{shiftcubic} yields $[\Qscr(-1)] = [\Qscr] - [\Pscr]$ and using Proposition \ref{basis2} one computes $\chi(\Iscr,\Qscr(-1)) = 0$. Since Serre duality gives $\Ext^{2}_{X}(\Iscr, \Qscr(-1)) \cong \Hom_{X}(\Qscr(3),\Iscr)' = 0$ we conclude $\Iscr \perp \Qscr(-1)$ if and only if $\Hom_{X}(\Iscr,\Qscr(-1)) = 0$.

Finally we prove the equivalence between $(4)$ and $(5)$. Applying $\Hom_{X}(-,\Iscr)$ to $0 \ra \Oscr(1) \ra \Oscr(3) \ra \Qscr(3) \ra 0$
gives a long exact sequence of $k$-vector spaces
\[
0 \ra \Ext^{1}_{X}(\Qscr(3),\Iscr) \ra M_{-3} \xrightarrow{f} M_{-1} \ra \Ext^{2}_{X}(\Qscr(3),\Iscr) \ra 0 
\]
where we have used Theorem \ref{cohomology}. As the middle map $f$ is given by \eqref{eqf} we deduce $f$ is an isomorphism if and only if $\Ext^{1}_{X}(\Qscr(3),\Iscr) = 0 = \Ext^{2}_{X}(\Qscr(3),\Iscr)$. Invoking Serre duality on $X$ the latter is equivalent with $\Iscr \perp \Qscr(-1)$. 
\end{proof} 
\begin{remark} \label{envelopingorthQ}
In case $A = H_c$ is the enveloping algebra we recover the property $M^{0}(Z_{-2})$ being an isomorphism (Theorem \ref{firstdescrHeis}), as for the conic object $\Qscr = \pi(H_c/zH_c)$
\[
\RHom_{X}(\Iscr,\Qscr(-1)) = \RHom_{X}(\Iscr,u_{\ast}\Oscr_{\Delta}(-1))
\cong \RHom_{\Delta}(Lu^{\ast}\Iscr, \Oscr_{\Delta}(-1))  
\]
and since $Lu^{\ast}\Iscr = \Oscr_{\Delta}$ we obtain $\Hom_{X}(\Iscr,\Qscr(-1)) \cong \Hom_{\Delta}(\Oscr_{\Delta}, \Oscr_{\Delta}(-1)) = 0$. As a consequence the restriction of the representations in $\Cscr_{(n_e,n_o)}(\Gamma)$ to $\Gamma^{0}$ are $\theta$-semistable for some $\theta \in \ZZ^{3}$. Since $\chi(-,\underline{\dim}Q^{0}) = - \cdot (-1,0,1)$ we may take $\theta = (-1,0,1)$.
\end{remark}
Inspired by the previous remark one might try to find, for all elliptic cubic AS-algebras $A$, a conic object $\Qscr$ on $X$ for which $\Hom_{X}(\Iscr, \Qscr(-1))$ is zero for all $\Iscr \in \Rscr_{(n_e,n_o)}(X)$. We did not manage to find such a conic object independent of $\Iscr$. However, we are able to prove that for a fixed normalized line bundle $\Iscr$ on $X$ there is at least one conic object $\Qscr$ (which depends on $\Iscr$) for which $\Hom_{X}(\Iscr, \Qscr(-1)) = 0$. We will then show how this leads to a proof that the representations in $\Cscr_{(n_e,n_o)}(X)$ restricted to $\Gamma^{0}$ are stable.
\begin{proposition} \label{curve} 
Assume $k$ is uncountable and $\sigma$ has infinite order. Let $(n_e,n_o) \in N$ such that $(n_e-1,n_o-1) \in N$. Let $\Iscr \in \Rscr_{(n_e,n_o)}(X)$. Then the set of conic objects $\Qscr$ for which $\Hom_{X}(\Iscr,\Qscr(-1))\neq 0$ is a curve of degree $n_o$ in $\PP(A_{2})$. In particular this set is non-empty. 
\end{proposition}
\begin{proof}
By Proposition \ref{lemcubic} we have $\Hom_{X}(\Iscr,\Qscr(-1)) \neq 0$ if and only if $\det f = 0$. This is a homogeneous equation in $(\alpha,\beta,\gamma,\delta)$ of degree $n_o$ and we have to show it is not identically zero, i.e. we have to show there is at least one $\Qscr$ for which $\Hom_{X}(\Iscr,\Qscr(-1))=0$. This follows from Lemma \ref{atmost} and Lemma \ref{infinitecubic} below. 
\end{proof}
\begin{lemma} \label{atmost} 
Assume $k$ is uncountable and $\sigma$ has infinite order. Let $(n_e,n_o) \in N$ and $l \geq 0$ as in \eqref{uniquel}. Let $\Iscr \in \Rscr_{(n_e,n_o)}(X)$. Let $p,p' \in C$ such that $p \neq \sigma^m p'$ for all integers $m$. Modulo zero-dimensional objects, there exist at most $l$ different critical conic objects $\Qscr$ on $X$ for which $\Hom_{X}(\Iscr,\Qscr(-1)) \neq 0$ and $\Hom_{X}(\Qscr,\Nscr_p) \neq 0$, $\Hom_{X}(\Qscr,\Nscr_{p'}) \neq 0$. 
\end{lemma}
\begin{proof}
That $p \neq \sigma^m p'$ for all integers $m$ assures $\Hom_X(\Nscr_p,\Nscr_{p'}(m)) = 0$ and \linebreak $\Hom_X(\Nscr_{p'}(m),\Nscr_p) = 0$ for all integers $m$, which we will use throughout this proof.

We prove the statement by induction on $l$. First let $l = 0$ and assume by contradiction there is a non-zero map $f: \Iscr \ra \Qscr(-1)$. Let $\Iscr'(-2)$ be the kernel of $f$. By Lemma \ref{subobjectscritical} the image of $f$ is a shifted conic object $\Qscr'(-d)$ where $d \geq 1$. Using \eqref{shiftcubic} one computes 
$[\Iscr'] = [\Oscr] - 2(n_e-n_o)[\Sscr] + (n_e-n_o)[\Qscr]-(n_o-d)[\Pscr]$. It follows that $\Iscr'$ is a normalized line bundle 
on $X$ with invariants $(n_e-d,n_o-d)$. Since $(n_e-d,n_o-d) \not\in N$ this yields a contradiction with Proposition \ref{nonemptycubic}.

Let $l > 0$. Let $(\Qscr_i)_{i=1,\ldots,m}$ be different critical conic objects (modulo zero-dimensional objects) satisfying $\Hom_{X}(\Iscr,\Qscr_{i}(-1))\neq 0$ and $\Hom_{X}(\Qscr_i,\Nscr_p) \neq 0$, $\Hom_{X}(\Qscr_i,\Nscr_{p'}) \neq 0$. We will show $m\leq l$. If $m = 0$ then we are done. So assume $m > 0$. Let $\Qscr'_{i}(-1)$ be the kernel of a non-trivial map $\Qscr_{i} \ra \Nscr_p$. By Lemma \ref{subobjectscritical} $\Qscr'_{i}$ is a critical conic object, and we have an exact sequence
\[
0 \ra \Qscr'_{i}(-1) \ra \Qscr_{i} \ra \Nscr_p \ra 0.
\]
Applying $\Hom_{X}(-,\Nscr_{p'})$ we find $\Hom_{X}(\Qscr'_i(-2),\Nscr_{p'}(-1)) \neq 0$. Lemma \ref{incidence}(2) implies such a map factors through $\Qscr_{i}(-1)$. 
Let $\Qscr''_{i}(-3)$ be the kernel of a non-trivial map $\Qscr'_{i}(-2) \ra \Nscr_{p'}(-1)$. Again by Lemma \ref{subobjectscritical} $\Qscr''_{i}$ is a critical conic object, and 
\begin{equation} \label{sesconic''}
0 \ra \Qscr''_{i}(-3) \ra \Qscr'_{i}(-2) \xrightarrow{\pi} \Nscr_{p'}(-1) \ra 0.
\end{equation}
Applying $\Hom_{X}(-,\Nscr_p(-2))$ yields $\Hom_{X}(\Qscr''_i,\Nscr_p(1)) \neq 0$. Furthermore, as \eqref{sesconic''} is non-split, Serre duality (Theorem \ref{Serreduality}) implies $0 \neq \Ext^{1}_{X}(\Nscr_{p'}(-1),\Qscr''_i(-3)) \cong \Ext^{1}_X(\Qscr''_i,\Nscr_{p'}(-2))'$. By $\chi(\Qscr''_i,\Nscr_{p'}(-2)) = 0$ and again by Serre duality \linebreak $\Ext^{2}_X(\Qscr''_i,\Nscr_{p'}(-2)) \cong \Hom_X(\Nscr_{p'},\Qscr''_i(-1))' = 0$ we deduce $\Hom_{X}(\Qscr''_i,\Nscr_{p'}(-2)) \neq 0$.

Let $\Iscr'(-2)$ be the kernel of a non-trivial map $\iota:\Iscr \ra \Qscr_1(-1)$. As in first part of the proof one may show that $\Iscr'$ is a normalized line bundle on $X$ with invariants $(n_e - d,n_o - d)$ for some $d \geq 1$. 

Since $\Nscr_p(-1) = u_\ast \Oscr_{p'}$ for some point $p' \in C$ it follows by adjointness \linebreak $\dim_{k}\Hom_{X}(\Iscr,\Nscr_p(-1)) = \dim_{k}\Hom_{C}(u^{\ast}\Iscr,\Oscr_{p'})=1$. Hence for all $i$ the composition $a_i: \Iscr \ra \Qscr_i(-1) \ra \Nscr_p(-1)$ is a scalar multiple of $a_1$. Thus for all $i$ the map $a_i \circ \iota$ is a scalar multiple of $a_1 \circ \iota = 0$. Hence the composition $\Iscr'(-2) \ra \Iscr \ra \Qscr_i(-1)$ maps $\Iscr'(-2)$ to $\Qscr'_i(-2)$. 

As pointed out above the map $\pi$ in \eqref{sesconic''} factors through $\Qscr_i(-1)$. Thus the composition $\Iscr'(-2) \ra \Qscr'_i(-2) \ra \Nscr_{p'}(-1)$ is the same as the composition $b_i: \Iscr'(-2) \ra \Iscr \ra \Qscr_i(-1) \ra \Nscr_{p'}(-1)$. Same reasoning as above shows $b_i = 0$ for all $i$ hence the composition $\Iscr'(-2) \ra \Qscr'_i(-2)$ maps $\Iscr'(-2)$ to $\Qscr''_i(-3)$. 

We claim this map must be non-zero for $i > 1$. If not then there is a non-trivial map $\Iscr/\Iscr'(-2) \ra  \Qscr_i(-1)$ and since $\Iscr/\Iscr'(-2)$ is also a subobject of $\Qscr_1(-1)$ it follows that $\Qscr_{1}$ and $\Qscr_{i}$ have a common subobject. By Lemma \ref{equivalentsubobject} this contradicts the assumption $\Qscr_1$ and $\Qscr_{i}$ being different modulo zero dimensional objects. 

Hence $\Hom_{X}(\Iscr',\Qscr''_{i}(-1)) \neq 0$ for $i = 2,\ldots,m$. Since the $\Qscr''_i$ are still different modulo zero dimensional objects and $\Hom_{X}(\Qscr''_i,\Nscr_p(1)) \neq 0$, $\Hom_{X}(\Qscr''_i,\Nscr_{p'}(-2)) \neq 0$ we obtain by induction hypothesis $m-1 \leq l - d \leq l-1$ and hence $m \leq l$.
\end{proof}
The following lemma and its proof is inspired on the proof of \cite[Theorem 5.5.4]{DV1}.
\begin{lemma} \label{stablesubrepcubic}
Put $\theta = (-1,0,1)$. Let $V \in \mod(\Gamma^{0})$ and assume the forms $- \cdot \theta$ and $\chi(-,\underline{\dim}V)$ are proportional. Let $F \in \mod(\Gamma^{0})$ for which $F \perp V$. Then
\begin{enumerate}
\item
If $F' \subset F$ such that $\underline{\dim}F' \cdot \theta = 0$ then $F' \perp V$ and $F/F' \perp V$ 
\item
$\Hom_{\Gamma^{0}}(F,\Res p) = \Hom_{\Gamma^{0}}(\Res p, F) = 0$ for all $p \in C$ for which $\Res p$ is not perpendicular to $V$.
\end{enumerate}
\end{lemma}
\begin{proof}
Using the Euler form \eqref{ref-5.4-44cubic} on $\Gamma$ it is easy to see $\underline{\dim} V = (3l,2l,l)$ for some $l \in \NN$. We may assume $V \neq 0$. Let $F \perp V$. This implies that $F$ is $\theta$-semistable and $\underline{\dim}F = (n_o,n_e,n_o)$ for some $n_e,n_o \in \NN$. Since the result trivially holds for $F = 0$ we may assume $(n_e,n_o) \neq (0,0)$. 

To prove $(1)$, let $F' \subset F$ for which $\underline{\dim}F' \cdot \theta = 0$. Thus $\underline{\dim}F' = (m_o,m_e,m_o)$ and $\underline{\dim}F/F' = (n_o-m_o,n_e-m_e,n_o-m_o)$ for some $m_o \leq n_o$, $m_e \leq n_e$. Applying $\Hom_{\Gamma^{0}}(-,V)$ to $0 \r F' \r F \r F/F' \r 0$ gives the long exact sequence
\begin{align*}
0 & \r \Hom_{\Gamma^{0}}(F/F',V) \r \Hom_{\Gamma^{0}}(F,V) \r \Hom_{\Gamma^{0}}(F',V) \\
& \r \Ext^{1}_{\Gamma^{0}}(F/F',V) \r \Ext^{1}_{\Gamma^{0}}(F,V) \r \Ext^{1}_{\Gamma^{0}}(F',V) \r 0
\end{align*}
Since $F \perp V$ we deduce $\Hom_{\Gamma^{0}}(F/F',V) = 0$ and $\Ext^{1}_{\Gamma^{0}}(F',V) = 0$. Computations show $\chi(F',V) = \chi(F/F',V) = 0$ which yield $\Ext^{1}_{\Gamma^{0}}(F/F',V) = 0$ and $\Hom_{\Gamma^{0}}(F',V) = 0$. We conclude $F' \perp V$ and $F/F' \perp V$.

To prove $(2)$, let $p \in C$ for which $\Res p$ is not perpendicular to $V$. Recall from Lemma \ref{stablepointcubic} that $\Res p$ is $\theta$-stable. If by contradiction $\Hom_{\Gamma^{0}}(F,\Res p) \neq 0$ then by Lemma \ref{stablesurjectivecubic} there is an epimorphism $F \r \Res p$. By the first part we obtain $\Res p \perp V$, contradiction. 
Similarly $\Hom_{\Gamma^{0}}(\Res p,F) \neq 0$, finishing the proof.
\end{proof}
\begin{lemma} \label{ref-5.5.3-57cubic} 
Assume $\sigma$ has infinite order. Let $N \in \mod(\Gamma^{0})$ with dimension vector $(n_o,n_e,n_o)$ and assume $n_e \neq 0$. If $\Hom_{\Gamma^{0}}(N,\Res p)=\Hom_{\Gamma^{0}}(\Res p,N) = 0$ for all $p \in C$ then $\dim_{k}(\Ind N)_0 \leq n_e-1$. 
\end{lemma}
\begin{proof}
This is the same as the proof of \cite[Lemma 5.5.3]{DV1}. 
\end{proof}
We now come to the main result of this section. 
\begin{theorem} \label{allstablecubic}
Assume $k$ is uncountable. Let $A$ be an elliptic cubic Artin-Schelter algebra for which $\sigma$ has infinite order. Let $(n_e,n_o) \in N \setminus \{(0,0)\}$. If $M \in \Cscr_{(n_e,n_o)}(\Gamma)$ then $M^0$ is $\theta$-stable for $\theta = (-1,0,1)$. 
\end{theorem}
\begin{proof}
Let $\Iscr \in \Rscr_{(n_e,n_o)}(X)$ such that $M[-1]=\RHom_X(\Escr,\Iscr)$ and write $F = M^0$. By Propositions \ref{lemcubic} and \ref{curve} there exists a conic object $\Qscr$ on $X$ for which $F \perp Q^0$, where $Q[-1]=\RHom_X(\Escr,\Qscr(-1))$. This shows $F$ is $\theta$-semistable for $\theta=(-1,0,1)$. Hence there is a representation $F' \subsetneq F$ such that $F/F'$ is $\theta$-stable. We will prove $F'$ is necessarily zero, from which the theorem will follow. 

Assume by contradiction $F' \neq 0$. Since $F/F'$ is $\theta$-stable we have $\theta \cdot \underline{\dim}F/F' = 0$ thus we may put $\underline{\dim}F/F' = (n_o-m_o,n_e-m_e,n_o-m_o)$ for some $m_o \leq n_o$, $m_e \leq n_e$ for which $(m_e,m_o) \neq (n_e,n_o)$. 
Note $F'$ is $\theta$-semistable. We now claim that for all $p \in C$ we have
\begin{eqnarray} \label{claimstablethm}
\begin{aligned}
&\Hom_{\Gamma^0}(F',\Res p) = \Hom_{\Gamma^0}(\Res p,F') = 0, \\
&\Hom_{\Gamma^0}(F/F',\Res p) = \Hom_{\Gamma^0}(\Res p,F/F') = 0. 
\end{aligned}
\end{eqnarray}
To show \eqref{claimstablethm}, let $p \in C$. Lemmas \ref{infinitecubic} and \ref{atmost} imply there exists a conic object $\Qscr'$ on $X$ for which $\Hom_{X}(\Qscr',\Nscr_p(-3)) \neq 0$ and $\Hom_{X}(\Iscr,\Qscr'(-1)) = 0$. Writing $Q'[-1] = \RHom_{X}(\Escr,\Qscr'(-1))$ and using Theorem \ref{Serreduality} yields $\Ext^{1}_{\Gamma}(p,Q') \neq 0$ and $F \perp Q'$ by Proposition \ref{lemcubic}. In particular $\Res p$ is not perpendicular to $\Res Q'$. The claim \eqref{claimstablethm} now follows from Lemma \ref{stablesubrepcubic}. 

Combining \eqref{claimstablethm} with Lemma \ref{ref-5.5.3-57cubic} 
we find $\dim_{k} (\Ind F/F')_{0} \leq n_e-m_e-1$ and $\dim_{k}(\Ind F')_{0} \leq m_e -1$. Application of the right exact functor $\Ind$ on 
$0 \ra F' \ra F \ra F/F' \ra 0$ gives a long exact sequence in $\mod(\Gamma)$
\[
\dots \ra \Ind F' \ra \Ind F \ra \Ind F/F' \ra 0
\]
and counting dimensions yields
\begin{align*}
n_e-1 = \dim_{k}(\Ind F)_{0} & \leq \dim_{k}(\Ind F')_{0} + \dim_{k}(\Ind F/F')_{0} \\
& \leq (m_e-1) + (n_e-m_e-1) = n_e-2 
\end{align*}
which is absurd. Thus $F' = 0$ hence $F$ is $\theta$-stable, finishing the proof.  
\end{proof}
\begin{remark} \label{hypothesisk}
In case $A$ is of generic type A and $\sigma$ has infinite order we do not need the hypothesis $k$ is uncountable in Theorem \ref{allstablecubic}. This is because we may prove Proposition \ref{curve} without the additional hypothesis on $k$, using the proof of Lemma \ref{atmost} and Remark \ref{remarkuncountable}.
\end{remark}
\begin{remark} \label{hypothesiskenv}
In case $A = H_c$ is the enveloping algebra we may choose $\Qscr = \pi(H_c/zH_c)$ in the proof of Theorem \ref{allstablecubic}. See also Remark \ref{envelopingorthQ}. Again we do not need the hypothesis $k$ is uncountable. 
\end{remark}

\subsection{Second description of $\Rscr_{(n_e,n_o)}(X)$}

For $(n_e,n_o) \in N \setminus \{ (0,0) \}$ we denote $\Dscr_{(n_e,n_o)}(\Gamma^{0})$ for the image of $\Cscr_{(n_e,n_o)}(\Gamma)$ under the equivalence \eqref{equivalenceresubic}. 
\begin{theorem} \label{seconddescription}
Assume $k$ is uncountable. Let $A$ be an elliptic cubic AS-algebra for which $\sigma$ has infinite order. Let $\theta = (-1,0,1)$ and $(n_e,n_o) \in N \setminus \{ (0,0),(1,1) \}$. Then there is an equivalence of categories
\begin{diagram}[heads=LaTeX]
\Cscr_{(n_e,n_o)}(\Gamma) 
& \pile{
\rTo^{\Res} \\
\vspace{0.3cm} \\
\lTo_{\Ind}
}
& \Dscr_{(n_e,n_o)}(\Gamma^{0})
\end{diagram}
where 
\begin{multline*} 
\Dscr_{(n_e,n_o)}(\Gamma^{0}) = \{ F \in \mod (\Gamma^{0}) \mid \underline{\dim} F = (n_o,n_e,n_o), 
F \text{ is $\theta$-stable, } \\ \dim_{k}(\Ind F)_{0} \ge n_e-1 \}. 
\end{multline*} 
\end{theorem}
\begin{proof}
Below we often use the equivalence $\Cscr_{(n_e,n_0)}(\Gamma) \cong \Rscr_{(n_e,n_o)}(X)$ from Theorem \ref{firstdescrcubic}. We break the proof into five steps.
\setcounter{step}{0}
\begin{step} 
$\Res(\Cscr_{(n_e,n_o)}(\Gamma)) \subset \Dscr_{(n_e,n_o)}(\Gamma^{0})$. This follows from Theorem \ref{allstablecubic} and Lemma \ref{ref-5.5.3-57cubic}.
\end{step} 
\begin{step} \label{stap2} 
If $F \in \Dscr_{(n_e,n_o)}(\Gamma^{0})$ then $\Hom_{\Gamma^{0}}(F,\Res p) = \Hom_{\Gamma^{0}}(\Res p,F) = 0$ for all $p \in C$. Indeed, by Lemma \ref{stablepointcubic} and Lemma \ref{stablesurjectivecubic} any non-zero morphism would yield an isomorphism $F \cong \Res p$, contradicting the assumption $(n_e,n_o) \neq (1,1)$.
\end{step}
\begin{step} \label{stap3} 
If $F \in \Dscr_{(n_e,n_o)}(\Gamma^{0})$ then $\Hom_{\Gamma}(\Ind F,p) = \Hom_{\Gamma}(p,\Ind F) = 0$ for all $p \in E$. This follows from $0 = \Hom_{\Gamma^{0}}(F, \Res p) = \Hom_{\Gamma}(\Ind F, p)$ and 
\[
0 = \Hom_{\Gamma^{0}}(\Res p, F) = \Hom_{\Gamma^{0}}(\Res p, \Res \Ind F) = \Hom_{\Gamma}(\Ind \Res p, \Ind F)
\]
where we have used Step \ref{stap2} and $\Ind \Res p = p$ by Lemma \ref{pointperp}.
\end{step}
\begin{step}
$\Ind(\Dscr_{(n_e,n_o)}(\Gamma^{0})) \subset \Cscr_{(n_e,n_o)}(\Gamma)$. Let $F \in \Dscr_{(n_e,n_o)}(\Gamma^{0})$. Combining Step \ref{stap2} with Lemma \ref{ref-5.5.3-57cubic} gives $\underline{\dim}\Ind F = (n_o,n_e,n_o,n_e-1)$. Now Step \ref{stap3} shows $\Ind F \in \Cscr_{(n_e,n_o)}(\Gamma)$.
\end{step}
\begin{step} 
$\Ind$ and $\Res$ are inverses to each other. To prove this we only need to show $\Ind \Res F = F$ for $F \in \Cscr_{(n_e,n_o)}(\Gamma)$. This follows from Lemma \ref{ref-5.1.2-46cubic}.
\qed
\end{step} 
\def\qed{}
\end{proof}
For $(n_e,n_o) \in N \setminus \{(0,0),(1,1)\}$ let $\alpha = (n_o,n_e,n_o)$ and put 
\begin{equation} \label{ref-5.9-64cubic} 
\begin{split}
\widetilde{D}_{(n_e,n_o)} & = \{F\in \Rep_{\alpha}(\Gamma^0)\mid F \in 
\Dscr_{(n_e,n_o)}(\Gamma^{0})\} \\
& =\{F\in\Rep_{\alpha}(\Gamma^0)\mid F \text{ is $\theta$-stable, } \dim_k(\Ind 
F)_0\ge n_e-1\}. 
\end{split}
\end{equation} 
As $\widetilde{D}_{(n_e,n_o)}$ is a closed subset of the dense open subset of $\Rep_{\alpha}(\Gamma^0)$ consisting of all $\theta$-stable representations we obtain that $\widetilde{D}_{(n_e,n_o)}$ is locally closed. 

Denote $\Gl_{\alpha}(k) = \Gl_{n_e}(k) \times \Gl_{n_o}(k) \times \Gl_{n_e}(k)$. Put $D_{(n_e,n_o)} = \widetilde{D}_{(n_e,n_o)} \quot \Gl_{\alpha}(k)$. The next theorem provides the first part of Theorem \ref{theorem1} from the introduction. Our proof is, up to some minor computations, completely analogous to the proof of \cite[Theorem 5.5.5]{DV1}. For convenience of the reader we have included the proof. 
\begin{theorem} \label{ref-5.5.5-65cubic} 
Assume $k$ is uncountable. Let $A$ be an elliptic cubic AS-algebra for which $\sigma$ has infinite order. Then for $(n_e,n_o) \in N$ there exists a smooth locally closed variety $D_{(n_e,n_o)}$ of dimension $2\left(n_e - (n_e - n_o)^2\right)$ such that the isomorphism classes in $\Dscr_{(n_e,n_o)}(\Gamma^{0})$ (and hence in $\Rscr_{(n_e,n_o)}(X)$) are in natural bijection with the points in $D_{(n_e,n_o)}$. 
\end{theorem}
\begin{proof}
For $(n_e,n_o) = (0,0)$ or $(1,1)$ we refer to Corollaries \ref{invzero}, \ref{smallinvcubic} to see that we may take a point for $D_{(0,0)}$ and $(\PP^{1}\times\PP^1) - C$ for $D_{(1,1)}$. So we may assume $(n_e,n_o) \in N \setminus \{ (0,0),(1,1) \}$ throughout this proof. 

Since all representations in $\widetilde{D}_{(n_{e},n_{o})}$ are stable, all $\Gl(\alpha)$-orbits on $\widetilde{D}_{(n_{e},n_{o})}$ are closed and so $D_{(n_{e},n_{o})}$ is really the orbit space for the $\Gl(\alpha)$ action on $\widetilde{D}_{(n_{e},n_{o})}$. This proves that the isomorphism classes in $\Dscr_{(n_{e},n_{o})}(\Gamma)$ are in natural bijection with the points in $D_{(n_{e},n_{o})}$. 

To prove $D_{(n_{e},n_{o})}$ is smooth it suffices to prove $\widetilde{D}_{(n_{e},n_{o})}$ is smooth \cite{Luna}. We first estimate the dimension of $\widetilde{D}_{(n_{e},n_{o})}$. Write the equations of $A$ in the usual form \eqref{relations}. For $n_{e} \times n_{o}$-matrices $X$, $Y$ and $n_{o} \times n_{e}$-matrices $X'$, $Y'$ let $M_A^t(X',Y',X,Y)$ be obtained from $M_A^t$ by replacing $x^{2},xy,yx,y^{2}$ by $X'X$, $Y'X$, $X'Y$, $Y'Y$ (thus $M_A^t(X',Y',X,Y)$ is a $2n_{o}\times 2n_{o}$-matrix). Then $\widetilde{D}(n_{e},n_{o})$ has the following alternative description: 
\begin{multline*}
\widetilde{D}_{(n_{e},n_{o})} = \{ (X,Y,X',Y')\in M_{n_{e} \times n_{o}}(k)^2 \times M_{n_{o} \times n_{e}}(k)^2 \mid 
(X,Y,X',Y') \text{ is $\theta$-stable} \\
\text{ and } \rank M_{A}(X',Y',X,Y) \le 2n_{o}-(n_{e}-1) \}. 
\end{multline*} 
By Proposition \ref{nonemptycubic}, $\widetilde{D}_{(n_{e},n_{o})}$ is non-empty. The $(X,Y,X',Y')$ for which the associated representation is stable are a dense open subset of $M_{n_{e} \times n_{o}}(k)^2 \times M_{n_{o} \times n_{e}}(k)^2 $ and hence they represent a quasi-variety of dimension $4n_{e}n_{o}$. Imposing \linebreak $M_A(X',Y',X,Y)$ should have corank $\ge n_{e}-1$ represents 
$(2n_{o} - (2n_{o} - (n_{e}-1)))^{2} = (n_{e}-1)^{2}$ independent conditions. So the irreducible components of $\widetilde{D}_{(n_{e},n_{o})}$ have dimension $\ge 4n_{e}n_{o}-(n_{e}-1)^2$. Define $\widetilde{C}_{(n_{e},n_{o})}$ by 
\[ 
\{ G \in \Rep(\Gamma,\widetilde{\alpha}) \mid  G \cong \Ind \Res G, \Res G \in \widetilde{D}_{(n_{e},n_{o})} \} 
\] 
where $\widetilde{\alpha}=(n_{o},n_{e},n_{o},n_{e}-1)$ (as usual we assume the points of $\Rep(\Gamma,\widetilde{\alpha})$ to satisfy the relation imposed on $\Gamma$). To extend $F \in \widetilde{D}_{(n_{e},n_{o})}$ to a point in $\widetilde{C}_{(n_{e},n_{o})}$ we need to choose a basis in $(\Ind F)_0$. Thus $\widetilde{C}_{(n_{e},n_{o})}$ is a principal 
$\Gl_{n_{e}-1}(k)$ fiber bundle over $\widetilde{D}_{(n_{e},n_{o})}$. In particular $\widetilde{C}_{(n_{e},n_{o})}$ is smooth if and only $\widetilde{D}_{(n_{e},n_{o})}$ is smooth and the irreducible components of $\widetilde{C}_{(n_{e},n_{o})}$ have dimension $\ge 4n_{e}n_{o}-(n_{e}-1)^2 + (n_{e}-1)^2 = 4n_{e}n_{o}$. By the description of $\Cscr_{(n_{e},n_{o})}$ in Theorem \ref{firstdescrcubic} it follows that $\widetilde{C}_{(n_{e},n_{o})}$ is an open subset of $\Rep(\Gamma,\widetilde{\alpha})$. 

Let $x\in \widetilde{C}_{(n_{e},n_{o})}$. The stabilizer of $x$ consists of scalars thus if we put ${{G}}=\Gl(\widetilde{\alpha})/k^\ast$ then we have inclusions 
$\Lie(G) \subset T_{x}(\widetilde{C}_{(n_{e},n_{o})}) = T_{x}(\Rep(\Gamma,\widetilde{\alpha}))$. Next there is a natural inclusion $T_{x}(\Rep(\Gamma,\widetilde{\alpha}))/\Lie(G) \hookrightarrow \Ext^{1}_{\Gamma}(x,x)$.
Now $x$ corresponds to some normalized line bundle $\Iscr$ 
on $X$ and we have $\Ext^{1}_{\Gamma}(x,x) = \Ext^{1}_X(\Iscr,\Iscr)$. Lemma \ref{dimensionExt1} implies $\dim_{k}\Ext^{1}_{X}(\Iscr,\Iscr) = 2(n_{e} - (n_{e} - n_{o})^{2})$. Hence we obtain $4n_{e}n_{o} \le \dim T_{x}(\widetilde{C}_{(n_{e},n_{o})}) \le \dim_{k}\Ext^{1}_{\Gamma}(x,x) + \dim G$ and the right-hand is equal to 
$2(n_{e} - (n_{e} - n_{o})^{2}) + (n_{o}^{2} + n_{e}^{2} + n_{o}^{2} + (n_{e}-1)^{2} - 1) = 4n_{e}n_{o}$. Thus $\dim T_{x}(\widetilde{C}_n) = 4n_{e}n_{o}$ is constant and hence $\widetilde{C}_{(n_{e},n_{o})}$ is smooth. We also obtain $\dim \widetilde{D}_{(n_{e},n_{o})} = 4n_{e}n_{o} - (n_{e}-1)^2$. The dimension of $D(n_{e},n_{o})$ is equal to 
\begin{align*}
\dim \widetilde{D}_{(n_{e},n_{o})} - \dim \Gl(\alpha) + 1 & = 4n_{e}n_{o}-(n_{e}-1)^2 - (n_{o}^{2} + n_{e}^{2} + n_{o}^{2})+1 \\
& = 2(n_{e} - (n_{e} - n_{o})^{2}).
\end{align*}
This finishes the proof. 
\end{proof}

\subsection{
Descriptions of the varieties $D_{(n_e,n_o)}$ for the enveloping algebra} 

In case of the enveloping algebra we may further simplify the description of $D_{(n_e,n_o)}$. 
\begin{theorem} \label{thirddescrenv}
Let $A = H_c$ be the enveloping algebra. Let $(n_e,n_o) \in N$. The isomorphism classes in $R_{(n_e,n_o)}(A)$ are in natural bijection with the points in the smooth affine variety $D_{(n_e,n_o)}$ of dimension $2(n_e - (n_e - n_o)^2)$ where
\begin{multline*}
D_{(n_e,n_o)} = \{ (X,Y,X',Y')\in M_{n_{e} \times n_{o}}(k)^2 \times M_{n_{o} \times n_{e}}(k)^2 \mid Y'X - X'Y \text{ isomorphism,}   \\ 
\rank 
\left(
\begin{array}{ccc}
Y'Y &  X'Y - 2Y'X\\  
Y'X - 2X'Y & X'X
\end{array}
\right) 
\leq 2n_o-(n_e-1) \}/\Gl_{\alpha}(k)
\end{multline*}
\end{theorem}
\begin{proof}
For $(n_e,n_o) = (0,0)$ or $(1,1)$ we refer to Corollaries \ref{invzero}, \ref{smallinvcubic} to see that $D_{(n_e,n_o)}$ has the description as in the statement of the current theorem. So we may assume $(n_e,n_o) \in N \setminus \{ (0,0),(1,1) \}$ throughout this proof. 

Consider the conic object $\Qscr = \pi(A/zA)$ on $X$ where $z = xy - yx$. Write $Q = \Ext^{1}_{X}(\Escr,\Qscr(-1)) \in \mod(\Gamma)$ (Lemma \ref{conicreps}) and put $V = \Res Q$. 

It is sufficient to show $\Dscr_{(n_e,n_o)}(\Gamma^{0})$ has the alternative description
\begin{multline*}
\Dscr'_{(n_e,n_o)}(\Gamma^{0}) := \{F \in \mod (\Gamma^{0}) \mid \underline{\dim} F=(n_o,n_e,n_o), \\
F \perp V, \dim_{k}(\Ind F)_{0} \ge n_e-1 \}.
\end{multline*}
Indeed, if $\Dscr_{(n_e,n_o)}(\Gamma^{0}) = \Dscr'_{(n_e,n_o)}(\Gamma^{0})$ we then have 
\begin{equation*} 
\begin{split}
\widetilde{D}_{(n_e,n_o)} & = \{F\in \Rep_{\alpha}(\Gamma^0)\mid F\in 
\Dscr_{(n_e,n_o)}(\Gamma^{0})\} \\
& =\{F\in\Rep_{\alpha}(\Gamma^0)\mid \phi_V(F)\neq 0, \dim_k(\Ind 
F)_0\ge n_e-1\}
\end{split}
\end{equation*} 
from which is clear that $\widetilde{D}_{(n_e,n_o)}$ is a closed subset of $\{\phi_V\neq 0\}$ so in particular $\widetilde{D}_{(n_e,n_o)}$ is affine. This means $D_{(n_e,n_o)} = \widetilde{D}_{(n_e,n_o)} / \Gl(\alpha)$ is an affine variety. Theorem \ref{ref-5.5.5-65cubic} further implies $D_{(n_e,n_o)}$ is smooth of dimension $2(n_e - (n_e - n_o)^2)$ which points are in natural bijection with the isomorphism classes in $\Rscr_{(n_e,n_o)}(X)$ whence in $R_{(n_e,n_o)}(A)$ by \S\ref{Normalized line bundles}. Moreover, as in the proof of Theorem \ref{ref-5.5.5-65cubic}, $\widetilde{D}_{(n_e,n_o)}$ has the alternative description
\begin{multline*}
\widetilde{D}_{(n_{e},n_{o})} = \{ F = ((X,Y),(X',Y')) \in \Rep_{(n_o,n_e,n_o)}(\Gamma^{0}) \mid F \perp V \\
\text{ and } \rank M_{A}(X',Y',X,Y) \le 2n_{o}-(n_{e}-1) \}.
\end{multline*} 
As shown in Proposition \ref{lemcubic} the condition $(X,Y,X',Y') \perp V$ is equivalent with saying $Y'X-X'Y$ is an isomorphism. Explicitely writing down $M_A$ by \eqref{relations}, \eqref{relenvweyl} yields the desired description of $D_{(n_e,n_o)}$.

So to prove the current theorem it remains to prove $\Dscr_{(n_e,n_o)}(\Gamma^{0}) = \Dscr'_{(n_e,n_o)}(\Gamma^{0})$. We will do this by showing that the functors $\Res$ and $\Ind$ define inverse equivalences between $\Cscr_{(n_e,n_o)}(\Gamma)$ and $\Dscr'_{(n_e,n_o)}(\Gamma^{0})$.
\setcounter{step}{0}
\begin{step} \label{ref-1-60}  
$\Res(\Cscr_{(n_e,n_o)}(\Gamma)) \subset \Dscr'_{(n_e,n_o)}(\Gamma^{0})$. Let $M \in \Cscr_{(n_e,n_o)}(\Gamma)$. That $\Res M \perp V$ follows from Proposition \ref{lemcubic} and Remark \ref{envelopingorthQ}. Further, since $\Ind \Res M = M$ by Lemma \ref{ref-5.1.2-46cubic} we find $\dim_k(\Ind \Res M)_0 = n_e-1$ by Theorem \ref{firstdescrcubic}.
\end{step} 
\begin{step}
If $p \in C$ then $\Res p \in \mod(\Gamma^{0})$ is not perpendicular to $V$. Indeed, by the equivalences \eqref{equivalenceresubic} and \eqref{Bondalcubic} one finds $\RHom_{\Gamma^{0}}(\Res p,V) = \RHom_{\Gamma}(p,Q) = \RHom_{X}(\Nscr_p,\Qscr[1])$. Thus $\Ext^{1}_{\Gamma_{0}}(\Res p,V) = \Ext^{2}_{X}(\Nscr_p,\Qscr)$. Further, Serre duality (Theorem \ref{Serreduality}) yields $\Ext^{2}_{X}(\Nscr_p,\Qscr) \cong \Hom_{X}(\Qscr,\Nscr_{\sigma^{4}p})' = \Hom_{X}(\Qscr,\Nscr_{p})'$, which is non-zero by Lemma \ref{incidence}(1). This proves Step 2.
\end{step}
\begin{step} \label{ref-4-63} 
$\Ind(\Dscr'_{(n_e,n_o)}(\Gamma^{0})) \subset \Cscr_{(n_e,n_o)}(\Gamma)$. Let $F \in \Dscr'_{(n_e,n_o)}(\Gamma^{0})$. Combining Step 2 with Lemmas \ref{ref-5.5.1-55cubic}, \ref{stablesubrepcubic} and \ref{ref-5.5.3-57cubic} we obtain $\dim_k(\Ind F)_{0} = n_e-1$. It remains to show $\Hom_{\Gamma}(\Ind F,p) = \Hom_{\Gamma}(p,\Ind F) = 0$ for $p \in C$. By Lemma \ref{stablepointcubic} we have $p = \Ind \Res p$. Thus  $\Hom_{\Gamma}(\Ind F,p) = \Hom_{\Gamma^{0}}(F,\Res p)=0$ and similarly
\[ 
\Hom_{\Gamma}(p,\Ind F) = \Hom_{\Gamma^{0}}(\Res p,\Res \Ind F) = \Hom_{\Gamma^{0}}(\Res p, F) = 0
\]
where we have used Lemma \ref{stablesubrepcubic} again. 
\end{step} 
\begin{step}  
$\Ind$ and $\Res$ are inverses to each other. 
This follows from Lemma \ref{ref-5.1.2-46cubic}. 
\qed
\end{step} 
\def\qed{} 
\end{proof}
We further simplify the description of $D_{(n_e,n_o)}$ as
\begin{theorem} \label{thirddescrenv'}
Let $A = H_c$ be the enveloping algebra. Let $(n_e,n_o) \in N$. The isomorphism classes in $R_{(n_e,n_o)}(A)$ are in natural bijection with the points in the smooth affine variety $D_{(n_e,n_o)}$ of dimension $2(n_e - (n_e - n_o)^2)$ where
\begin{multline*}
D_{(n_e,n_o)} = \{ (\XX,\YY,\XX',\YY') \in M_{n_{e} \times n_{o}}(k)^2 \times M_{n_{o} \times n_{e}}(k)^2 \mid 
\YY'\XX - \XX'\YY = \II \text{ and }   \\ 
\rank (\YY\XX'-\XX\YY'-\II)\leq 1 \}/\Gl_{n_e}(k) \times \Gl_{n_o}(k) 
\end{multline*}
\end{theorem}
\begin{proof}
For $(n_e,n_o) = (0,0)$ or $(1,1)$ we refer to Corollaries \ref{invzero}, \ref{smallinvcubic} to see that $D_{(n_e,n_o)}$ has the description as in the statement of the current theorem. So we may assume $(n_e,n_o) \in N \setminus \{ (0,0),(1,1) \}$. 

Similarly as in Theorem \ref{firstdescrHeis} we define for any $F \in \mod(\Gamma^{0})$ the linear map 
\[
F(Z_{-3}) = F(Y_{-2})F(X_{-3}) - F(X_{-2})F(Y_{-3})
\]
In order to prove the current theorem, it is sufficient to show  
$\Dscr_{(n_e,n_o)}(\Gamma^{0})$ has the alternative description
\begin{multline*}
\Dscr''_{(n_e,n_o)}(\Gamma^{0}) := \{F \in \mod (\Gamma^{0}) \mid \underline{\dim} F=(n_o,n_e,n_o), 
F(Z_{-3}) \text{ isomorphism}, \\
\rank(F(Y_{-3})F(Z_{-3})^{-1}F(X_{-2}) - F(X_{-3})F(Z_{-3})^{-1}F(Y_{-2}) - \Id) \leq 1 \}.
\end{multline*}
Indeed, if $\Dscr_{(n_e,n_o)}(\Gamma^{0}) = \Dscr''_{(n_e,n_o)}(\Gamma^{0})$ then 
\begin{align*}
\widetilde{D}_{(n_e,n_o)} & = \{F \in \Rep_{\alpha}(\Gamma^0) \mid F \in \Dscr_{(n_e,n_o)}(\Gamma^{0})\} \\
& = \{ (X,Y,X',Y') \in \Rep_{\alpha}(\Gamma^0) \mid Z:=Y'X - X'Y \text{ isomorphism}, \\
& \hspace{5cm} \rank (YZ^{-1}X'-XZ^{-1}Y'-\Id)\leq 1  \}
\end{align*}
and by $D_{(n_e,n_o)} = \widetilde{D}_{(n_e,n_o)}/\Gl_{\alpha}(k)$
the statement of the current theorem will follow.

What remains to prove is $\Dscr_{(n_e,n_o)}(\Gamma^{0}) = \Dscr''_{(n_e,n_o)}(\Gamma^{0})$. We will do this by showing that the functors $\Res$ and $\Ind$ define inverse equivalences between $\Cscr_{(n_e,n_o)}(\Gamma)$ and $\Dscr''_{(n_e,n_o)}(\Gamma^{0})$. This is done in the following three steps.
\setcounter{step}{0}
\begin{step} 
$\Res(\Cscr_{(n_e,n_o)}(\Gamma)) \subset \Dscr''_{(n_e,n_o)}(\Gamma^{0})$. Let $M \in \Cscr_{(n_e,n_o)}(\Gamma)$ and put $F = \Res M$. For convenience we denote 
$X = M(X_{-3})$, $X' = M(X_{-2})$, $X'' = M(X_{-1})$ (similarly for $Y$). Theorem \ref{firstdescrHeis} already implies $Z=Y'X - X'Y$ is an isomorphism and $Z'=Y''X' - X''Y'$ is surjective. Thus to show Step 1, what remains to prove is $\rank(YZ^{-1}X' - XZ^{-1}Y' - \Id) \leq 1$. From \eqref{relationsGammacubic} we deduce
\begin{equation} \label{relsM}
\begin{aligned} 
Y''Y'X - 2Y''X'Y + X''Y'Y & = 0 \\
X''X'Y - 2X''Y'X + Y''X'X & = 0
\end{aligned}
\end{equation}
and these equations may be written as $X''Z = Z'X$, $Y''Z = Z'Y$. Since $Z$ is an isomorphism we find $X'' = Z'XZ^{-1}$ and $Y'' = Z'YZ^{-1}$. Substitution yields
\[
Z' = Y''X' - X''Y' = Z'(YZ^{-1}X'-XZ^{-1}Y')
\]
thus $Z'(YZ^{-1}X'-XZ^{-1}Y' - \Id) = 0$. As $Z'$ is surjective, it has a one dimensional kernel, completing proof of Step 1.
\end{step}
\begin{step}
$\Ind(\Dscr''_{(n_e,n_o)}(\Gamma^{0})) \subset \Cscr_{(n_e,n_o)}(\Gamma)$. To prove so, let $F \in \Dscr''_{(n_e,n_o)}(\Gamma^{0})$. We will construct a representation $M$ of $\Gamma$ for which $M \in \Cscr_{(n_e,n_o)}(\Gamma)$ and $\Res M = F$. For then, $F \in \Dscr_{(n_e,n_o)}(\Gamma^{0})$ by Theorem \ref{seconddescription} hence $\Ind F = M \in \Cscr_{(n_e,n_o)}(\Gamma)$. 

For simplicity we denote $X = F(X_{-3})$, $X' = F(X_{-2})$ (similarly for $Y$). Put $Z = Y'X - X'Y$. Let $Z'$ denote the projection $F_{-2} \r F_{-2}/\im (YZ^{-1}X'-XZ^{-1}Y' - \Id)$. Define the linear maps $X'' = Z'XZ^{-1}$, $Y''$ = $Z'YZ^{-1}$. We now define $M$ as
\begin{equation*}
\begin{diagram}[heads=LaTeX]
F_{-3}
& \pile{
\rTo^{X} \\
\rTo^{Y} 
}
& F_{-2} & \pile{
\rTo^{X'} \\
\rTo^{Y'} 
}
& F_{-1} & \pile{
\rTo^{X''} \\
\rTo^{Y''} 
}
& \im Z'
\end{diagram}
\end{equation*}
\end{step}
In fact $\dim_{k} M_0 = n_e - 1$, as otherwise $\id = YZ^{-1}X'-XZ^{-1}Y'$ and by taking traces we find 
$n_e = \tr(YZ^{-1}X'-XZ^{-1}Y') = \tr(-Z^{-1}(Y'X - X'Y)) = \tr(-Z^{-1}Z) = -n_e$
whence $n_e = 0$. By definition \eqref{defN} of $N$ this leads to $n_o = 0$, contradiction the assumption $(n_e,n_o) \neq (0,0)$. 

To prove $M \in \mod(\Gamma)$ we need to check the relations \eqref{relsM}. This is easy to do. One also checks $Z' = Y''X' - X''Y'$. Now Theorem \ref{firstdescrHeis} implies $M \in \Cscr_{(n_e,n_o)}(\Gamma)$. 

By the construction of $M$ we have $\Res M = F$. This proves Step 2.

\begin{step}  
$\Ind$ and $\Res$ are inverses to each other. 
This follows from Lemma \ref{ref-5.1.2-46cubic}. 
\qed
\end{step} 
\def\qed{} 
\end{proof}

\subsection{Description of the varieties $D_{(n_e,n_o)}$ for generic type A} 

In Theorem \ref{thirddescrenv} we have simplified the description of the varieties $D_{(n_e,n_o)}$ for the enveloping algebra. That such a simplification is possible is due to the fact that there exists a conic object $\Qscr$ on $X$ for which $M^{0} \perp Q^0$ for all $M \in \Cscr_{(n_e,n_o)}(\Gamma)$. 

It is therefore natural to ask if, for general cubic AS-algebras $A$, there exists a conic object $\Qscr$ on $X$ for which $M^{0} \perp Q^0$ for all $M \in \Cscr_{(n_e,n_o)}(\Gamma)$ where $\Qscr$ is independent of $M$. This is unknown (and probably unlikely). However there is another interpretation. In case of the enveloping algebra we relied on fact $u^{\ast}\Mscr \cong \Oscr_{\Delta}$ for all normalized line bundles $\Mscr$ on $X$. For generic $A$ we have
\begin{lemma}  \label{ref-5.5.1-55cubic} 
Let $A$ be a cubic AS-algebra of generic type A for which $\sigma$ has infinite order. There exists $V \in \mod(\Gamma^{0})$ with $\underline{\dim} V = (6,4,2)$ for which 
\begin{enumerate} 
\item 
for all $M \in \Cscr_{(n_e,n_o)}(\Gamma)$ we have $M^{0} \perp V$, and
\item 
if $p \in C$ then 
$\Res p$ is not perpendicular to $V$.
\end{enumerate} 
\end{lemma} 
\begin{proof}
Obtained by repeating the arguments in \cite[Lemma 5.5.1]{DV1} where, in the current setting, we pick a degree zero line bundle $\Uscr$ on $C$ which is not of the form $\Oscr((o)-(2(n_e+n_o)\xi))$ for $n_e,n_o \in \NN$, see Proposition \ref{ref-4.3-39cubic}.
\end{proof}
\begin{theorem} \label{ref-5.5.4-59cubic}
Let $A$ be a cubic AS-algebra of generic type A for which $\sigma$ has infinite order. Let $V \in \mod(\Gamma^0)$ be as in Lemma \ref{ref-5.5.1-55cubic}. Let $(n_e,n_o) \in N$. The isomorphism classes in $R{(n_e,n_o)}(A)$ are in natural bijection with the points in the smooth affine variety $D_{(n_e,n_o)}$ of dimension $2(n_e - (n_e - n_o)^2)$ where
\begin{multline*}
D_{(n_e,n_o)} = \{ F = (X,Y,X',Y') \in M_{n_{e} \times n_{o}}(k)^2 \times M_{n_{o} \times n_{e}}(k)^2 \mid F \perp V,  \\ 
\rank 
\left(
\begin{array}{ccc}
aY'Y + cX'X &  bX'Y+aY'X\\
bY'X + aX'Y & aX'X+cY'Y
\end{array}
\right) 
\leq 2n_o-(n_e-1) \}/\Gl_{\alpha}(k)
\end{multline*}
\end{theorem}
\begin{proof}
Analogous to the proof of Theorem \ref{thirddescrenv} where one uses Lemma \ref{ref-5.5.1-55cubic}.
\end{proof}

\section{Invariant ring of the first Weyl algebra and proof of Theorem \ref{theorem3}} \label{Ideals of an invariantring of the first Weyl algebra}

In this section we show how application of the previous results for the enveloping algebra $H_c$ may be used to classify the right ideals of an invariant ring of the first Weyl algebra.

Let $A = H_c$ denote the enveloping algebra and write $z = xy - yx \in A_2$. In the notations of \S\ref{Geometric data associated to a cubic AS-algebra} the canonical normalizing element $g$ is given by $z^2 \in A_4$ and $h = z$ is central. Consider the graded algebra $\Lambda = A[h^{-1}]$, the localisation of $A$ at the powers of $h = z$, and its subalgebra $\Lambda_0$ of elements of degree zero. It is shown in \cite[Theorem 8.20]{ATV2} that $\Lambda_0 = A_1^{\langle \varphi \rangle}$, the algebra of invariants of the first Weyl algebra $A_1 = k<x,y>/(xy-yx-1)$ under the automorphism $\varphi$ defined by $\varphi(x)=-x$, $\varphi(y)=-y$. 

For any positive integer $l$, let $V_l$ be the $k$-linear space spanned by the set $\{x^iy^j \mid i+j \text{ even and } i+j \leq 2l  \}$. Then $k = V_0 \subset V_1 \subset \dots$ endow $A_1^{\langle \varphi \rangle}$ with a positive filtration. 

The associated Rees ring $\Rees(A_1^{\langle \varphi \rangle})= \bigoplus_{l \in \NN}V_l$ is identified with the subring $\bigoplus_{l \in \mathbb{N}}V_lt^l$ of the ring of Laurent polynomials $A_1^{\langle \varphi \rangle}[t,t^{-1}]$, and identifying $t = h$ we see the Rees ring of $A_1^{\langle \varphi \rangle}$ is isomorphic to $A^{(2)}$, the $2$-Veronese of $A$. The associated graded algebra $\gr(A_1^{\langle \varphi \rangle}) = \bigoplus_{l \in \NN}V_l/V_{l-1}$ is isomorphic to $A^{(2)}/hA^{(2)} = k[x,y]^{(2)}$, the $2$-Veronese of the commutative polynomial ring $k[x,y]$. 

Similarly, for a filtered $A_1^{\langle \varphi \rangle}$-module $M$ we write $\Rees(M) = \bigoplus_{l \in \NN}M_l$, which is isomorphic to an object in $\GrMod(A^{(2)})$ and $\gr(M) = \bigoplus_{l \in \NN}M_l/M_{l-1}$ for the associated graded module, identified with an object of $\Mod(k[x,y]^{(2)})$.

Write $\Filt(A_1^{\langle \varphi \rangle})$ for the category which objects are the filtered right $A_1^{\langle \varphi \rangle}$-modules and morphisms the $A_1^{\langle \varphi \rangle}$-morphisms $f: M \r N$ which are strict i.e. $N_n \cap \im(f) = f(M_n)$ for all $n$.
Write $\GrMod(A^{(2)})_{h}$ for the full subcategory of $\GrMod(A^{(2)})$ consisting of the $h$-torsion free modules. The exact functor
\[
\Rees(-): \Filt(A_1^{\langle \varphi \rangle}) \r \GrMod(A^{(2)})_{h}
\]
is an equivalence and $(\Rees(M)[h^{-1}])_0 \cong M$ for all $M \in \Filt(A_1^{\langle \varphi \rangle})$. 

Let $R(A_1^{\langle \varphi \rangle})$ denote the set of isomorphism classes of right $A_1^{\langle \varphi \rangle}$-ideals. Note 
\[
R(A_1^{\langle \varphi \rangle}) = \{ M \in \mod(A_{1}^{\langle \varphi \rangle}) \mid  M \text{ torsion free of rank one} \} / \text{ iso }
\]
Performing a similar treatment as in \cite[\S 4]{BW} yields
\begin{proposition} \label{ideauxinvariants1}
The set $R(A_1^{\langle \varphi \rangle})$ is in natural bijection with the isomorphism classes in the full subcategory of $\coh(\Delta)$ with objects 
\begin{equation} \label{setenv}
\{\Mscr \in \coh(X) \mid u^{\ast}\Mscr \cong \Oscr_{\Delta} \}
\end{equation}
\end{proposition}
\begin{proof}
We will make use of the following commutative diagram
\begin{equation} \label{diagramenv}
\begin{diagram}[heads=LaTeX]
\filt(A_1^{\langle \varphi \rangle})
& 
\rTo^{\Rees} 
&  
\grmod(A^{(2)})
&
\rTo^{\pi}
&
\tails(A^{(2)})
&
\rTo^{\cong}
&
\tails(A) \\
\dTo^{\cong}
&&&&&&
\dTo_{u^{\ast}} \\
\filt(A_1^{\langle \varphi \rangle})
& 
\rTo^{\gr} 
&  
\grmod(k[x,y]^{(2)})
&
\rTo^{\pi}
&
\tails(k[x,y]^{(2)})
&
\rTo^{\cong}
&
\tails(k[x,y]) 
\end{diagram}
\end{equation}
Let $\Mscr \in \coh(X)$ for which $u^{\ast}\Mscr \cong \Oscr_{\Delta}$. It follows from Remark \ref{RemarkHeis} and \S\ref{Normalized line bundles} that $\Mscr = \pi M$ for some $M \in \grmod(A)$ which is torsion free of rank one (and therefore critical). Thus $M[h^{-1}]_0$ is a critical rank one object i.e. $M[h^{-1}]_0 \in R(A_{1}^{\langle \varphi \rangle})$. To show this correspondence $\Mscr \mapsto M[h^{-1}]_0$ is bijective it suffices to give its inverse.

Let $M \in \mod(A_{1}^{\langle \varphi \rangle})$ be torsion free of rank one and let us fix, temporarily, an embedding of $M$ as an ideal of $A_{1}^{\langle \varphi \rangle}$. The filtration $V_l$ on $A_1^{\langle \varphi \rangle}$ induces a filtration $M_l = M \cap V_l$ on $M$.
Let us still write $M$ for the associated object in $\filt(A_{1}^{\langle \varphi \rangle})$. 

Arguing on the bottom half of \eqref{diagramenv} it follows that $\gr(M) \subset k[x,y]^{(2)}$ is a homogeneous ideal and writing $\overline{\Mscr} \in \tails(k[x,y]) = \coh(\Delta)$ for the image of $\pi\gr(M)$ we deduce $\overline{\Mscr} \subset \Oscr_{\Delta}$. As the quotient $\Oscr_{\Delta}/\overline{\Mscr}$ has rank zero it is finite dimensional i.e. of the form $\Oscr_{D}$ for some divisor $D$ on $\Delta$ of degree $d = \deg D \geq 0$. Thus $\overline{\Mscr} \cong \Oscr_{\Delta}(-D)$. Redefine the filtration on $M$ by $M_{l}^{0} := M_{l+d}$ and write $M^{0}$ for the associated object in $\filt(A_{1}^{\langle \varphi \rangle})$. Repeating the arguments we now find $\overline{\Mscr^{0}} \cong \Oscr_{\Delta}$. By the commutativity of the diagram \eqref{diagramenv}, $M^{0}$ now corresponds to an object $\Mscr^{0} \in \tails(A)$ for which $u^{\ast} \Mscr^{0} \cong \Oscr_{\Delta}$. This finishes the proof.
\end{proof}
We may now prove Theorem \ref{theorem3}.
\begin{proof}[Proof of Theorem \ref{theorem3}]
Combine Theorem \ref{thirddescrenv'}, Proposition \ref{ideauxinvariants1} and Remark \ref{RemarkHeis}. 
\end{proof}

\section{Filtrations of line bundles and proof of Theorem \ref{theorem4}} \label{filtrationscubic}

Let $A$ be an elliptic cubic AS-algebra for which $\sigma$ has infinite order. The following analogue of \cite[Lemma 5.6.5]{DV1} shows how to reduce the invariants of a line bundle. 
\begin{lemma} \label{ref-5.6.5-70cubic} 
Assume $k$ is uncountable and $\sigma$ has infinite order. Let $(n_e,n_o) \in N$ such that $(n_e-1,n_o-1) \in N$. Let $\Iscr \in \Rscr_{(n_e,n_o)}(X)$. Then there exists a conic object $\Qscr$ on $X$ for which $\Ext^{1}_{X}(\Qscr(1),\Iscr(-2))\neq 0$. If $\Jscr = \pi J$ is the middle term of a corresponding non-trivial extension and $\Jscr^{\ast\ast} = \pi J^{\ast\ast}$ then $\Jscr^{\ast\ast} \in \Rscr_{(m_e,m_o)}(X)$ with $m_e < n_e$, $m_o < n_o$. Furthermore $\Jscr^{\ast\ast}/\Iscr(-2)$ is a shifted conic object on $X$. 
\end{lemma} 
\begin{proof} 
We have $\Ext^{1}_{X}(\Qscr(1),\Iscr(-2)) \cong \Ext^{1}_{X}(\Iscr(-2),\Qscr(-3))'= \Ext^{1}_{X}(\Iscr,\Qscr(-1))'$ and $\Ext^{2}_{X}(\Iscr,\Qscr(-1)) \cong \Hom_{X}(\Qscr(-1),\Iscr(-4))' = 0$ by Theorem \ref{Serreduality} (Serre duality). Thus $\chi(\Iscr,\Qscr(-1)) = 0$ shows $\dim_{k} \Hom_{X}(\Iscr,\Qscr(-1)) = \dim_{k}\Ext^{1}_{X}(\Iscr,\Qscr(-1))$. Hence it follows from Proposition \ref{curve} there exist a conic object $\Qscr$ for which $\Ext^1_{X}(\Qscr(1),\Iscr(-2))\neq 0$. 

Let $\Jscr = \pi J$ be the middle term of a non-trivial extension of $\Iscr(-2)$ by $\Qscr(1)$. It is easy to see $\Jscr$ is torsion free, see for example \cite[Lemma 5.6.5]{DV1}. A computation yields $[\Iscr(-2)] = [\Oscr]- 2(n_e - n_o)[\Sscr]+ (n_e - n_o - 1)[\Qscr] - n_o[\Pscr]$
hence $[\Jscr] = [\Oscr]- 2(n_e - n_o)[\Sscr]+ (n_e - n_o)[\Qscr] - (n_o-1)[\Pscr]$. Thus $\Jscr$ is normalized with invariants $(n_e-1,n_o-1)$.

By \cite[Corollary 4.2]{ATV2} $\gkdim J^{\ast\ast}/J \le 1$. As $\sigma$ has infinite order, any zero dimensional object on $X$ admits a filtration by shifted point objects \cite{ATV2}. Hence $[\Jscr^{\ast\ast}/\Jscr] = c[\Pscr]$ for some $c \geq 0$ and therefore $[\Jscr^{\ast \ast}] = [\Oscr]- 2(n_e - n_o)[\Sscr]+ (n_e - n_o)[\Qscr] - (n_o - c -1)[\Pscr]$. Thus $\Jscr^{\ast \ast} \in \Rscr_{(m_e,m_o)}(X)$ where $m_o = n_o - c - 1 < n_o$ and $m_e = n_e - c - 1$. Let $\Nscr = \Jscr^{\ast\ast}/\Iscr(-2)$. Then $\Nscr$ is pure and furthermore we have $[\Nscr] = [\Qscr] + (c+1)[\Pscr]$. Thus $e(\Nscr) = 1$. Moreover $\Qscr \subset \Nscr$ and $\Nscr / \Qscr$ is zero dimensional. By Lemma \ref{extentionsconic} $\Nscr$ is a shifted conic object on $X$.
\end{proof} 
\begin{theorem} \label{ref-5.6.6-72cubic} 
Assume $k$ is uncountable. Let $A$ be an elliptic cubic AS-algebra and assume $\sigma$ has infinite order. Let $(n_e,n_o) \in N$ and $l$ as in \eqref{uniquel}. Let $\Iscr \in \Rscr_{(n_e,n_o)}(X)$. Then there exists an integer $m$, $0 \leq m \leq l$ together with a filtration of line bundles $\Iscr_0 \supset \Iscr_1\supset \cdots \supset \Iscr_m = \Iscr(-2l)$ on $X$ with the property that the $\Iscr_i/\Iscr_{i+1}$ are shifted conic objects on $X$ and $\Iscr_0$ has invariants $(m_e,m_o) = (n_e - l,n_o - l)$.
\end{theorem}
\begin{proof}
By Lemma \ref{ref-5.6.5-70cubic}, \eqref{uniquel} and downwards induction on $l$. 
\end{proof} 
We end with the proof of Theorem \ref{theorem4} from the introduction.
\begin{proof}[Proof of Theorem \ref{theorem4}]
The first part of Theorem \ref{theorem4} is due to Theorem \ref{ref-5.6.6-72cubic} and the equivalence $R(A) = \coprod_{(n_e,n_o) \in N}\Rscr_{(n_e,n_o)}$ from \S\ref{Normalized line bundles}. Proposition \ref{minimalresunique} implies $I_0 = \omega \Iscr_0$ having a minimal resolution of the form \eqref{commall}. By Proposition \ref{minimalresunique} and Remark \ref{minimalresuniqueremark}, $I_0$ and $\Iscr_0$ are up to isomorphism uniquely determined by $(m_e,m_o)$, and therefore by $(n_e,n_o)$. This finishes the proof.
\end{proof}
\begin{remark}
By Remarks \ref{hypothesisk} and \ref{hypothesiskenv} we do not need the hypothesis $k$ is uncountable in Theorem \ref{ref-5.6.6-72cubic} in case $A = H_c$ is the enveloping algebra or $A$ is of generic type A and $\sigma$ has infinite order. Furthermore it follows from Proposition \ref{commfreecubic} that Theorem \ref{ref-5.6.6-72cubic} is (trivially) true in case $A$ is a linear cubic AS-algebra, again without the hypothesis $k$ is uncountable.
\end{remark}

\appendix 
\section{Hilbert series of ideals up to invariants $(3,3)$} \label{ref-B-67} 
Let $A$ be a cubic Artin-Schelter algebra, and let $I$ be a normalized rank one torsion free graded right $A$-module of projective dimension one with invariants $(n_e,n_o)$. According to Theorem \ref{theorem2} the Hilbert series of $I$ has the form 
\[
h_{I}(t) = \frac{1}{(1-t)^2(1-t^2)} - \frac{s_I(t)}{1-t^2}
\]
where $s_I(t)$ is a Castelnuovo polynomial of even weight $n_e$ and odd weight $n_o$. For $n_e \leq 3$, $n_o \leq 3$ we list the possible Hilbert series for $I$, the corresponding Castelnuovo polynomial, the possible minimal resolutions and $\dim_{k} \Ext^{1}_A(I,I)$.
\def\mystrut{\vrule width 0em height 1em} 
{\small 
\[ 
\begin{array}{|c|l|} 
\hline 
(n_e,n_o) = (0,0) 
& 
h_{I}(t)  = 1 + 2t + 4t^{2} + 6t^{3} + 9t^{4} + 12t^{5} + \ldots 
\mystrut\\ 
& s_{I}(t)  = 0 \\ 
& \dim_{k}  \Ext^{1}_{A}(I,I) = 0 \\ 
& 0  \r A \r I \r 0\\ 
\hline 
(n_e,n_o) = (1,0)& 
h_{I}(t)  = 2t + 3t^{2} + 6t^{3} + 8t^{4} + 12t^{5} + 15t^{6} + \ldots 
\mystrut\\ 
\parbox[c]{0.5cm}{
\unitlength 1mm
\begin{picture}(5.00,0.00)(0,0)

\linethickness{0.15mm}
\put(0.00,-5.00){\line(1,0){5.00}}
\put(0.00,-5.00){\line(0,1){5.00}}
\put(5.00,-5.00){\line(0,1){5.00}}
\put(0.00,0.00){\line(1,0){5.00}}

\linethickness{0.15mm}
\put(0.00,-5.00){\rule{5.00\unitlength}{5.00\unitlength}}

\end{picture}
}
& 
s_{I}(t)  = 1\\ 
&\dim_{k} \Ext^{1}_{A}(I,I) = 0 \\ 
&0  \r A(-2) \r A(-1)^{2} \r I \r 0\\ 
\hline
(n_e,n_o) = (1,1)& 
h_{I}(t)  = t + 3t^{2} + 5t^{3} + 8t^{4} + 11t^{5} + 15t^{6} + \ldots 
\mystrut\\ 
\parbox[c]{1cm}{
\unitlength 1mm
\begin{picture}(10.00,0.00)(0,0)

\linethickness{0.15mm}
\put(0.00,-5.00){\line(1,0){5.00}}
\put(0.00,-5.00){\line(0,1){5.00}}
\put(5.00,-5.00){\line(0,1){5.00}}
\put(0.00,0.00){\line(1,0){5.00}}

\linethickness{0.15mm}
\put(5.00,-5.00){\line(1,0){5.00}}
\put(5.00,-5.00){\line(0,1){5.00}}
\put(10.00,-5.00){\line(0,1){5.00}}
\put(5.00,0.00){\line(1,0){5.00}}

\linethickness{0.15mm}
\put(0.00,-5.00){\rule{5.00\unitlength}{5.00\unitlength}}

\end{picture}
}
&s_{I}(t)  = 1 + t 
 \\ 
&\dim_{k} \Ext^{1}_{A}(I,I) = 2 \\ 
&0  \r A(-3) \r A(-1) \oplus A(-2) \r I \r 0 \\ 
\hline
(n_e,n_o) = (1,2) 
&h_{I}(t)  = 3t^{2} + 4t^{3} + 8t^{4} + 10t^{5} + 15t^{6} + \ldots \mystrut\\ 
\parbox[c]{1cm}{
\unitlength 1mm
\begin{picture}(10.00,5.00)(0,0)

\linethickness{0.15mm}
\put(0.00,-5.00){\line(1,0){5.00}}
\put(0.00,-5.00){\line(0,1){5.00}}
\put(5.00,-5.00){\line(0,1){5.00}}
\put(0.00,0.00){\line(1,0){5.00}}

\linethickness{0.15mm}
\put(5.00,-5.00){\line(1,0){5.00}}
\put(5.00,-5.00){\line(0,1){5.00}}
\put(10.00,-5.00){\line(0,1){5.00}}
\put(5.00,0.00){\line(1,0){5.00}}

\linethickness{0.15mm}
\put(5.00,0.00){\line(1,0){5.00}}
\put(5.00,0.00){\line(0,1){5.00}}
\put(10.00,0.00){\line(0,1){5.00}}
\put(5.00,5.00){\line(1,0){5.00}}

\linethickness{0.15mm}
\put(0.00,-5.00){\rule{5.00\unitlength}{5.00\unitlength}}

\end{picture}
}
&s_{I}(t)  = 1 + 2t 
  \\ 
&\dim_{k}  \Ext^{1}_{A}(I,I) = 0 \\ 
&0  \r A(-3)^{2} \r A(-2)^{3} \r I \r 0 \\
\hline
(n_e,n_o) = (2,1) 
& h_{I}(t)  = t + 2t^{2} + 5t^{3} + 7t^{4} + 11t^{5} + 14t^{6} + \ldots 
\mystrut\\ 
\parbox[c]{1.5cm}{
\unitlength 1mm
\begin{picture}(15.00,0.00)(0,0)

\linethickness{0.15mm}
\put(0.00,-5.00){\line(1,0){5.00}}
\put(0.00,-5.00){\line(0,1){5.00}}
\put(5.00,-5.00){\line(0,1){5.00}}
\put(0.00,0.00){\line(1,0){5.00}}

\linethickness{0.15mm}
\put(5.00,-5.00){\line(1,0){5.00}}
\put(5.00,-5.00){\line(0,1){5.00}}
\put(10.00,-5.00){\line(0,1){5.00}}
\put(5.00,0.00){\line(1,0){5.00}}

\linethickness{0.15mm}
\put(10.00,-5.00){\line(1,0){5.00}}
\put(10.00,-5.00){\line(0,1){5.00}}
\put(15.00,-5.00){\line(0,1){5.00}}
\put(10.00,0.00){\line(1,0){5.00}}

\linethickness{0.15mm}
\put(0.00,-5.00){\rule{5.00\unitlength}{5.00\unitlength}}

\linethickness{0.15mm}
\put(10.00,-5.00){\rule{5.00\unitlength}{5.00\unitlength}}

\end{picture}
}
&s_{I}(t)  = 1 + t + t^{2} 
\quad \\ 
&\dim_{k} \Ext^{1}_{A}(I,I) = 2 \\ 
&0  \r A(-4) \r A(-1) \oplus A(-3) \r I \r 0\\ 
\hline 
\end{array}
\] 
\[
\begin{array}{|c|l|} 
\hline
(n_e,n_o) = (2,2)& 
h_{I}(t)  = 2t^{2} + 4t^{3} + 7t^{4} + 10t^{5} + 14t^{6} + \ldots \mystrut\\ 
\parbox[c]{1.5cm}{
\unitlength 1mm
\begin{picture}(15.00,0.00)(0,0)

\linethickness{0.15mm}
\put(0.00,-10.00){\line(1,0){5.00}}
\put(0.00,-10.00){\line(0,1){5.00}}
\put(5.00,-10.00){\line(0,1){5.00}}
\put(0.00,-5.00){\line(1,0){5.00}}

\linethickness{0.15mm}
\put(5.00,-10.00){\line(1,0){5.00}}
\put(5.00,-10.00){\line(0,1){5.00}}
\put(10.00,-10.00){\line(0,1){5.00}}
\put(5.00,-5.00){\line(1,0){5.00}}

\linethickness{0.15mm}
\put(10.00,-10.00){\line(1,0){5.00}}
\put(10.00,-10.00){\line(0,1){5.00}}
\put(15.00,-10.00){\line(0,1){5.00}}
\put(10.00,-5.00){\line(1,0){5.00}}

\linethickness{0.15mm}
\put(5.00,-5.00){\line(1,0){5.00}}
\put(5.00,-5.00){\line(0,1){5.00}}
\put(10.00,-5.00){\line(0,1){5.00}}
\put(5.00,0.00){\line(1,0){5.00}}

\linethickness{0.15mm}
\put(0.00,-10.00){\rule{5.00\unitlength}{5.00\unitlength}}

\linethickness{0.15mm}
\put(10.00,-10.00){\rule{5.00\unitlength}{5.00\unitlength}}

\end{picture}
}
&s_{I}(t)  = 1 + 2t + t^{2} 
\\ 
&\dim_{k} \Ext^{1}_{A}(I,I) = 4 \\ 
&0  \r A(-4) \r A(-2)^{2} \r I \r 0\\ 
&0  \r A(-3)\oplus A(-4) \r A(-2)^{2} \oplus A(-3)\r I \r 0\\ 
\cline{2-2} 
&h_{I}(t)  = t + 2t^{2} + 4t^{3} + 7t^{4} + 10t^{5} + 14t^{6} + \ldots 
\mystrut\\ 
\parbox[c]{2cm}{
\unitlength 1mm
\begin{picture}(20.00,0.00)(0,0)

\linethickness{0.15mm}
\put(0.00,-5.00){\line(1,0){5.00}}
\put(0.00,-5.00){\line(0,1){5.00}}
\put(5.00,-5.00){\line(0,1){5.00}}
\put(0.00,0.00){\line(1,0){5.00}}

\linethickness{0.15mm}
\put(5.00,-5.00){\line(1,0){5.00}}
\put(5.00,-5.00){\line(0,1){5.00}}
\put(10.00,-5.00){\line(0,1){5.00}}
\put(5.00,0.00){\line(1,0){5.00}}

\linethickness{0.15mm}
\put(10.00,-5.00){\line(1,0){5.00}}
\put(10.00,-5.00){\line(0,1){5.00}}
\put(15.00,-5.00){\line(0,1){5.00}}
\put(10.00,0.00){\line(1,0){5.00}}

\linethickness{0.15mm}
\put(15.00,-5.00){\line(1,0){5.00}}
\put(15.00,-5.00){\line(0,1){5.00}}
\put(20.00,-5.00){\line(0,1){5.00}}
\put(15.00,0.00){\line(1,0){5.00}}

\linethickness{0.15mm}
\put(0.00,-5.00){\rule{5.00\unitlength}{5.00\unitlength}}

\linethickness{0.15mm}
\put(10.00,-5.00){\rule{5.00\unitlength}{5.00\unitlength}}

\end{picture}
}
& s_{I}(t)  = 1 + t + t^{2} + t^{3}\\ 
&\dim_{k} \Ext^{1}_{A}(I,I) = 3 \\ 
&0  \r A(-5) \r A(-1) \oplus A(-4) \r I \r 0\\ 
\hline 
(n_e,n_o) = (2,3) 
&h_{I}(t)  = 2t^{2} + 3t^{3} + 7t^{4} + 9t^{5} + 14t^{6} + \ldots \mystrut\\ 
\parbox[c]{1.5cm}{
\unitlength 1mm
\begin{picture}(20.00,5.00)(0,0)

\linethickness{0.15mm}
\put(0.00,-5.00){\line(1,0){5.00}}
\put(0.00,-5.00){\line(0,1){5.00}}
\put(5.00,-5.00){\line(0,1){5.00}}
\put(0.00,0.00){\line(1,0){5.00}}

\linethickness{0.15mm}
\put(5.00,-5.00){\line(1,0){5.00}}
\put(5.00,-5.00){\line(0,1){5.00}}
\put(10.00,-5.00){\line(0,1){5.00}}
\put(5.00,0.00){\line(1,0){5.00}}

\linethickness{0.15mm}
\put(10.00,-5.00){\line(1,0){5.00}}
\put(10.00,-5.00){\line(0,1){5.00}}
\put(15.00,-5.00){\line(0,1){5.00}}
\put(10.00,0.00){\line(1,0){5.00}}

\linethickness{0.15mm}
\put(5.00,0.00){\line(1,0){5.00}}
\put(5.00,0.00){\line(0,1){5.00}}
\put(10.00,0.00){\line(0,1){5.00}}
\put(5.00,5.00){\line(1,0){5.00}}

\linethickness{0.15mm}
\put(15.00,-5.00){\line(1,0){5.00}}
\put(15.00,-5.00){\line(0,1){5.00}}
\put(20.00,-5.00){\line(0,1){5.00}}
\put(15.00,0.00){\line(1,0){5.00}}

\linethickness{0.15mm}
\put(0.00,-5.00){\rule{5.00\unitlength}{5.00\unitlength}}

\linethickness{0.15mm}
\put(10.00,-5.00){\rule{5.00\unitlength}{5.00\unitlength}}

\end{picture}
}
&s_{I}(t)  = 1 + 2t + t^{2} + t^3 
\\ 
&\dim_{k} \Ext_A^{1}(I,I) = 2 \\ 
&0  \r A(-3) \oplus A(-5) \r A(-2)^{2} \oplus A(-4) \r I \r 0\\
\hline 
(n_e,n_o) = (3,2) 
&h_{I}(t)  = t^{2} + 4t^{3} + 6t^{4} + 10t^{5} + 13t^{6} + \ldots \mystrut\\ 
\parbox[c]{1.5cm}{
\unitlength 1mm
\begin{picture}(15.00,5.00)(0,0)

\linethickness{0.15mm}
\put(0.00,-5.00){\line(1,0){5.00}}
\put(0.00,-5.00){\line(0,1){5.00}}
\put(5.00,-5.00){\line(0,1){5.00}}
\put(0.00,0.00){\line(1,0){5.00}}

\linethickness{0.15mm}
\put(5.00,-5.00){\line(1,0){5.00}}
\put(5.00,-5.00){\line(0,1){5.00}}
\put(10.00,-5.00){\line(0,1){5.00}}
\put(5.00,0.00){\line(1,0){5.00}}

\linethickness{0.15mm}
\put(10.00,-5.00){\line(1,0){5.00}}
\put(10.00,-5.00){\line(0,1){5.00}}
\put(15.00,-5.00){\line(0,1){5.00}}
\put(10.00,0.00){\line(1,0){5.00}}

\linethickness{0.15mm}
\put(5.00,0.00){\line(1,0){5.00}}
\put(5.00,0.00){\line(0,1){5.00}}
\put(10.00,0.00){\line(0,1){5.00}}
\put(5.00,5.00){\line(1,0){5.00}}

\linethickness{0.15mm}
\put(10.00,0.00){\line(1,0){5.00}}
\put(10.00,0.00){\line(0,1){5.00}}
\put(15.00,0.00){\line(0,1){5.00}}
\put(10.00,5.00){\line(1,0){5.00}}

\linethickness{0.15mm}
\put(0.00,-5.00){\rule{5.00\unitlength}{5.00\unitlength}}

\linethickness{0.15mm}
\put(10.00,0.00){\rule{5.00\unitlength}{5.00\unitlength}}

\linethickness{0.15mm}
\put(10.00,-5.00){\rule{5.00\unitlength}{5.00\unitlength}}

\end{picture}
}
&s_{I}(t)  = 1 + 2t + 2t^{2} 
\\ 
&\dim_{k} \Ext_A^{1}(I,I) = 4 \\ 
&0  \r A(-4)^{2} \r A(-2) \oplus A(-3)^{2} \r I \r 0\\ 
\cline{2-2} 
&h_{I}(t)  = t + 2t^{2} + 4t^{3} + 6t^{4} + 10t^{5} + 13t^{6} + 
\ldots 
\mystrut\\ 
\parbox[c]{2.5cm}{
\unitlength 1mm
\begin{picture}(25.00,0.00)(0,0)

\linethickness{0.15mm}
\put(0.00,-5.00){\line(1,0){5.00}}
\put(0.00,-5.00){\line(0,1){5.00}}
\put(5.00,-5.00){\line(0,1){5.00}}
\put(0.00,0.00){\line(1,0){5.00}}

\linethickness{0.15mm}
\put(5.00,-5.00){\line(1,0){5.00}}
\put(5.00,-5.00){\line(0,1){5.00}}
\put(10.00,-5.00){\line(0,1){5.00}}
\put(5.00,0.00){\line(1,0){5.00}}

\linethickness{0.15mm}
\put(10.00,-5.00){\line(1,0){5.00}}
\put(10.00,-5.00){\line(0,1){5.00}}
\put(15.00,-5.00){\line(0,1){5.00}}
\put(10.00,0.00){\line(1,0){5.00}}

\linethickness{0.15mm}
\put(15.00,-5.00){\line(1,0){5.00}}
\put(15.00,-5.00){\line(0,1){5.00}}
\put(20.00,-5.00){\line(0,1){5.00}}
\put(15.00,0.00){\line(1,0){5.00}}

\linethickness{0.15mm}
\put(20.00,-5.00){\line(1,0){5.00}}
\put(20.00,-5.00){\line(0,1){5.00}}
\put(25.00,-5.00){\line(0,1){5.00}}
\put(20.00,0.00){\line(1,0){5.00}}

\linethickness{0.15mm}
\put(0.00,-5.00){\rule{5.00\unitlength}{5.00\unitlength}}

\linethickness{0.15mm}
\put(10.00,-5.00){\rule{5.00\unitlength}{5.00\unitlength}}

\linethickness{0.15mm}
\put(20.00,-5.00){\rule{5.00\unitlength}{5.00\unitlength}}

\end{picture}
}
&s_{I}(t)  = 1 + t + t^{2} + t^{3} + t^{4} \\ 
&\dim_{k} \Ext_A^{1}(I,I) = 3 \\ 
&0  \r A(-6) \r A(-1) \oplus A(-5) \r I \r 0\\ 
\hline 
(n_e,n_o) = (3,3)  
&h_{I}(t)  = t^{2} + 3t^{3} + 6t^{4} + 9t^{5} + 13t^{6} + \ldots\mystrut \\ 
\parbox[c]{2cm}{
\unitlength 1mm
\begin{picture}(20.00,0.00)(0,0)

\linethickness{0.15mm}
\put(0.00,-10.00){\line(1,0){5.00}}
\put(0.00,-10.00){\line(0,1){5.00}}
\put(5.00,-10.00){\line(0,1){5.00}}
\put(0.00,-5.00){\line(1,0){5.00}}

\linethickness{0.15mm}
\put(5.00,-10.00){\line(1,0){5.00}}
\put(5.00,-10.00){\line(0,1){5.00}}
\put(10.00,-10.00){\line(0,1){5.00}}
\put(5.00,-5.00){\line(1,0){5.00}}

\linethickness{0.15mm}
\put(10.00,-10.00){\line(1,0){5.00}}
\put(10.00,-10.00){\line(0,1){5.00}}
\put(15.00,-10.00){\line(0,1){5.00}}
\put(10.00,-5.00){\line(1,0){5.00}}

\linethickness{0.15mm}
\put(5.00,-5.00){\line(1,0){5.00}}
\put(5.00,-5.00){\line(0,1){5.00}}
\put(10.00,-5.00){\line(0,1){5.00}}
\put(5.00,0.00){\line(1,0){5.00}}

\linethickness{0.15mm}
\put(10.00,-5.00){\line(1,0){5.00}}
\put(10.00,-5.00){\line(0,1){5.00}}
\put(15.00,-5.00){\line(0,1){5.00}}
\put(10.00,0.00){\line(1,0){5.00}}

\linethickness{0.15mm}
\put(15.00,-10.00){\line(1,0){5.00}}
\put(15.00,-10.00){\line(0,1){5.00}}
\put(20.00,-10.00){\line(0,1){5.00}}
\put(15.00,-5.00){\line(1,0){5.00}}

\linethickness{0.15mm}
\put(0.00,-10.00){\rule{5.00\unitlength}{5.00\unitlength}}

\linethickness{0.15mm}
\put(10.00,-10.00){\rule{5.00\unitlength}{10.00\unitlength}}

\end{picture}
}
& s_{I}(t)  = 1 + 2t + 2t^{2} + t^{3}\\ 
&\dim_{k} \Ext_A^{1}(I,I) = 6 \\ 
&0  \r A(-5) \r A(-2) \oplus A(-3) \r I \r 0\\ 
&0  \r A(-4)\oplus A(-5) \r A(-2) \oplus A(-3)\oplus A(-4) \r I \r 0\\ 
\cline{2-2} 
&h_{I}(t)  = 2t^{2} + 3t^{3} + 6t^{4} + 9t^{5} + 13t^{6} + \ldots \mystrut\\ 
\parbox[c]{2.5cm}{
\unitlength 1mm
\begin{picture}(25.00,5.00)(0,0)

\linethickness{0.15mm}
\put(0.00,-5.00){\line(1,0){5.00}}
\put(0.00,-5.00){\line(0,1){5.00}}
\put(5.00,-5.00){\line(0,1){5.00}}
\put(0.00,0.00){\line(1,0){5.00}}

\linethickness{0.15mm}
\put(5.00,-5.00){\line(1,0){5.00}}
\put(5.00,-5.00){\line(0,1){5.00}}
\put(10.00,-5.00){\line(0,1){5.00}}
\put(5.00,0.00){\line(1,0){5.00}}

\linethickness{0.15mm}
\put(10.00,-5.00){\line(1,0){5.00}}
\put(10.00,-5.00){\line(0,1){5.00}}
\put(15.00,-5.00){\line(0,1){5.00}}
\put(10.00,0.00){\line(1,0){5.00}}

\linethickness{0.15mm}
\put(5.00,0.00){\line(1,0){5.00}}
\put(5.00,0.00){\line(0,1){5.00}}
\put(10.00,0.00){\line(0,1){5.00}}
\put(5.00,5.00){\line(1,0){5.00}}

\linethickness{0.15mm}
\put(15.00,-5.00){\line(1,0){5.00}}
\put(15.00,-5.00){\line(0,1){5.00}}
\put(20.00,-5.00){\line(0,1){5.00}}
\put(15.00,0.00){\line(1,0){5.00}}

\linethickness{0.15mm}
\put(20.00,-5.00){\line(1,0){5.00}}
\put(20.00,-5.00){\line(0,1){5.00}}
\put(25.00,-5.00){\line(0,1){5.00}}
\put(20.00,0.00){\line(1,0){5.00}}

\linethickness{0.15mm}
\put(0.00,-5.00){\rule{5.00\unitlength}{5.00\unitlength}}

\linethickness{0.15mm}
\put(10.00,-5.00){\rule{5.00\unitlength}{5.00\unitlength}}

\linethickness{0.15mm}
\put(20.00,-5.00){\rule{5.00\unitlength}{5.00\unitlength}}

\end{picture}
}
&s_{I}(t)  = 1 + 2t + t^{2} + t^{3} + t^{4}\\ 
&\dim_{k} \Ext_A^{1}(I,I) = 4 \\ 
&0  \r A(-3) \oplus A(-6) \r A(-2)^{2} \oplus A(-5) \r I \r 0\\ 
\cline{2-2} 
&h_{I}(t)  = t + 2t^{2} + 4t^{3} + 6t^{4} + 9t^{5} + 13t^{6} +  \ldots 
\mystrut\\ 
\parbox[c]{3cm}{
\unitlength 1mm
\begin{picture}(30.00,0.00)(0,0)

\linethickness{0.15mm}
\put(0.00,-5.00){\line(1,0){5.00}}
\put(0.00,-5.00){\line(0,1){5.00}}
\put(5.00,-5.00){\line(0,1){5.00}}
\put(0.00,0.00){\line(1,0){5.00}}

\linethickness{0.15mm}
\put(5.00,-5.00){\line(1,0){5.00}}
\put(5.00,-5.00){\line(0,1){5.00}}
\put(10.00,-5.00){\line(0,1){5.00}}
\put(5.00,0.00){\line(1,0){5.00}}

\linethickness{0.15mm}
\put(10.00,-5.00){\line(1,0){5.00}}
\put(10.00,-5.00){\line(0,1){5.00}}
\put(15.00,-5.00){\line(0,1){5.00}}
\put(10.00,0.00){\line(1,0){5.00}}

\linethickness{0.15mm}
\put(15.00,-5.00){\line(1,0){5.00}}
\put(15.00,-5.00){\line(0,1){5.00}}
\put(20.00,-5.00){\line(0,1){5.00}}
\put(15.00,0.00){\line(1,0){5.00}}

\linethickness{0.15mm}
\put(20.00,-5.00){\line(1,0){5.00}}
\put(20.00,-5.00){\line(0,1){5.00}}
\put(25.00,-5.00){\line(0,1){5.00}}
\put(20.00,0.00){\line(1,0){5.00}}

\linethickness{0.15mm}
\put(25.00,-5.00){\line(1,0){5.00}}
\put(25.00,-5.00){\line(0,1){5.00}}
\put(30.00,-5.00){\line(0,1){5.00}}
\put(25.00,0.00){\line(1,0){5.00}}

\linethickness{0.15mm}
\put(0.00,-5.00){\rule{5.00\unitlength}{5.00\unitlength}}

\linethickness{0.15mm}
\put(10.00,-5.00){\rule{5.00\unitlength}{5.00\unitlength}}

\linethickness{0.15mm}
\put(20.00,0.00){\rule{0.00\unitlength}{0.00\unitlength}}

\linethickness{0.15mm}
\put(20.00,-5.00){\rule{5.00\unitlength}{5.00\unitlength}}

\end{picture}
}
&s_{I}(t)  = 1 + t + t^{2} + t^{3} + t^{4} + t^{5} 
 \\ 
&\dim_{k} \Ext_A^{1}(I,I) = 4 \\ 
&0  \r A(-7) \r A(-1) \oplus A(-6) \r I \r 0\\ 
\hline 
\end{array} 
\]
} 
 
\end{document}